\documentclass[11pt,a4paper]{article}

\usepackage[english]{babel}
\usepackage[utf8]{inputenc}
\usepackage{algorithm}
\usepackage{algpseudocode}
\usepackage{amsmath,amssymb,amsbsy,amsfonts}
\usepackage{bm}
\usepackage{stmaryrd}
\usepackage[scr=rsfs]{mathalpha}
\usepackage{amsthm,mathtools}
\numberwithin{equation}{section}
\usepackage[paperwidth=192mm,
		    paperheight=262mm,
		    vmargin={19mm,19mm},
		    hmargin={13.7mm,13.7mm},
		    headsep=12pt,
		    footskip=12pt]{geometry}

\usepackage{hyperref}
\usepackage{orcidlink}

\usepackage{siunitx}
\usepackage{enumitem}

\usepackage{graphicx}
\usepackage{float}
\usepackage{caption}
\usepackage{subcaption}
\usepackage{tabularx}
\usepackage{booktabs}

\usepackage{tikz}
\usetikzlibrary{trees}
\usepackage{pgfplots}
\usetikzlibrary{external}
\usetikzlibrary{positioning}

\newtheorem{theorem}{Theorem}[section]
\newtheorem{definition}[theorem]{Definition}
\newtheorem{remark}{Remark}[section]

\algnewcommand\algorithmicinput{\textbf{Input:}}
\algnewcommand\algorithmicoutput{\textbf{Output:}}
\algnewcommand\Given{\item[\algorithmicinput]}

\DeclareMathOperator{\grad}{\nabla}
\DeclareMathOperator{\dive}{\nabla\cdot}

\newcommand{\ave}[1]{\left\{\!\!\left\{#1\right\}\!\!\right\}}
\newcommand{\jump}[1]{\left\llbracket#1\right\rrbracket}

\newcommand{\pad}[2]{\frac{\partial{#1}}{\partial{#2}}}
\newcommand{\odif}[2]{\frac{\mathrm{d}{#1}}{\mathrm{d}{#2}}}

\newcommand{\rpth}[1]{\left(#1\right)}
\newcommand{\spth}[1]{\left[#1\right]}
\newcommand{\cpth}[1]{\left\{#1\right\}}

\newcommand{\partition}{\mathscr{T}_{h}}
\newcommand{\facesinternal}{\mathscr{F}^\mathrm{I}_{h}}
\newcommand{\faces}{\mathscr{F}_{h}}
\newcommand{\facesN}{\mathscr{F}_{h}^{\mathrm{N}}}
\newcommand{\facesD}{\mathscr{F}_{h}^{\mathrm{D}}}
\newcommand{\facesboundary}{\mathscr{F}^{\mathrm{B}_{h}}}
\newcommand{\Vh}{\mathbf{V}_{h}^{\mathrm{DG}}}

\newcommand{\timesuperscript}{n}
\newcommand{\ncomp}{N}
\newcommand{\dimDG}{M}

\newcommand{\Verhulstcomp}{u}
\newcommand{\FKcomp}{u}

\newcommand{\GSfirstcomp}[1][]{u_{1#1}}
\newcommand{\GSsecondcomp}[1][]{u_{2#1}}

\newcommand{\orderthreeLDp}{IMEX$(3,5)$-$\text{LDp}$}
\newcommand{\orderthreeLDs}{IMEX$(3,5)$-$\text{LDs}$}
\newcommand{\orderthreeLDsone}{IMEX$(3,5)$-$\text{LDs}_{1}$}
\newcommand{\orderthreeLDstwo}{IMEX$(3,5)$-$\text{LDs}_{2}$}

\newcommand{\orderfourLDp}{IMEX$(4,6)$-$\text{LDp}$}
\newcommand{\orderfourLDpone}{IMEX$(4,6)$-$\text{LDp}_{1}$}
\newcommand{\orderfourLDptwo}{IMEX$(4,6)$-$\text{LDp}_{2}$}
\newcommand{\orderfourLDpthree}{IMEX$(4,6)$-$\text{LDp}_{3}$}

\newcommand{\boscarinoduemilasedici}{IMEX-SSP$(4,3,3)$}
\newcommand{\boscarinoduemilaventidue}{IMEX-RK$(4,3,3)$}
\newcommand{\boscarinoduemilanove}{BHR$(5,5,3)$}
\newcommand{\carpenterkennedyduemilatre}{ARK3(2)4L[2]SA}
\newcommand{\nostroordinetreprimo}{\orderthreeLDsone}
\newcommand{\nostroordinetresecondo}{\orderthreeLDp}
\newcommand{\nostroordinetreterzo}{\orderthreeLDstwo}
\newcommand{\calvoduemilauno}{ARS$(5,5,4)$}
\newcommand{\kennedycarpenterduemilatre}{ARK4(3)6L[2]SA}\newcommand{\kennedycarpenterduemiladiciannove}{ARK4(3)7L[2]SA$_{1}$}
\newcommand{\nostroordinequattroprimo}{\orderfourLDpone}
\newcommand{\nostroordinequattrosecondo}{\orderfourLDptwo}
\newcommand{\nostroordinequattroterzo}{\orderfourLDpthree}

\begin{document}

\pagenumbering{arabic}

\title{Optimized high-order IMEX-RK schemes for degenerate diffusion-reaction problems with application to travelling waves phenomena\footnote{\textbf{Fundings:} This work is partially funded by the European Union (ERC SyG, NEMESIS, project number 101115663). Views and opinions expressed are, however, those of the authors only and do not necessarily reflect those of the European Union or the European Research Council Executive Agency. Neither the European Union nor the granting authority can be held responsible for them. The present research is part of the activities of the Dipartimento di Eccellenza 2023-2027 grant, funded by MUR. MC acknowledges “INdAM - GNCS Project”, codice CUP E53C25002010001. PFA, MC, and GO are members of INdAM-GNCS.}}

\date{}

\author{Paola F. Antonietti\,\orcidlink{0000-0002-2138-3878}\thanks{MOX-Dipartimento di Matematica, Politecnico di Milano, Piazza Leonardo da Vinci 32, Milan, 20133, Italy (\href{mailto:paola.antonietti@polimi.it}{paola.antonietti@polimi.it}, \href{mailto:mattia.corti@polimi.it}{mattia.corti@polimi.it})} 
\and
Mattia Corti\,\orcidlink{0000-0002-7014-972X}\footnotemark[2]
\and
Giuseppe Orlando\,\orcidlink{0000-0002-7119-4231}\thanks{CMAP, CNRS, \'{E}cole polytechnique, Institut Polytechnique de Paris, Route de Saclay, 91120 Palaiseau, France (\href{mailto:giuseppe.orlando@polytechnique.edu}{giuseppe.orlando@polytechnique.edu})}
}       

\maketitle


%%%%%%%%%%%%%%%%%%%%%%%%%%%%%%%%%%%%%%
%%%%%%%%%%%% Abstract %%%%%%%%%%%%%%%%
%%%%%%%%%%%%%%%%%%%%%%%%%%%%%%%%%%%%%%
\begin{abstract}
We study a class of IMplicit-EXplicit Runge--Kutta (IMEX-RK) schemes for the numerical approximation of reaction and diffusion-reaction problems arising in a variety of biological and physical applications. Such models may admit travelling wave solutions, with the Fisher--Kolmogorov equation representing a prototypical example. Motivated by this feature, the proposed time integration schemes are designed to accurately capture sharp propagating fronts. We also investigate a less standard use of IMEX-RK methods that circumvents a splitting of reaction terms into linear and nonlinear components, while still requiring the solution of linear systems at each stage. This semi-implicit formulation, referred to as SI-IMEX-RK, enables a targeted treatment of stiffness by isolating its relevant contributions. The time discretization is coupled with a high-order polygonal discontinuous Galerkin method for space discretization, resulting in a flexible and robust framework for the treatment of multiscale dynamics in complex geometries. A comprehensive validation strategy is presented to assess the accuracy and stability properties of the proposed schemes across a hierarchy of increasingly challenging test problems.
\end{abstract}
\setlength{\leftskip}{0pt}
\setlength{\rightskip}{0pt}

%%%%%%%%%%%%%%%%%%%%%%%%%%%%%%%%%%%%%%
%%%%%%%%% Introduction %%%%%%%%%%%%%%%
%%%%%%%%%%%%%%%%%%%%%%%%%%%%%%%%%%%%%%
\section{Introduction}
\label{sec:intro} \indent

Diffusion-reaction systems constitute a standard mathematical framework for modelling a wide range of phenomena arising, for instance, in ecology, combustion, and biology~\cite{deluca:2026, duarte:2012, dumont:2013, fisher:1937, peters:2009, weickenmeier:2019}. A common numerical approach for approximating such systems is based on the method of lines. Within this framework, the spatial discretization of the governing equations is performed first, typically leading to a system of ordinary differential equations (ODEs) of the form
\begin{equation*}
    \odif{\bm{y}}{t} = \mathbf{F}_{\mathrm{D}}(\bm{y}) + \mathbf{F}_{\mathrm{R}}(\bm{y},t),
\end{equation*}
where $\mathbf{F}_{\mathrm{D}}$ and $\mathbf{F}_{\mathrm{R}}$ denote the discrete diffusion and reaction operators, respectively. Several approaches have been proposed in the literature for the numerical solution of this problem, including exponential integrators~\cite{caliari:2007}, fully implicit schemes~\cite{sheu:2000}, operator splitting techniques~\cite{duarte:2012, dumont:2013}, and IMplicit-EXplicit (IMEX) Runge--Kutta (RK) methods~\cite{abdulle:2013, deluca:2026}. In this work, we focus on the latter class of methods. IMEX-RK methods have proven to be particularly effective for a wide range of multiscale problems, see, e.g.,~\cite{boscarino:2024, boscarino:2009_IMEX_relaxation, boscheri:2023, boscheri:2024, dumbser:2016, orlando:2022, pareschi:2005}. Our interest is specifically directed toward regimes in which the reaction term becomes dominant and therefore potentially stiff, making a fully explicit treatment prohibitively restrictive. In the singular limit, the differential system degenerates into an algebraic constraint, and IMEX-RK methods must satisfy additional order conditions \cite{boscarino:2007,boscarino:2009_third_order} to preserve the designed order of accuracy for the algebraic variable. At the same time, reaction operators are often (highly) nonlinear, so that a fully implicit discretization would require the solution of nonlinear systems at each time step or stage. Despite the absence of spatial coupling that allows decomposing the reaction term into independent local systems, such an approach may still be computationally inefficient. As far as the diffusion operator $\mathbf{F}_{\mathrm{D}}$ is concerned, we employ an implicit discretization. Nevertheless, we mention that stabilized explicit RK schemes based on Orthogonal Runge--Kutta--Chebyshev (ROCK) methods~\cite{abdulle:2002, abdulle:2001} may represent an attractive alternative for the explicit treatment of diffusion-dominated problems. However, despite being formally explicit RK methods, ROCK schemes typically involve a very large number of stages, although they are designed to be low-storage and have a very low cost per stage. Moreover, ROCK methods are not naturally expressed in terms of classical Butcher tableaux, which makes their integration into standard RK frameworks less straightforward.
\\~\\
A commonly adopted strategy for problems involving nonlinear reaction terms consists in treating the linear part implicitly while discretizing the nonlinear contribution explicitly~\cite{calvo:2001, corti:2024_IMEX}. This approach represents the first class of methods considered in this work. However, in order to avoid a rigid splitting between linear and nonlinear components of the reaction operator, we further develop a class of semi-implicit time discretizations based on IMEX Runge--Kutta schemes. Following the partitioned formulation proposed in~\cite{boscarino:2016}, these methods --- hereafter referred to as Semi-Implicit IMEX Runge--Kutta (SI-IMEX-RK) methods --- allow one to identify and isolate the (possibly linear) source of stiffness within the system. As will be shown, this framework significantly enlarges the applicability of IMEX-RK methods and leads to efficient and flexible schemes that can be naturally extended to more complex nonlinear dependencies. Finally, we recall that diffusion-reaction systems may admit travelling wave solutions~\cite{antonietti:2026, corti:2024_IMEX, francois:2020, leimer_saglio:2026}. Accurate numerical schemes for such problems should therefore provide a good approximation of the propagation speed and shape of these waves, while minimising artificial numerical dissipation and dispersion. This work tries to construct high–order IMEX-RK schemes that provide an improved numerical description of wave propagation properties. In this context, the main novelty of the proposed work lies in the development of schemes aimed at improving the numerical approximation of wave propagation and the algebraic accuracy. In addition, the SI-IMEX-RK framework is exploited to provide a more flexible treatment of stiff reaction mechanisms within the overall formulation.
\\~\\
The remainder of the paper is organized as follows. In Section~\ref{sec:model}, we introduce the mathematical model of interest together with the abstract system of ODEs considered throughout the paper. Section~\ref{sec:IMEX} is devoted to the time discretization strategy based on IMEX-RK methods. After recalling the properties of interest for the model under consideration, we first discuss a classical IMEX-RK approach based on an implicit treatment of the linear reaction term and an explicit discretization of the nonlinear contribution. We then introduce the novel SI-IMEX-RK methodology inspired by~\cite{boscarino:2016}. In Section~\ref{sec:optmized_IMEX}, we derive new high-order IMEX-RK schemes optimized with respect to the properties discussed in Section~\ref{sec:IMEX}. In Section~\ref{sec:polydg}, we describe the chosen space discretization, which is based on a polytopal discontinuous Galerkin (polydG) method~\cite{cangiani:2017, cangiani:2014}. Sections~\ref{sec:numerical_results_ODE} and~\ref{sec:numerical_results_PDE} present a collection of numerical experiments that validate the proposed methods on various ODEs and PDEs, respectively. Finally, conclusions and perspectives for future work are outlined in Section~\ref{sec:conclu}.

%%%%%%%%%%%%%%%%%%%%%%%%%%%%%%%%%%%%%
%%%%%%%%%%%%%%%%% Model %%%%%%%%%%%%%
%%%%%%%%%%%%%%%%%%%%%%%%%%%%%%%%%%%%%
\section{Mathematical model}
\label{sec:model}  \indent

In this section, we briefly outline the mathematical model. Let $\Omega \subset \mathbb{R}^{d}, d=2,3$, be a connected open bounded set with a sufficiently smooth boundary $\partial\Omega$ and denote by $\bm{x}$ the spatial coordinates and by $t$ the temporal coordinate. The boundary can be divided into two subsets $\Gamma_{\mathrm{D}}$ and $\Gamma_{\mathrm{N}}$, where we impose Dirichlet and Neumann boundary conditions, respectively, and such that $\Gamma_{\mathrm{D}} \cup \Gamma_{\mathrm{N}} = \partial\Omega$ and $\Gamma_{\mathrm{D}} \cap \Gamma_{\mathrm{N}} = \emptyset$. We consider as model problem a possibly degenerate diffusion-reaction system with nonlinear reaction term: find $\bm{u} = \bm{u}(\bm{x},t) \in \mathbb{R}^{\ncomp}$, with $\ncomp$ denoting the number of components, such that
\begin{subequations}\label{eq:model}
\begin{alignat}{3}
    \pad{\bm{u}}{t} &= \dive\rpth{\mathbf{\Sigma}\grad\bm{u}} + \mathbf{F}(\bm{u},t), && \quad \mathrm{in}\,\Omega \times (0,T], \label{eq:model_equation} \\[3pt]
    \bm{u} &= \bm{u}_{\mathrm{D}}, && \quad \mathrm{on}\, \Gamma_{\mathrm{D}} \times (0,T], \\[3pt]
    \mathbf{\Sigma}\grad\bm{u} \cdot \bm{n}_{\Omega} &= \bm{g}_{\mathrm{N}}, && \quad \mathrm{on}\,\Gamma_{\mathrm{N}} \times (0,T], \\[3pt]
    \bm{u}(\cdot,0) &= \bm{u}_{0}, && \quad \mathrm{in}\,\Omega,
\end{alignat}
\end{subequations}
where $\mathbf{\Sigma} = \mathrm{diag}(\sigma_{1},\dots,\sigma_{\ncomp})$ is the diagonal matrix containing the diffusion coefficients ($\sigma_{j} \geq 0$ for $j=1,\dots,\ncomp$) and $\mathbf{F}:\mathbb{R}^{\ncomp} \rightarrow \mathbb{R}^{\ncomp}$ is the nonlinear reaction term. Moreover, $\bm{u}_{\mathrm{D}} \in H^{1/2}(\Gamma_{\mathrm{D}})$ and $\bm{g}_\mathrm{N}\in H^{-1/2}(\Gamma_{\mathrm{N}})$ are the Dirichlet and Neumann boundary conditions, respectively, while $\bm{n}_{\Omega}$ is the outward unit normal from the domain $\Omega$. Finally, $\bm{u}_{0} \in L^{2}(\Omega)$ is the initial condition. We denote by $u_{j}$ the $j$-th component of the vectorial variable $\bm{u}$. As mentioned in the Introduction, a notable example of this class of problems is the Fisher--Kolmogorov (FK) equation~\cite{fisher:1937}, which will be introduced in Section~\ref{ssec:travelling_wave_FK} and will serve as relevant numerical testbed for validating our approach.
\\~\\
In the limit of dominant reaction, Equation~\eqref{eq:model} could be interpreted as a differential-algebraic equation (DAE), see, e.g.,~\cite{boscarino:2007, boscarino:2009_IMEX_relaxation}. Within the notation of~\eqref{eq:model}, this corresponds to the presence of rapidly relaxing reaction dynamics that drive the solution toward a low-dimensional invariant manifold defined by $\mathbf{F}(\bm{u},t) = \bm{0}$. More precisely, we assume that the system exhibits an \emph{index-1} DAE structure, in the sense that the algebraic constraint $\mathbf{F}(\bm{u},t) = \bm{0}$ defines locally a smooth manifold. In particular, we assume that the Jacobian matrix $\partial_{\bm{u}}\mathbf{F}(\bm{u},t)$ is nonsingular when evaluated at the corresponding equilibrium solution, which guarantees that the constraint can be locally resolved with respect to $\bm{u}$ by the implicit function theorem. This intrinsic stiffness and multiscale nature poses severe challenges for standard time integration methods, often leading to order reduction and loss of accuracy unless specifically designed strategies are employed. In particular, the coexistence of diffusive and (highly) stiff reactive dynamics calls for time discretizations that combine robustness, accuracy, and stability in a non-standard and non-trivial way, motivating the use of the time discretization strategies that will be depicted in Section~\ref{sec:IMEX}. An additional structural feature of systems of the form~\eqref{eq:model} is that they may admit travelling wave solutions~\cite{antonietti:2026, corti:2024_IMEX, corti:2024_SP, francois:2020}, which play a fundamental role in the qualitative behaviour of the underlying dynamics. This observation motivates the development and investigation of novel time discretization schemes designed to improve the description of key wave propagation properties, as discussed in Sections~\ref{sec:IMEX} and~\ref{sec:optmized_IMEX}.

%%%%%%%%%%%%%%%%%%%%%%%%%%%%%%%%%%%%%%
%%%%%%%%%%% IMEX-RK schemes %%%%%%%%%%
%%%%%%%%%%%%%%%%%%%%%%%%%%%%%%%%%%%%%%
\section{Time discretization: IMEX-RK methods}
\label{sec:IMEX}

In this work, we first focus on the time discretization of the governing problem~\eqref{eq:model}. Regardless of the particular space discretization employed, the resulting semi-discrete problem can be written in the abstract form
\begin{subequations}\label{eq:algebraic}
\begin{alignat}{3}
    \mathbf{M}\dot{\boldsymbol{U}}(t) + \mathbf{K}\boldsymbol{U}(t) &= \boldsymbol{G}\rpth{\boldsymbol{U}(t),t},
    && \qquad t \in (0,T], \label{eq:algebraic_ode} \\
    \boldsymbol{U}(0) &= \boldsymbol{U}_{0}.
\end{alignat}
\end{subequations}
The specific form of the matrices, vectors, and discrete operators appearing in~\eqref{eq:algebraic} depends on the chosen space discretization, and their expressions and interpretation will be detailed in Section~\ref{sec:polydg}. Since the time-integration procedure can be described independently of these space discretization details, we now discuss the temporal discretization of~\eqref{eq:algebraic} by means of IMEX-RK methods~\cite{boscarino:2024, kennedy:2003}. In particular, we recall some concepts related to IMEX-RK methods that will be particularly relevant in Section~\ref{sec:optmized_IMEX}. IMEX-RK methods are represented compactly by the following two Butcher tableaux~\cite{butcher:2008}
\begin{center}
    \begin{tabular}{c|c}
	$\widehat{\mathbf{c}}$ & $\widehat{\mathbf{A}}$ \\
	\hline \\ [-0.3cm]
	& $\widehat{\mathbf{b}}^{\top}$
    \end{tabular}
    \qquad
    \begin{tabular}{c|c}
	$\mathbf{c}$ & $\mathbf{A}$ \\
	\hline \\ [-0.3cm]
	& $\mathbf{b}^{\top}$
    \end{tabular}
\end{center}
with $\widehat{\mathbf{A}} = \cpth{\hat{a}_{lm}}_{l,m = 1,\dots,s}, \widehat{\mathbf{b}} = \cpth{\hat{b}_{l}}_{l = 1,\dots,s}$, and $\widehat{\mathbf{c}} = \cpth{\hat{c}_{l}}_{l = 1,\dots,s}$ denoting the explicit companion method of the IMEX-RK method, and $\mathbf{A} = \cpth{a_{lm}}_{l = 1,\dots,s}, \mathbf{b} = \cpth{b_{l}}_{l = 1,\dots,s}, \mathbf{c} = \cpth{c_{l}}_{l = 1,\dots,s}$ denoting the implicit one. Here $s$ denotes the number of stages of both explicit and implicit methods. Hence, an IMEX-RK method can be interpreted as the combination of two RK methods for which the coefficients $\hat{a}_{lm}, a_{lm}, \hat{c}_{l}, c_{l}, \hat{b}_{l}$, and $b_{l}$ are determined so that the method is consistent of a given order. A detailed discussion on the coupling conditions for high-order IMEX-RK methods can be found in~\cite{pareschi:2005}. In particular, the following relation has to be satisfied~\cite{kennedy:2003}
\begin{equation*}\label{eq:first_order_consistency_RK}
    \sum_{l=1}^{s}b_{l} = \sum_{l=1}^{s}\hat{b}_{l} = 1.
\end{equation*}
Moreover, the classical relations of RK methods for the quadrature nodes $\mathbf{c}$ and $\widehat{\mathbf{c}}$~\cite{butcher:2008, kennedy:2016}, i.e.
\begin{equation}\label{eq:simplyfing_nodes_conditions_RK}
    \sum_{m=1}^{s}a_{lm} = c_{l}, \qquad \sum_{m=1}^{s}\hat{a}_{lm} = \hat{c}_{l}, \qquad l=1,\dots,s,
\end{equation}
are assumed. %Following~\cite{boscarino:2024}, we adopt the following

\begin{definition}[IMEX-RK methods of type I/II/ARS~\cite{boscarino:2024}]
\label{def:IMEX_type_I_II}
    An IMEX-RK method is said to be of \textbf{type I}~\cite{pareschi:2005} if the matrix $\mathbf{A}$ is invertible. It is said to be of \textbf{type II}~\cite{kennedy:2003} if the matrix $\mathbf{A}$ can be written in the form
    $$\mathbf{A} = 
    \begin{pmatrix}
	0 & 0 \\
	\mathbf{a} & \boldsymbol{\mathcal{A}}
    \end{pmatrix},$$
    with $\mathbf{a} = (a_{21},\dots,a_{s1})^{\top} \in \mathbb{R}^{s-1}$ and the matrix $\boldsymbol{\mathcal{A}} \in \mathbb{R}^{s-1 \times s-1}$ is invertible. In the special case $\mathbf{a} = \mathbf{0}$, $b_{1} = 0$, the method is said of \textbf{type ARS} (see~\cite{ascher:1997}) and it is reducible to a method using $s-1$ stages.
\end{definition}

The formal accuracy imposed by the order conditions is typically not preserved in the case of very stiff problems~\cite{boscarino:2017, hairer:1996, pareschi:2005} and one can incur in the so-called order reduction phenomenon~\cite{biswas:2025, hairer:1996}. Hence, other properties have to be imposed to IMEX-RK methods so as to preserve full accuracy in a very stiff regime. 

\begin{definition}[Implicitly and global stiffly accurate methods~\cite{boscarino:2017}]
    An IMEX-RK method is said to be \textbf{Implicitly Stiffly Accurate} or simply \textbf{Stiffly Accurate} (\textbf{SA}) if the implicit companion method is stiffly accurate, i.e. $a_{si} = b_{i}, i = 1,\dots,s$. An IMEX-RK method is said \textbf{Globally Stiffly Accurate} if it is Implicitly Stiffly Accurate and, in addition, the explicit companion method is First Same As Last, i.e. $\hat{a}_{si} = \hat{b}_{i}, i = 1,\dots,s$.
\end{definition}

Next, we recall some relevant properties that a time discretization method should obey for the problem of interest and their extension to IMEX-RK methods. The minimal requirement for the implicit method is to be \texttt{A}-stable. A RK method with a Butcher tableau given by $\rpth{\mathbf{A}, \mathbf{b}, \mathbf{c}}$ is called \texttt{A}-stable if $\left|R(z)\right| \le 1$ for $\text{Re}(z) \le 0$~\cite{kennedy:2016}, where $R(z)$ is the stability function
\begin{equation*}\label{eq:RK_stability_function}
  R(z) = 1 + z\mathbf{b}^{\top}\rpth{\mathbf{I} - z\mathbf{A}}^{-1}\mathbf{e},
\end{equation*}
with $\mathbf{e} = \rpth{1,1,\dots,1} \in \mathbb{R}^{s}$. The asymptotic properties of IMEX-RK methods are strongly related to the \texttt{L}-stability of the implicit part of the scheme. An implicit Runge--Kutta scheme is said to be \texttt{L}-stable~\cite{hairer:1996} if it is \texttt{A}-stable and $R(z) \to 0$ as $z \to \infty$. Since the seminal work of~\cite{pareschi:2005}, the use of \texttt{L}-stable methods in stiff nonlinear regimes has been widely regarded as the appropriate choice, as it endows the resulting scheme with intrinsic robustness against high-frequency disturbances that may arise. An \texttt{L}-stable method typically results from the combination of an \texttt{A}-stable method with a SA method~\cite{hairer:1996}. However, for methods of type II where $\mathbf{A}$ is not invertible, this combination does not necessarily lead to an \texttt{L}-stable method~\cite{boscarino:2009_IMEX_relaxation}. Then, the following supplementary condition is required to obtain the \texttt{L}-stability~\cite{boscarino:2009_IMEX_relaxation}, i.e.
\begin{equation}\label{eq:Lstab_type_II}
    \mathbf{e}^{\top}_{s}\boldsymbol{\mathcal{A}}^{-1}\mathbf{a}  = \sum_{m=2}^{s}w_{sm}{a}_{m1} = 0,
\end{equation}
where $\mathbf{e}^{\top}_{s} = \rpth{0,\dots,0,1}^{\top}$ and $w_{lm}$ denotes the elements of the inverse of $\boldsymbol{\mathcal{A}}$. In addition, \texttt{L}-stability alone is not sufficient to guarantee the desired order of accuracy when the method is applied to algebraic systems. In this case, supplementary conditions must be satisfied. Specifically, for type II schemes, one has to impose~\cite{boscarino:2009_IMEX_relaxation, kennedy:2019}
\begin{equation}\label{eq:supplementary_algebraic_condition}
    \sum_{m=2}^{s}b_{l}w_{lm}w_{m2} = 0, \qquad l=2,\dots,s,
\end{equation}
so as to ensure at least second-order accuracy for algebraic systems. For SA schemes, this condition reduces to $w_{s2}=0$.
\par
Another relevant property is related to the so-called \emph{stage order}, i.e., the order of accuracy of each internal stage. As discussed in~\cite{ketcheson:2009}, the stage order is the minimum order of accuracy when a scheme is applied to stiff problems. It is defined as the maximum $\hat{p}$ for which the following condition holds~\cite{boscarino:2007}
\begin{equation}\label{eq:stage_order}
    \sum_{m=1}^{s}a_{lm}c_{m}^{k - 1} = \frac{1}{k}c_{l}^{k}, \qquad 1 \le k \le \hat{p} \text{ and for all } l.
\end{equation}

%%%%%%%%%%%%%%% Dissipation and dispersion %%%%%%%%%%%%%%
\subsubsection*{Numerical dissipation and dispersion}
\label{ssec:dissipation_dispersion}

In this section, we recall several concepts related to dissipation and dispersion in the context of IMEX-RK methods. When dealing with wave phenomena, especially highly oscillatory and travelling waves, it is crucial to ensure low numerical dissipation in order to accurately capture wave amplitudes and low numerical dispersion to correctly reproduce phase speed and front propagation. A phase-lag analysis of RK methods is typically based on linear ODEs with oscillatory solutions (pure imaginary eigenvalues). More specifically, one considers the following test equation~\cite{franco:1997, vanderhouwen:1989}
\begin{equation*}\label{eq:test_equation_wave}
    \odif{Y}{t} = i\omega Y + \delta\exp\rpth{i\omega_{p}t},
\end{equation*}
where $i$ denotes the imaginary unit, whereas the natural frequency $\omega$, the forcing frequency $\omega_{p} \neq \omega$, and the forcing amplitude $\delta$ are real numbers. To characterize dissipation properties, let $R(z) = P(z)/Q(z)$ and define the $E$-polynomial~\cite{norsett:1975}
\begin{equation}\label{eq:E_polynomial}
    E(y) = \left|Q(iy)\right|^{2} - \left|P(iy)\right|^{2}, \qquad y \in \mathbb{R}.
\end{equation}

\begin{definition}\label{def:RK_dissipation_dispersion}
    A Runge--Kutta method is said to be dispersive of order $q$ ($q$ even) if $\phi(\omega\Delta t) = \omega\Delta t - \arg\spth{R(i\omega\Delta t)} = \mathcal{O}\rpth{\rpth{\omega\Delta t}^{q+1}}$, with $\Delta t > 0$ denoting the time step, or if, equivalently,
    $$\sum_{i=0}^{j}\frac{(-1)^{i}r_{i}}{\rpth{j-1}!} = 0, \qquad j = 1,3,\dots,q-1,$$
    where $r_{i}$ are the coefficients of the Taylor expansion of $R(z)$, i.e. $R(z) = \sum\limits_{i=0}^{\infty}\,r_{i}z^{i}$.
    A Runge--Kutta method is said to be dissipative of order $r$ ($r$ odd) if $d(\omega\Delta t) = 1 - \left|R(i\omega\Delta t)\right| = \mathcal{O}\rpth{\rpth{\omega\Delta t}^{r+1}},$   or if, equivalently, $E_{2j} = 0$ for $0 \le j \le \frac{1}{2}(r - 1)$, where $E_{2j}$ denotes the coefficients of the \texttt{E}-polynomial~\eqref{eq:E_polynomial}.
\end{definition}

%%%%%% Semi-discrete IMEX-RK time discretization %%%%%%%
\subsection{IMEX-RK time discretization for diffusion-reaction PDEs}
\label{ssec:IMEX_RK_model}

We are now in a position to introduce the semi-discrete IMEX-RK time discretization for problem~\eqref{eq:algebraic}. IMEX-RK methods are widely employed for time-dependent problems that can be decomposed into the sum of a stiff and a non-stiff component. Alternatively, they are commonly formulated by splitting the right-hand side of the resulting ODE into a linear term, treated implicitly, and a nonlinear term, discretized explicitly~\cite{calvo:2001, corti:2024_IMEX}. In either formulation, such IMEX-RK methods are constructed in an \emph{additive} fashion and therefore fall within the broader class of additive RK schemes~\cite{kennedy:2003}.
\par
Consider therefore the time-dependent problem $\odif{\bm{y}}{t} = \bm{f}_{\mathrm{I}}(\bm{y},t) + \bm{f}_{\mathrm{E}}(\bm{y},t)$, where $\mathrm{I}$ and $\mathrm{E}$ subscripts denote the components of the system discretized implicitly and explicitly, respectively. If $\bm{v}^{\timesuperscript} \approx \bm{y}(t^{\timesuperscript})$, the generic $l-$stage IMEX-RK method reads as
\begin{equation}\label{eq:IMEX_RK_generic_stage}
    \bm{v}^{(\timesuperscript,l)} = \bm{v}^{\timesuperscript} + \Delta t \spth{\sum_{m=1}^{l-1} \hat{a}_{lm}\bm{f}_{\mathrm{E}}\rpth{\bm{v}^{(\timesuperscript,m)}, t + \hat{c}_{m}\Delta t} + \sum_{m=1}^{s} a_{lm} \bm{f}_{\mathrm{I}}\rpth{\bm{v}^{(\timesuperscript,m)}, t + c_{m}\Delta t}}
\end{equation}
where $l=1,\dots,s$. After the computation of the intermediate stages, the final update $\bm{v}^{\timesuperscript+1}$ is computed as
\begin{equation}\label{eq:IMEX_RK_update}
  \bm{v}^{\timesuperscript+1} = \bm{v}^{\timesuperscript} + \Delta t \sum_{l=1}^{s}\spth{\hat{b}_{l}\bm{f}_{\mathrm{E}}\rpth{\bm{v}^{(\timesuperscript,l)}, t + \hat{c}_{l}\Delta t} + b_{l}\bm{f}_{\mathrm{I}}\rpth{\bm{v}^{(\timesuperscript,l)}, t + c_{l}\Delta t}}.
\end{equation}
\par
We now apply~\eqref{eq:IMEX_RK_generic_stage}-\eqref{eq:IMEX_RK_update} to~\eqref{eq:model_equation}, then assuming therefore that $\mathbf{F}(\bm{u},t) = \mathbf{F}_{\mathrm{E}}(\bm{u},t) + \mathbf{F}_{\mathrm{I}}(\bm{u},t)$, it yields the following expression for the generic $l$-stage
\begin{alignat}{2}\label{eq:IMEX_RK_generic_stage_model_equation}
    \bm{u}^{(\timesuperscript,l)} = \bm{u}^{\timesuperscript} &+ \Delta t\sum_{m=1}^{s}a_{lm}\dive\rpth{\mathbf{\Sigma}\grad\bm{u}^{(\timesuperscript,m)}} \nonumber \\
    &+ \Delta t\spth{\sum_{m=1}^{l-1}\hat{a}_{lm}\mathbf{F}_{\mathrm{E}}\rpth{\bm{u}^{(\timesuperscript,m)}, t + \hat{c}_{m}\Delta t} + \sum_{m=1}^{s}a_{lm}\mathbf{F}_{\mathrm{I}}\rpth{\bm{u}^{(\timesuperscript,m)}, t + c_{m}\Delta t}},
\end{alignat}
whereas the final update reads
\begin{equation}\label{eq:IMEX_RK_update_model_equation}
    \bm{u}^{\timesuperscript+1} = \bm{u}^{\timesuperscript} + \Delta t\sum_{l=1}^{s}\spth{b_{l}\dive\rpth{\mathbf{\Sigma}\grad\bm{u}^{(\timesuperscript,l)}} + \hat{b}_{l}\mathbf{F}_{\mathrm{E}}\rpth{\bm{u}^{(\timesuperscript,l)}, t + \hat{c}_{l}\Delta t} + b_{l}\mathbf{F}_{\mathrm{I}}\rpth{\bm{u}^{(\timesuperscript,l)}, t + c_{l}\Delta t}}.
\end{equation}
We observe that the diffusion term is treated implicitly. Throughout this work, we assume that the diffusion coefficients $\sigma_{j}$ are constant or, at least, independent of $\bm{u}$. We point out that, if the diffusion coefficients depend on space, the same algorithm remains applicable. If they depend on time, it is  sufficient to replace $\mathbf{\Sigma}$ with $\mathbf{\Sigma}^{(n,m)}$ in~\eqref{eq:IMEX_RK_generic_stage_model_equation} and with $\mathbf{\Sigma}^{(n,l)}$ in~\eqref{eq:IMEX_RK_update_model_equation}. We will briefly comment on the assumption of independence from $\bm{u}$ further later in the manuscript.

\begin{remark}[Application of IMEX-RK to Equation~\eqref{eq:algebraic_ode}] Assuming that $\mathbf{M}$ is invertible, the semidiscrete problem~\eqref{eq:algebraic_ode} can be written in the standard additive form as
$$\dot{\bm{U}}(t) = \bm{f}_{\mathrm{E}}(\bm{U},t) + \bm{f}_{\mathrm{I}}(\bm{U},t),$$
where
$$\bm{f}_{\mathrm{E}}(\bm{U},t) = \mathbf{M}^{-1}\bm{G}_{\mathrm{E}}(\bm{U},t), \qquad \bm{f}_{\mathrm{I}}(\bm{U},t) = \mathbf{M}^{-1} \spth{-\mathbf{K}\bm{U} + \bm{G}_{\mathrm{I}}(\bm{U},t)}.$$
%.
Hence, by applying the relations~\eqref{eq:IMEX_RK_generic_stage_model_equation}-\eqref{eq:IMEX_RK_update_model_equation} to Equation~\eqref{eq:algebraic_ode}, the fully discrete generic $l$-stage reads
\begin{alignat}{2}
    \boldsymbol{U}^{(\timesuperscript,l)} = \boldsymbol{U}^{\timesuperscript} &- \mathbf{M}^{-1}\Delta t\sum_{m=1}^{s}a_{lm}\mathbf{K}\boldsymbol{U}^{(\timesuperscript,m)} \nonumber \\
    &+ \mathbf{M}^{-1}\Delta t\spth{\sum_{m=1}^{l-1}\hat{a}_{lm}\boldsymbol{G}_{\mathrm{E}}\rpth{\boldsymbol{U}^{(\timesuperscript,m)}, t + \hat{c}_{m}\Delta t} + \sum_{m=1}^{s}a_{lm}\boldsymbol{G}_{\mathrm{I}}\rpth{\boldsymbol{U}^{(\timesuperscript,m)}, t + c_{m}\Delta t}},
\end{alignat}
whereas the final update reads
\begin{alignat}{2}\label{eq:IMEX_RK_update_fully_discrete_model_equation}
    \boldsymbol{U}^{\timesuperscript+1} = \boldsymbol{U}^{\timesuperscript} &- \mathbf{M}^{-1}\Delta t\sum_{l=1}^{s}b_{l}\mathbf{K}\boldsymbol{U}^{(\timesuperscript,l)} \nonumber \\
    &+ \mathbf{M}^{-1}\Delta t \sum_{l=1}^{s}\spth{\hat{b}_{l}\boldsymbol{G}_{\mathrm{E}}\rpth{\boldsymbol{U}^{(\timesuperscript,l)}, t + \hat{c}_{l}\Delta t} + b_{l}\boldsymbol{G}_{\mathrm{I}}\rpth{\boldsymbol{U}^{(\timesuperscript,l)}, t + c_{l}\Delta t}}.
\end{alignat}
Here $\boldsymbol{U}$ denotes the vector of the degrees of freedom, whereas $\boldsymbol{G}_{\mathrm{E}}$ and $\boldsymbol{G}_{\mathrm{I}}$ denote the spatial approximation of $\mathbf{F}_{\mathrm{E}}$ and $\mathbf{F}_{\mathrm{I}}$, respectively.
\end{remark}

%%%%%% Semi-discrete SI-IMEX-RK time discretization %%%%
\subsection{SI-IMEX-RK time discretization for diffusion-reaction PDEs}
\label{ssec:SI_IMEX_RK_method}

Additive RK methods are commonly adopted for systems admitting a decomposition into the sum of stiff and non-stiff terms. However, the stiff dependence can be (highly) nonlinear, thus requiring the solution of a nonlinear system at each stage. In other cases, the stiffness may be associated only to some variables~\cite{boscarino:2016}. Let us consider a \emph{partitioned} time-dependent system
\begin{subequations}
\begin{alignat}{2}
    \odif{\bm{y}_{\mathrm{E}}}{t} &= \bm{f}_{\mathrm{E}}(\bm{y}_{\mathrm{E}},\bm{y}_{\mathrm{I}},t), \\
    \odif{\bm{y}_{\mathrm{I}}}{t} &= \frac{1}{\varepsilon}\bm{f}_{\mathrm{I}}(\bm{y}_{\mathrm{E}},\bm{y}_{\mathrm{I}},t), \label{eq:ODE_stiff_partitioned}
\end{alignat}
\end{subequations}
where $\bm{y}_{\mathrm{E}}$ and $\bm{y}_{\mathrm{I}}$ may be vectors of different dimensions and $\bm{y}_{\mathrm{E}}(0) = \bm{y}_{\mathrm{E},0}, \bm{y}_{\mathrm{I}}(0) = \bm{y}_{\mathrm{I},0}$ are the initial conditions. The stiffness of the system is associated to the variables $\bm{y}_{\mathrm{I}}$ under the assumption $0 < \varepsilon \ll 1$.

We now introduce a class of Semi-Implicit IMEX-RK (SI-IMEX-RK) methods for problem~\eqref{eq:model}, in a sense that will be made precise below. Following the approach proposed in~\cite{boscarino:2016}, we rewrite~\eqref{eq:model_equation} as the following partitioned system
\begin{subequations}\label{eq:partitioned_model}
\begin{alignat}{2}
    \pad{\bm{u}_{\mathrm{E}}}{t} &= \dive\rpth{\mathbf{\Sigma}\grad\bm{u}_{\mathrm{I}}} + \mathbf{F}\rpth{\bm{u}_{\mathrm{E}},\bm{u}_{\mathrm{I}},t}, \\
    \pad{\bm{u}_{\mathrm{I}}}{t} &= \dive\rpth{\mathbf{\Sigma}\grad\bm{u}_{\mathrm{I}}} + \mathbf{F}\rpth{\bm{u}_{\mathrm{E}},\bm{u}_{\mathrm{I}},t},
\end{alignat}
\end{subequations}
with initial conditions $\bm{u}_{\mathrm{E}} = \bm{u}_{0}, \bm{u}_{\mathrm{I}} = \bm{u}_{0}$ and supplied with compatible boundary conditions. Subscripts $\mathrm{E}$ and $\mathrm{I}$ in~\eqref{eq:partitioned_model} indicate the explicit and (semi-)implicit
treatment of the first and the second term of $\mathbf{F}$, respectively. In this setting, the solution of~\eqref{eq:partitioned_model} satisfies $\bm{u}_{\mathrm{E}} = \bm{u}_{\mathrm{I}}$ for any $t \ge 0$ and is also a solution of~\eqref{eq:model}.

The SI-IMEX-RK framework is obtained by applying a partitioned IMEX Runge--Kutta method to system~\eqref{eq:partitioned_model}. For simplicity, and following~\cite{boscarino:2016}, we assume that the system is autonomous, i.e. $\mathbf{F}\rpth{\bm{u}_{\mathrm{E}},\bm{u}_{\mathrm{I}},t} = \mathbf{F}\rpth{\bm{u}_{\mathrm{E}},\bm{u}_{\mathrm{I}}}$. Extensions to non-autonomous problems, as well as to the case $\mathbf{c} \neq \widehat{\mathbf{c}}$, would unnecessarily complicate the presentation of the core ideas; we refer to~\cite{boscarino:2016} for a detailed discussion. A crucial assumption is instead that the weights of the implicit and explicit tableaux coincide, i.e. $\mathbf{b}^{\top} = \widehat{\mathbf{b}}^{\top}$. The relevance of this condition will become clear shortly. However, we emphasize that this requirement is important independently of the SI-IMEX-RK construction, as it is necessary to ensure the preservation of linear invariants~\cite{giraldo:2013}. Under these assumptions, the generic $l$-stage of the SI-IMEX-RK method applied to~\eqref{eq:model} reads as in Algorithm \ref{alg:si_imex_rk}.

\begin{algorithm}[t]
\caption{SI-IMEX-RK scheme for~\eqref{eq:model}}
\label{alg:si_imex_rk}
\begin{algorithmic}
    \Given $\bm{u}^{0} = \bm{u}_{0}$, final time $T$, time step $\Delta t$, Butcher tableaux $(a_{lm}, \hat{a}_{lm}, b_l)$
    \For{$l = 1, \dots, s$}
        \State \textbf{Step 1 (explicit stage).} Compute
        \begin{equation}\label{eq:SI_IMEX_RK_generic_stage_model_equation_1}
            \bm{u}_{\mathrm{E}}^{(\timesuperscript,l)} = \bm{u}^{\timesuperscript} + \Delta t \sum_{m=1}^{l-1} \hat{a}_{lm}\,\dive\rpth{\mathbf{\Sigma}\nabla\bm{u}_{\mathrm{I}}^{(\timesuperscript,m)}} + \Delta t \sum_{m=1}^{l-1} \hat{a}_{lm}\,\mathbf{F}\!\rpth{\bm{u}_{\mathrm{E}}^{(\timesuperscript,m)},\bm{u}_{\mathrm{I}}^{(\timesuperscript,m)}}.
        \end{equation}
        \State \textbf{Step 2 (implicit stage).} Solve for $\bm{u}_{\mathrm{I}}^{(n,l)}$:
        \begin{equation}\label{eq:SI_IMEX_RK_generic_stage_model_equation_2}
            \bm{u}_{\mathrm{I}}^{(\timesuperscript,l)} = \bm{u}^{\timesuperscript} + \Delta t \sum_{m=1}^{s} a_{lm}\,\dive\rpth{\mathbf{\Sigma}\nabla\bm{u}_{\mathrm{I}}^{(\timesuperscript,m)}} + \Delta t \sum_{m=1}^{s} a_{lm}\,\mathbf{F}\rpth{\bm{u}_{\mathrm{E}}^{(\timesuperscript,m)},\bm{u}_{\mathrm{I}}^{(\timesuperscript,m)}}.
        \end{equation}
    \EndFor
    \State \textbf{Step 3 (update).} Set
        \begin{equation}\label{eq:SI_IMEX_RK_generic_stage_model_equation_3}
            \bm{u}^{\timesuperscript+1} = \bm{u}^{\timesuperscript} + \Delta t \sum_{l=1}^{s} b_{l}\,\dive\rpth{\mathbf{\Sigma}\nabla\bm{u}_{\mathrm{I}}^{(\timesuperscript,l)}} + \Delta t \sum_{l=1}^{s} b_{l}\,\mathbf{F}\rpth{\bm{u}_{\mathrm{E}}^{(\timesuperscript,l)},\bm{u}_{\mathrm{I}}^{(\timesuperscript,l)}}.
        \end{equation}
    \end{algorithmic}
\end{algorithm}

\begin{remark}\label{rem:SI_IMEX_DIRK}
If the implicit tableau $\mathbf{A}$ is not lower triangular in Algorithm~\ref{alg:si_imex_rk}, the stages can no longer be computed sequentially: $\boldsymbol{U}_{\mathrm{I}}^{(n,l)}$ depends on $\boldsymbol{U}_{\mathrm{E}}^{(n,m)}$ for $m>l$, which in turn depend on $\boldsymbol{U}_{\mathrm{I}}^{(n,l)}$. All explicit and implicit stages thus become coupled into a single nonlinear system. Diagonally Implicit Runge--Kutta (DIRK) methods, such as those discussed in Section~\ref{sec:optmized_IMEX} and for which $\mathbf{A}$ is lower triangular, have to be employed to preserve the sequential evaluation.
\end{remark}

Owing to the condition $\mathbf{b}^{\top} = \widehat{\mathbf{b}}^{\top}$, we obtain $\bm{u}_{E}^{n+1} = \bm{u}_{I}^{n+1} = \bm{u}^{n+1}$ for all $n$, showing that the duplication of variables is purely formal.

The main advantage of this approach lies in the identification of a suitable (possibly linear) implicit dependence, which avoids the solution of (highly) nonlinear systems or the need for stage-by-stage linearization. Similarly to the discussion in Section 3.1, the algorithm can also be readily applied when $\mathbf{\Sigma}$ depends on space and time. In the presence of nonlinear diffusion, as encountered, for example, in cross-diffusion systems \cite{gomez:2024}, where $\mathbf{\Sigma} = \mathbf{\Sigma}(\bm{u})$, the same framework can be naturally extended by evaluating the diffusion coefficient explicitly as $\mathbf{\Sigma}(\bm{u}_{\mathrm{E}})$, thereby preserving a linear implicit solve at each stage in the case of DIRK methods.

\begin{remark}[Application of IMEX-RK to Equation~\eqref{eq:algebraic_ode}] The $l$-stage of the SI-IMEX-RK method in Algorithm \ref{alg:si_imex_rk} for problem~\eqref{eq:algebraic_ode} reads
\begin{subequations}\label{eq:SI_IMEX_RK_generic_stage_model_equation_discrete}
\begin{alignat}{2}
    \boldsymbol{U}_{E}^{(\timesuperscript,l)} &= \boldsymbol{U}^{\timesuperscript} - \mathbf{M}^{-1}\Delta t\sum_{m=1}^{l-1}\hat{a}_{lm}\mathbf{K}\boldsymbol{U}_{I}^{(\timesuperscript,m)} + \mathbf{M}^{-1}\Delta t \sum_{m=1}^{l-1}\hat{a}_{lm}\boldsymbol{G}\rpth{\boldsymbol{U}_{\mathrm{E}}^{(\timesuperscript,m)},\boldsymbol{U}_{\mathrm{I}}^{(\timesuperscript,m)}}, \\
    \boldsymbol{U}_{I}^{(\timesuperscript,l)} &= \boldsymbol{U}^{\timesuperscript} - \mathbf{M}^{-1}\Delta t\sum_{m=1}^{s}a_{lm}\mathbf{K}\boldsymbol{U}_{I}^{(\timesuperscript,m)} + \mathbf{M}^{-1}\Delta t \sum_{m=1}^{s}a_{lm}\boldsymbol{G}\rpth{\boldsymbol{U}_{\mathrm{E}}^{(\timesuperscript,m)},\boldsymbol{U}_{\mathrm{I}}^{(\timesuperscript,m)}},
\end{alignat}
\end{subequations}
whereas the final update reads
\begin{equation}
    \boldsymbol{U}^{\timesuperscript+1} = \boldsymbol{U}^{\timesuperscript} - \mathbf{M}^{-1}\Delta t\sum_{l=1}^{s}b_{l}\mathbf{K}\boldsymbol{U}_{I}^{(\timesuperscript,l)} + \mathbf{M}^{-1}\Delta t \sum_{l=1}^{s} b_{l}\boldsymbol{G}\rpth{\boldsymbol{U}_{\mathrm{E}}^{(\timesuperscript,m)},\boldsymbol{U}_{\mathrm{I}}^{(\timesuperscript,m)}}.
\end{equation}
Here $\boldsymbol{U}_{E}$ and $\boldsymbol{U}_{I}$ denote the vectors of the degrees of freedom for the explicit and (semi)-implicit components, respectively.
\end{remark}

%%%%%%%%%%%%%%%%%%%%%%%%%%%%%%%%%%%%%%%%%%%
%%%%%%%%%%% Optimized IMEX-RK %%%%%%%%%%%%%
%%%%%%%%%%%%%%%%%%%%%%%%%%%%%%%%%%%%%%%%%%%
\section{Optimized IMEX-RK for degenerate diffusion-reaction PDEs}
\label{sec:optmized_IMEX}

In this section, we propose third and fourth-order IMEX-RK schemes that optimize the properties discussed in Section~\ref{sec:IMEX} and are therefore particularly well suited for model~\eqref{eq:model}. The symbolic computations have been carried out using both \texttt{MATLAB} and \texttt{sympy}~\cite{meurer:2017} in combination with \texttt{ponio}~\cite{dubois:2025, massot:2024} for the symbolic derivation of the stability functions. For the novel schemes derived in this work, we adopt the notation IMEX$(p,s)$–X, where $p$ denotes the order of the method, and $s$ represents the number of stages of both the implicit and explicit schemes. The symbol X denotes a key distinguishing property of the method. As will become clear in the following, all the proposed methods achieve stage order 2, and their implicit part is always \texttt{L}-stable. Since these two features constitute essential requirements in the derivation, we do not include them in the nomenclature in order to keep it concise.
\par
First, we assume that the implicit scheme is a Explicit Singly Diagonally Implicit Runge--Kutta (ESDIRK) method, namely $\mathbf{A}$ is lower triangular, $a_{ll} = \gamma > 0$ with $l=2,\dots,s$, and the first row is zero. Hence, the first stage is only formal, unless necessary to achieve stage order 2 (see relation~\eqref{eq:stage_order}). Moreover, as discussed in Remark~\ref{rem:SI_IMEX_DIRK}, this is necessary to ensure sequential evaluation for the SI-IMEX-RK method and to obtain linear systems whenever the implicit dependence is linear. Finally, to guarantee full accuracy in the stiff regime, we consider a \texttt{L}-stable method, as discussed in Section~\ref{sec:IMEX}. A similar analysis for ESDIRK schemes has been conducted in~\cite{kvaerno:2004}. Then, starting from the implicit method, an explicit companion method is obtained in order to satisfy the order and coupling conditions~\cite{pareschi:2005}.

Several alternative IMEX strategies have been proposed in the literature for diffusion-reaction problems, (i.e.,~\cite{abdulle:2013, hundsdorfer:2007, izzo:2017}). In particular, we mention the PIROCK method~\cite{abdulle:2013}, where second-order SDIRK schemes are coupled with explicit second-order ROCK integrators~\cite{abdulle:2002, abdulle:2001}, facilitating the inclusion of non-trivial, possibly nonlinear operators. However, extending PIROCK to higher order is not straightforward. The semi-implicit approach in Section~\ref{ssec:SI_IMEX_RK_method} provides an alternative strategy.% A detailed comparison with the aforementioned approaches is far beyond the scope of the present work and is left for future investigation.

%%%%%%%%%%%% Third-order schemes %%%%%%%%%
\subsection{Third-order IMEX-RK schemes}
\label{ssec:third_order_schemes}

In this section, we derive some third-order time discretization schemes that satisfy the largest possible number of the properties described in Section~\ref{sec:IMEX}. All the derived schemes are assumed to be \texttt{L}-stable and of stage order 2 in their implicit part. Moreover, we impose the supplementary condition~\eqref{eq:supplementary_algebraic_condition}, together with high-order dissipation or dispersion requirements and additional order conditions for improved accuracy when solving algebraic systems. We recall that the considered IMEX-RK schemes satisfy condition~\eqref{eq:simplyfing_nodes_conditions_RK}. The conditions to achieve third-order accuracy for a RK method $\rpth{\mathbf{A}, \mathbf{b}, \mathbf{c}}$ read as follows
\begin{subequations}
\begin{alignat}{3}
    &\sum_{l=1}^{s}b_{l} = 1, \label{eq:order_conditions_third_order_partI} & & \qquad \sum_{l=1}^{s}b_{l}c_{l} = \frac{1}{2}, \\
    &\sum_{l=1}^{s}b_{l}c_{l}^{2} = \frac{1}{3}, & & \qquad  \sum_{l,m=1}^{s}b_{l}a_{lm}c_{m} = \frac{1}{6}. \label{eq:order_conditions_third_order_partII}
\end{alignat}
\end{subequations}
The coupling conditions to obtain a third order IMEX-RK method are
\begin{alignat*}{7}
    &\sum_{l=1}^{s}b_{l}\hat{c}_{l} = \frac{1}{2}, & & \qquad \sum_{l=1}^{s}\hat{b}_{l}c_{l} = \frac{1}{2}, & & \qquad \sum_{l=1}^{s}\hat{b}_{l}c_{l}^{2} = \frac{1}{3}, & & \qquad \sum_{l=1}^{s}\hat{b}_{l}\hat{c}_{l}c_{l} = \frac{1}{3}, \\
    & \sum_{l=1}^{s}b_{l}\hat{c}_{l}^{2} = \frac{1}{3}, & & \qquad \sum_{l=1}^{s}b_{l}\hat{c}_{l}c_{l} = \frac{1}{3}, & & \qquad
    \sum_{l,m=1}^{s}b_{l}\hat{a}_{lm}c_{m} = \frac{1}{6}, & & \qquad \sum_{l,m=1}^{s}b_{l}a_{lm}\hat{c}_{m} = \frac{1}{6}, \\
    & \sum_{l,m=1}^{s}b_{l}\hat{a}_{lm}\hat{c}_{m} = \frac{1}{6}, & & \qquad  \sum_{l,m=1}^{s}\hat{b}_{l}\hat{a}_{lm}c_{m} = \frac{1}{6}, & & \qquad \sum_{l,m=1}^{s}\hat{b}_{l}a_{lm}\hat{c}_{m} = \frac{1}{6}, & & \qquad \sum_{l,m=1}^{s}\hat{b}_{l}a_{lm}c_{m} = \frac{1}{6}.
\end{alignat*}
One can easily notice that $12$ coupling conditions are necessary in the general case $\mathbf{b}^{\top} \neq \widehat{\mathbf{b}}^{\top}, \mathbf{c} \neq \widehat{\mathbf{c}}$ which reduce to $3$ coupling conditions if $\mathbf{b}^{\top} = \widehat{\mathbf{b}}^{\top}$ or to $2$ coupling condition if $\mathbf{c} = \widehat{\mathbf{c}}$. If both $\mathbf{b}^{\top} = \widehat{\mathbf{b}}^{\top}$ and $\mathbf{c} = \widehat{\mathbf{c}}$, zero coupling conditions are necessary. Finally, the additional order conditions up to third order that guarantee better accuracy for an algebraic system are~\cite{boscarino:2009_IMEX_relaxation}
\begin{equation}\label{eq:index_1_DAE_conditions_order_3}
    \sum_{l,m=2}^{s}b_{l}w_{lm}\hat{c}_{m} = 1, \qquad \sum_{l,m=2}^{s}b_{l}w_{lm}\hat{c}_{m}^{2} = 1, \qquad \sum_{l,m,k=2}^{s}b_{l}w_{lm}\hat{a}_{mk}\hat{c}_{k} = \frac{1}{2},
\end{equation}
where $w_{lm}$ are the elements of the inverse of the matrix $\boldsymbol{\mathcal{A}}$ (see Definition~\ref{def:IMEX_type_I_II}).

Following the analysis in~\cite{boscarino:2009_third_order}, we focus on five-stage ($s = 5$) IMEX-RK schemes, which represent the minimum number of stages required to simultaneously achieve \texttt{L}-stability and satisfy the supplementary condition~\eqref{eq:supplementary_algebraic_condition}~\cite{boscarino:2009_third_order}. Five-stage stiffly-accurate ESDIRK methods have 10 degrees of freedom (see Table~\ref{tab:implicit_s5_p3}).

\begin{table}[H]
    \centering
    \begin{tabular}{c|ccccc}
	    $0$ & $0$ & $0$ & $0$ & $0$ & $0$ \\
	    $c_{2}$ & $c_{2} - \gamma$ & $\gamma$ & $0$ & $0$ & $0$ \\
        $c_{3}$ & $c_{3} - a_{32} - \gamma$ & $a_{32}$ & $\gamma$ & $0$ & $0$ \\
        $c_{4}$ & $c_{4} - a_{42} - a_{43} - \gamma$ & $a_{42}$ & $a_{43}$ & $\gamma$ & $0$ \\
        $1$ & $1 - a_{52} - a_{53} - a_{54} - \gamma$ & $a_{52}$ & $a_{53}$ & $a_{54}$ & $\gamma$ \\
	    \hline
        & $1 - a_{52} - a_{53} - a_{54} - \gamma$ & $a_{52}$ & $a_{53}$ & $a_{54}$ & $\gamma$
    \end{tabular}
    \caption{Butcher tableau of the implicit five-stage stiffly-accurate ESDIRK method.}
    \label{tab:implicit_s5_p3}
\end{table}
\noindent
The coefficients $c_{2}, a_{32}, a_{42}, a_{43}, a_{52}, a_{53},$ and $a_{54}$, together with the associated determining conditions, are reported in Table~\ref{tab:order_coefficients_fourth_order}.
\begin{table}[h!]
    \centering
    \begin{tabular}{|r@{\hspace{2pt}} l|l|}
        \hline
        \multicolumn{2}{|c|}{\textbf{Coefficient}} & \multicolumn{1}{c|}{\textbf{Condition used}} \\
        \hline
        \rule{0pt}{4.5ex} $a_{52} = $ & $\displaystyle \frac{\frac{1}{2} - \gamma - a_{53}c_{3} - a_{54}c_{4}}{c_{2}}$ & Order conditions~\eqref{eq:order_conditions_third_order_partI}-\eqref{eq:order_conditions_third_order_partII} \\[10pt]
        \hline
        \rule{0pt}{4.5ex} $a_{53} = $ & $\displaystyle \frac{\gamma\rpth{c_{2} - 1} - a_{54}c_{4}\rpth{c_{2} - c_{4}} + \frac{c_{2}}{2} - \frac{1}{3}}{c_{3}\rpth{c_{2} - c_{3}}}$ & Order conditions~\eqref{eq:order_conditions_third_order_partI}-\eqref{eq:order_conditions_third_order_partII} \\[10pt]
        \hline
        \rule{0pt}{4.5ex} $a_{42} = $ & $\displaystyle \frac{\gamma^{2} - \gamma - a_{32}a_{53}c_{2} - a_{43}a_{54}c_{3} + \frac{1}{6}}{a_{54}c_{2}}$ & Order conditions~\eqref{eq:order_conditions_third_order_partI}-\eqref{eq:order_conditions_third_order_partII} \\[10pt]
        \hline
        $c_{2} = $ & $\displaystyle 2\gamma$ & Stage order $2$~\eqref{eq:stage_order} \\
        \hline
        \rule{0pt}{4.5ex} $a_{32} = $ & $\displaystyle \frac{c_{3}\rpth{c_{3} - 2\gamma}}{4\gamma}$ & Stage order $2$~\eqref{eq:stage_order} \\[10pt]
        \hline
        \rule{0pt}{4.5ex} $a_{43} = $ & $\displaystyle \frac{\gamma\rpth{6\gamma^{3} - 18\gamma^{2} + 9\gamma - 1}}{3 a_{54}c_{3}\rpth{2\gamma - c_{3}}}$ & \texttt{L}-stability~\eqref{eq:Lstab_type_II} \\[10pt]
        \hline
        \rule{0pt}{4.5ex} $a_{54} = $ & $\displaystyle \frac{\gamma\rpth{2\gamma - 1 + c_{3}\rpth{1 - \gamma}}}{c_{4}\rpth{c_{3} - c_{4}}}$ & Supplementary algebraic condition~\eqref{eq:supplementary_algebraic_condition} \\[10pt]
       \hline
    \end{tabular}
    \caption{Expressions of the coefficients for the five-stage third-order implicit scheme obtained from order conditions, stage order $2$, \texttt{L}-stability and supplementary algebraic condition \eqref{eq:supplementary_algebraic_condition}.}
    \label{tab:order_coefficients_fourth_order}
\end{table}
For later use, we recall the \texttt{E}-polynomial \eqref{eq:E_polynomial}, which, as reported in~\cite{boscarino:2009_IMEX_relaxation, kennedy:2016}, reduces to
\begin{align}\label{eq:E_polynomial_s5_p3}
    E(y) = \gamma^{8}y^{8} & + \rpth{-\frac{1}{36} + \frac{2}{3}\gamma - 6\gamma^{2} + \frac{76}{3}\gamma^{3} - 52\gamma^{4} + 48\gamma^{5} - 12\gamma^{6}}y^{6} \nonumber \\
    &+ \rpth{\frac{1}{12} - \frac{4}{3}\gamma + 6\gamma^{2} - 8\gamma^{3} + 2\gamma^{4}}y^{4}.
\end{align}
Then, the region of \texttt{A}-stability (hence \texttt{L}-stability) is given by the following interval~\cite{boscarino:2009_IMEX_relaxation, hairer:1996, kennedy:2016}
\begin{equation}\label{eq:L_stability_s5_p3}
    0.2236478009341764510696898 \le \gamma \le 0.5728160624821348554080014.
\end{equation}
Moreover, one can easily verify that the first two index-1 DAE accuracy conditions~\eqref{eq:index_1_DAE_conditions_order_3} are identically satisfied. We are left with 3 free parameters, i.e. $c_{3}, c_{4}$, and $\gamma$. Several choices can be considered to determine $c_{3}$ and $c_{4}$. For instance, one may impose stage order 3 for the intermediate stages, or alternatively analyze the stability region of the internal stages, the so-called internal stability~\cite{kennedy:2016}. The investigation performed suggests that better results are obtained when stage order 3 is imposed on the third and fourth intermediate stages. Moreover, this choice preserves the high-order accuracy as much as possible, thereby avoiding degradation of the high-order accuracy of the spatial discretization (see Section~\ref{sec:polydg}). Finally, an alternative closure for $c_{3}$ will be discussed below when deriving the explicit method. Hence, imposing stage order 3 for the third intermediate stage leads to two admissible values for $c_{3}$
\begin{equation}\label{eq:c3_stage_order3_s5_p3}
    c_{3} = c_{3}^{(1\pm)} = \gamma\rpth{3 \pm \sqrt{3}},
\end{equation}
whereas imposing stage order 3 for the fourth (penultimate) intermediate stage yields two admissible values for $c_{4}$
\begin{equation}\label{eq:c4_stage_order3_s5_p3}
    c_{4} = c_{4}^{(1\pm)} = \scriptstyle \frac{-6\gamma^{3} + 6\gamma^{2}\rpth{c_{3} + 1} - 3\gamma\rpth{2c_{3} + 1} \pm \sqrt{\rpth{6\gamma^{2} - 2c_{3}\rpth{3\gamma - 1} - 1}\rpth{6\gamma^{4} - 6\gamma^{3}\rpth{c_{3} + 2} + 3\gamma^{2}\rpth{6c_{3} - 1} + 2\gamma\rpth{3 - 7c_{3}} + 2c_{3} - 1}} + 1}{2\rpth{\rpth{c_{3} - 2}\rpth{\gamma - 1} - 1}}.
\end{equation}

Finally, we are left with one free parameter $\gamma$. In this work, we consider $\gamma = 0.435866521508482$ used also for the ARS(3,4,3) scheme in~\cite{ascher:1997}, the ARK3(2)4L[2]SA scheme in~\cite{kennedy:2003}, and the MARS(3,4,3) and MARK3(2)4L[2]SA schemes in~\cite{boscarino:2007}, as well as the BHR(5,5,3) in~\cite{boscarino:2009_third_order}. Whenever schemes are taken or adapted from the literature, we adopt the notation of the related reference where available. As an alternative, since as already mentioned~\eqref{eq:model_equation} may be characterized by a travelling-wave solution, we maximize the order of dissipation and dispersion. Regarding the order of dissipation, owing to Definition~\eqref{def:RK_dissipation_dispersion}, we solve $E_{4} = 0$ in~\eqref{eq:E_polynomial_s5_p3} to obtain dissipation order $5$. This gives $\gamma = 0.572816062482135$, which is the limit of the \texttt{I}-stability, hence \texttt{L}-stability~\eqref{eq:L_stability_s5_p3}. A slightly different value, namely $\gamma = 0.57281606248208$ has been considered in~\cite{boscarino:2009_IMEX_relaxation}. Regarding the order of dispersion, one can easily verify that
$$D(j) = \sum_{i=0}^{j}\frac{(-1)^{i}r_{i}}{\rpth{j-1}!} = 0, \qquad \text{for } j = 1,3,$$
where we recall that $r_{i}$ are the coefficients of the Taylor expansion of $R(z)$ (see Definition~\eqref{def:RK_dissipation_dispersion}). Imposing $D(j) = 0$ for $j = 5$ so as to obtain an implicit RK-method with dispersion order equal to $6$ leads to the following equation
\begin{equation*}\label{eq:high_order_dispersion_s5_p3}
    4\gamma^{5} - 16\gamma^{4} + 14\gamma^{3} - \frac{14}{3}\gamma^{2} + \frac{2}{3}\gamma - \frac{1}{30} = 0.
\end{equation*}
The only solution that guarantees \texttt{L}-stability is $\gamma = 0.525721461435005$. Hence, if $\gamma = 0.525721461435005$ the implicit method has a dispersion order of $6$, otherwise its dispersion order is $4$.

Next, we consider an explicit companion method for the IMEX scheme. We consider $\widehat{\mathbf{b}} = \mathbf{b}$ and $\widehat{\mathbf{c}} = \mathbf{c}$, which results in 6 degrees of freedom. The corresponding Butcher tableau then reads as follows (Table~\ref{tab:explicit_s5_p3})

\begin{table}[H]
    \centering
    \begin{tabular}{c|ccccc}
	    $0$ & $0$ & $0$ & $0$ & $0$ & $0$ \\
	    $c_{2}$ & $c_{2}$ & $0$ & $0$ & $0$ & $0$ \\
        $c_{3}$ & $c_{3} - \hat{a}_{32}$ & $\hat{a}_{32}$ & $0$ & $0$ & $0$ \\
        $c_{4}$ & $c_{4} - \hat{a}_{42} - \hat{a}_{43}$ & $\hat{a}_{42}$ & $\hat{a}_{43}$ & $0$ & $0$ \\
        $1$ & $1 - \hat{a}_{52} - \hat{a}_{53} - \hat{a}_{54}$ & $\hat{a}_{52}$ & $\hat{a}_{53}$ & $\hat{a}_{54}$ & $0$ \\
	    \hline
        & $1 - a_{52} - a_{53} - a_{54} - \gamma$ & $a_{52}$ & $a_{53}$ & $a_{54}$ & $\gamma$
    \end{tabular}
    \caption{Butcher tableau of the explicit five-stage method.}
    \label{tab:explicit_s5_p3}
\end{table}

Since $\widehat{\mathbf{b}} = \mathbf{b}$ and $\widehat{\mathbf{c}} = \mathbf{c}$, no coupling condition with the implicit companion method is required to obtain a third order IMEX scheme (see the discussion in Section~\ref{sec:IMEX}). However, we are left with the last condition in~\eqref{eq:order_conditions_third_order_partII} to ensure that the explicit method is of the third order, from which we get
\begin{equation*}
    \hat{a}_{52} = \frac{-\frac{\gamma\hat{a}_{54}c_{4}}{2} - \gamma\rpth{\hat{a}_{32}a_{53} + \hat{a}_{42}a_{54}} - \frac{c_{3}\rpth{\gamma\hat{a}_{53} + \hat{a}_{43}a_{54}}}{2} + \frac{1}{12}}{\gamma^{2}}.
\end{equation*}
Next, we impose the condition~\eqref{eq:index_1_DAE_conditions_order_3} that guarantees third order accuracy for the algebraic variable, from which we obtain
\begin{equation*}
    \hat{a}_{43} = \frac{1 - 3\gamma\rpth{4\hat{a}_{32}a_{53} + 4\hat{a}_{42}a_{54} + 1}}{6a_{54}c_{3}}.
\end{equation*}
Moreover, one can easily verify that this condition is equivalent to imposing stage order 2 for the fifth (last) intermediate stage. Imposing stage order 2 for the third stage yields
\begin{equation}\label{eq:stage_order2_stage2_exp}
    {\hat{a}_{32} = \hat{a}_{32}^{(1)}} = \frac{c_{3}^{2}}{2c_{2}},
\end{equation}
while imposing stage order 2 for the fourth intermediate stage reduces to
\begin{equation}\label{eq:stage_order2_stage3_exp}
    {\hat{a}_{32} = \hat{a}_{32}^{(2)}} = \frac{1 - 3\gamma - 3b_{4}c_{4}^{2}}{12\gamma b_{3}}.
\end{equation}
The only admissible solution which guarantees that both~\eqref{eq:stage_order2_stage2_exp} and~\eqref{eq:stage_order2_stage3_exp} are satisfied is $a_{52} = 0$~\cite{boscarino:2009_IMEX_relaxation, kennedy:2019}, i.e.
\begin{equation}\label{eq:c3_stage_order2_exp}
    {c_{3} = c_{3}^{(2)}} = \frac{2\rpth{6\gamma^{2} - 6\gamma + 1}}{3\rpth{2\gamma^{2} - 4\gamma + 1}}.
\end{equation}
If such condition is not imposed, then it is possible to fix the stage order 2 only for one of the two intermediate stages, and in such case, we consider $\hat{a}_{32}^{(1)}$.

We are now left with three free parameters, i.e. $\hat{a}_{42}, \hat{a}_{53}$, and $\hat{a}_{54}$. Following~\cite{boscarino:2009_IMEX_relaxation}, we fix $\hat{a}_{53}$ such that
\begin{equation*}
    \sum_{l,m,k}b_{l}\hat{a}_{lm}\hat{a}_{mk}c_{k} = \frac{1}{24},
\end{equation*}
so as to obtain the stability region of a fourth order explicit method. Hence, we get
\begin{equation*}
    \hat{a}_{53} = \frac{\gamma\hat{a}_{54}\rpth{\gamma + b_{3}c_{3}^{2} - \frac{1}{3}} + \frac{1}{12}a_{54}\rpth{24\gamma\hat{a}_{42}a_{54}c_{3} + 2c_{3}\rpth{3\gamma + 3a_{53}c_{3}^{2} - 1} + 1}}{\gamma a_{54}c_{3}^{2}}.
\end{equation*}
The coefficient $\hat{a}_{54}$ is determined by analyzing the \texttt{E}-polynomial associated to the Butcher tableau~\ref{tab:explicit_s5_p3} so as to enhance the order of dissipation. After some tedious calculations, one can show that $E_{4} = 0$ and that
\begin{equation*}
    E_{6} = \frac{\gamma\hat{a}_{54}c_{3}\rpth{144\gamma\hat{a}_{42}a_{54} + 36\gamma  + 36a_{53}c_{3}^{2} - 12} + b_{4}}{72b_{4}}
\end{equation*}
so that, imposing $E_{6} = 0$ in order to achieve dissipation order $7$ (see Definition~\ref{def:RK_dissipation_dispersion}), we obtain
\begin{equation}\label{eq:ahat54_order3_low_dissipation}
    \hat{a}_{54} = \hat{a}_{54}^{(1)} = -\frac{b_{4}}{12\gamma c_{3}\rpth{12\gamma\hat{a}_{42}a_{54} + 3\gamma + 3a_{53}c_{3}^{2} - 1}}.
\end{equation}
If one instead desires dispersion order $6$, which is the maximum attainable~\cite{vanderhouwen:1989}, calculations analogous to those described for the implicit scheme show that
\begin{equation}\label{eq:ahat54_order3_low_dispersion}
    \hat{a}_{54} = \hat{a}_{54}^{(2)} = -\frac{b_{4}}{10\gamma c_{3}\rpth{12\gamma\hat{a}_{42}a_{54} + 3\gamma + 3a_{53}c_{3}^{2} - 1}}.
\end{equation}
Finally, the remaining parameter $\hat{a}_{42}$ may be determined by imposing stage order 3 for the fourth (penultimate) intermediate stage. However, following~\cite{boscarino:2009_IMEX_relaxation}, we consider the choice $\hat{a}_{42} = 0$, which our investigation shows to yield better results and will therefore be adopted in the following. In conclusion, following the nomenclature introduced at the beginning of Section~\ref{sec:optmized_IMEX}, the third-order schemes derived in this work will be denoted by \orderthreeLDs\ when $\hat{a}_{54}^{(1)}$ is employed, highlighting the Low Dissipation (LDs) of the explicit part, and by \orderthreeLDp\ when $\hat{a}_{54}^{(2)}$ is used, highlighting the Low Dispersion (LDp) of the explicit part. A subscript will be used to distinguish between the (sub-)variants. We emphasize that all the schemes are characterized by low dispersion in the implicit part. The coefficients can be found in Appendix~\ref{app:IMEX_coeffs} (Tables~\ref{tab:ACO_order3_LDs_1}-\ref{tab:ACO_order3_LDs_2}).

%%%%%%%%%%%%%% Fourth-order schemes %%%%%%%%%%%%%%%
\subsection{Fourth-order IMEX-RK schemes}
\label{ssec:fourth_order_schemes}

In this section, we consider fourth-order time discretization schemes. The conditions to achieve fourth-order accuracy for a Runge--Kutta method $\rpth{\mathbf{A}, \mathbf{b}, \mathbf{c}}$ read as follows
\begin{subequations}\label{eq:conditions_order_4}
\begin{alignat}{7}
    &\sum_{l=1}^{s}b_{l} = 1, & &  \qquad \sum_{l=1}^{s}b_{l}c_{l} = \frac{1}{2},  & & \qquad \sum_{l=1}^{s}b_{l}c_{l}^{2} = \frac{1}{3}, & & \quad \sum_{l,m=1}^{s}b_{l}a_{lm}c_{m} = \frac{1}{6}, \\ 
    &\sum_{l=1}^{s}b_{l}c_{l}^{3} = \frac{1}{4}, & & \qquad \sum_{l,m=1}^{s}b_{l}c_{l}a_{lm}c_{m} = \frac{1}{8}, & & \qquad \sum_{l,m=1}^{s}b_{l}a_{lm}c_{m}^{2} = \frac{1}{12}, & & \quad \sum_{l,m,k=1}^{s}b_{l}a_{lm}a_{mk}c_{k} = \frac{1}{24}.
\end{alignat}
\end{subequations}

In the general case $\mathbf{b}^{\top} \neq \widehat{\mathbf{b}}^{\top}, \mathbf{c} \neq \widehat{\mathbf{c}}$, $56$ coupling conditions are necessary to obtain a fourth order accurate IMEX-RK method~\cite{pareschi:2005}. The coupling conditions reduce to $21$ if $\mathbf{b}^{\top} = \widehat{\mathbf{b}}^{\top}$ or to $12$ if $\mathbf{c} = \widehat{\mathbf{c}}$~\cite{pareschi:2005}. We consider $\mathbf{b}^{\top} = \widehat{\mathbf{b}}^{\top}$ and $\mathbf{c} = \widehat{\mathbf{c}}$, so that only $2$ coupling conditions are necessary~\cite{ern:2023}, i.e.
\begin{equation}\label{eq:coupling_conditions_order_4}
    \sum_{l,m,k=1}^{s}b_{l}a_{lm}\hat{a}_{mk}c_{k} = \frac{1}{24}, \qquad \sum_{l,m,k=1}^{s}b_{l}\hat{a}_{lm}a_{mk}c_{k} = \frac{1}{24}.
\end{equation}
Finally, the additional order conditions up to fourth order that guarantee better accuracy for an index-1 DAE can be obtained by imposing that the numerical solution and the exact solution agree up to a certain order in the Hilbert expansion with respect to the small parameter $\varepsilon$ of~\eqref{eq:ODE_stiff_partitioned} \cite{boscarino:2009_IMEX_relaxation}. Such order conditions can be also obtained with the help of the so-called bicolour rooted trees~\cite{boscarino:2009_third_order, hairer:1996}. Tedious calculations (see Appendix~\ref{app:index1_DAE_order4}) show that the supplementary order conditions for a fourth-order method are
\begin{subequations}\label{eq:index_1_DAE_conditions_order_4}
\begin{alignat}{5}
    &\sum_{l,m=2}^{s}b_{l}w_{lm}\hat{c}_{m} = 1, & & \qquad \sum_{l,m=2}^{s}b_{l}w_{lm}\hat{c}_{m}^{2} = 1, & & \qquad \sum_{l,m,k=2}^{s}b_{l}w_{lm}\hat{a}_{mk}\hat{c}_{k} = \frac{1}{2}, \label{eq:index_1_DAE_conditions_order_4_part1} \\
    &\sum_{l,m=2}^{s}b_{l}w_{lm}\hat{c}_{m}^{3} = 1, & & \qquad \sum_{l,m,k=2}^{s}b_{l}w_{lm}\hat{a}_{mk}\hat{c}_{k}^{2} = \frac{1}{3}, & &\qquad \sum_{l,m,k,j=2}^{s}b_{l}w_{lm}\hat{a}_{mk}\hat{a}_{kj}\hat{c}_{j} = \frac{1}{6}, \label{eq:index_1_DAE_conditions_order_4_part2} \\ &\sum_{l,m,k=2}^{s}b_{l}w_{lm}\hat{c}_{m}\hat{a}_{mk}\hat{c}_{k} = \frac{1}{2}. \label{eq:index_1_DAE_conditions_order_4_part3}
\end{alignat}    
\end{subequations}

We immediately consider six-stage ($s = 6$) fourth-order time discretization schemes. Indeed, as discussed, e.g., in~\cite{kennedy:2016}, when considering a five-stage fourth-order time discretization scheme, \texttt{L}-stability is guaranteed only for a discrete set of values of $\gamma$. Hence, it is not possible to mimic the analysis performed in Section~\ref{ssec:third_order_schemes} in terms of optimizing the dispersion properties of the scheme. Six-stage stiffly-accurate ESDIRK methods have 16 degrees of freedom (Table~\ref{tab:implicit_s6_p4})

\begin{table}[H]
    \centering
    \begin{tabular}{c|cccccc}
	    $0$ & $0$ & $0$ & $0$ & $0$ & $0$ & $0$ \\
	    $c_{2}$ & $c_{2} - \gamma$ & $\gamma$ & $0$ & $0$ & $0$ & $0$ \\
        $c_{3}$ & $c_{3} - a_{32} - \gamma$ & $a_{32}$ & $\gamma$ & $0$ & $0$ & $0$ \\
        $c_{4}$ & $c_{4} - a_{42} - a_{43} - \gamma$ & $a_{42}$ & $a_{43}$ & $\gamma$ & $0$ & $0$ \\
        $c_{5}$ & $c_{5} - a_{51} - a_{52} - a_{53} - a_{54} - \gamma$ & $a_{51}$ & $a_{52}$ & $a_{53}$ & $\gamma$ & $0$ \\
        $1$ & $1 - a_{62} - a_{63} - a_{64} - a_{65} - \gamma$ & $a_{62}$ & $a_{63}$ & $a_{64}$ & $a_{65}$ & $\gamma$ \\
	    \hline
        & $1 - a_{62} - a_{63} - a_{64} - a_{65} - \gamma$ & $a_{62}$ & $a_{63}$ & $a_{64}$ & $a_{65}$ & $\gamma$
    \end{tabular}
    \caption{Butcher tableau of the implicit six-stage stiffly-accurate ESDIRK method.}
    \label{tab:implicit_s6_p4}
\end{table}

Considering the order conditions~\eqref{eq:conditions_order_4}, stage order equal to 2, results in 7 remaining degrees of freedom which we identify in $a_{53}, a_{54}, b_{5}, c_{3}, c_{4}, c_{5}$, and $\gamma$. Unfortunately, the expressions of the remaining coefficients with respect to these free parameters are too complex to be obtained symbolically. Requiring the scheme is \texttt{L}-stable results allows to compute, e.g., $a_{53}$, and the following restriction on $\gamma$
\begin{equation}
    0.247005962517487 \le \gamma \le 0.676042393226281.
\end{equation}
Notice that, while the upper bound is identical to that reported in~\cite{kennedy:2016, kennedy:2019}, the lower bound is slightly different. Following, the analysis in Section~\ref{ssec:third_order_schemes}, the value of $\gamma$ is determined imposing order of dispersion equal to 6, and results in $\gamma = 0.2780538841136452$. Finally, another free parameter, e.g., $a_{53}$ is determined imposing the algebraic condition~\eqref{eq:supplementary_algebraic_condition}. One can also verify that the supplementary conditions~\eqref{eq:index_1_DAE_conditions_order_4} that do not depend on the explicit companion method are satisfied. Hence, we are left with 4 free parameters, e.g., $b_{5}, c_{3}, c_{4}, c_{5}$. Following~\cite{boscarino:2009_third_order, kennedy:2019}, $b_{5}$ is determined setting $b_{2} = 0$. Several choices are possible for the remaining coefficients. We have analyzed both stage order and internal stability, and identified three schemes that exhibit superior overall accuracy (see Tables~\ref{tab:ACO_order4_LDp_1}-\ref{tab:ACO_order4_LDp_3}). We refer to the end of the section for a description of the properties of such schemes.

Now, we need to build an explicit companion method on top of it. Obviously, we consider $\mathbf{b} = \widehat{\mathbf{b}}$ and $\mathbf{c} = \widehat{\mathbf{c}}$, so that, by imposing the classical row-simplifying conditions~\eqref{eq:simplyfing_nodes_conditions_RK}, we are left with 10 degrees of freedom (Table~\ref{tab:explicit_s6_p4}).

\begin{table}[H]
    \centering
    \begin{tabular}{c|cccccc}
	    $0$ & $0$ & $0$ & $0$ & $0$ & $0$ & $0$ \\
	    $c_{2}$ & $c_{2}$ & $0$ & $0$ & $0$ & $0$ & $0$ \\
        $c_{3}$ & $c_{3} - \hat{a}_{32}$ & $\hat{a}_{32}$ & $0$ & $0$ & $0$ & $0$ \\
        $c_{4}$ & $c_{4} - \hat{a}_{42} - \hat{a}_{43}$ & $\hat{a}_{42}$ & $\hat{a}_{43}$ & $0$ & $0$ & $0$ \\
        $c_{5}$ & $c_{5} - \hat{a}_{52} - \hat{a}_{53} - \hat{a}_{54}$ & $\hat{a}_{52}$ & $\hat{a}_{53}$ & $\hat{a}_{54}$ & $0$ & $0$ \\
        $1$ & $1 - \hat{a}_{62} - \hat{a}_{63} - \hat{a}_{64} - \hat{a}_{65}$ & $\hat{a}_{62}$ & $\hat{a}_{63}$ & $\hat{a}_{64}$ & $\hat{a}_{65}$ & $0$ \\
	    \hline
        & $1 - a_{62} - a_{63} - a_{64} - a_{65} - \gamma$ & $a_{62}$ & $a_{63}$ & $a_{64}$ & $a_{65}$ & $\gamma$
    \end{tabular}
    \caption{Butcher tableau of the explicit six-stage method.}
    \label{tab:explicit_s6_p4}
\end{table}

The (remaining) order conditions~\eqref{eq:conditions_order_4}, the coupling conditions~\eqref{eq:coupling_conditions_order_4}, and the (remaining) supplementary order conditions~\eqref{eq:index_1_DAE_conditions_order_4} yield 10 conditions. However, one can verify that imposing the order conditions and the supplementary order conditions~\eqref{eq:index_1_DAE_conditions_order_4_part1}-\eqref{eq:index_1_DAE_conditions_order_4_part2}, and taking into account the implicit companion method, the condition~\eqref{eq:index_1_DAE_conditions_order_4_part3} is automatically satisfied. Hence, following~\cite{kennedy:2019}, and similarly to what we have done in Section~\ref{ssec:third_order_schemes}, the remaining parameter is fixed by enlarging the stability region of the explicit method and imposing
\begin{equation*}
\sum_{l,m,k,j}b_{l}\hat{a}_{lm}\hat{a}_{mk}\hat{a}_{kj}c_{j} = \frac{1}{135}.
\end{equation*}
Following the notation introduced in Section~\ref{sec:optmized_IMEX}, the fourth-order schemes derived in this work are denoted by \orderfourLDp, where LDp refers to the low-dispersion property of the implicit method. Specifically, we consider \orderfourLDpone\ (Table~\ref{tab:ACO_order4_LDp_1}), which, for the implicit companinion method, is characterized by stage order 3 at the third intermediate stage and internal stability at the fourth and fifth (penultimate) stages. Next, we have \orderfourLDptwo\ (Table~\ref{tab:ACO_order4_LDp_2}), which achieves stage order 3 at the third and fourth intermediate stages and internal stability at the fifth stage. Finally, \orderfourLDpthree\ (Table~\ref{tab:ACO_order4_LDp_3}) is characterized by stage order 3 at the third, fourth, and fifth intermediate stages.

The difficulty of pursuing a general analytical treatment, together with the need to rely on numerical solvers to determine all coefficients, suggests that for higher-order schemes ($\ell > 4$) alternative approaches to time integration may be preferable. In particular, one may consider space–time variational discretizations~\cite{antonietti:2020}, ADER methods, or deferred correction techniques~\cite{dumbser:2008, hanveiga:2021, micalizzi:2025, offner:2025, toro:2009}. However, these approaches are beyond the scope of the present work.

%%%%%%%%%%%%%%%%%%%%%%%%%%%%%%%%%%%%%%%%%%
%%%%%%%%% Space discretization %%%%%%%%%%%
%%%%%%%%%%%%%%%%%%%%%%%%%%%%%%%%%%%%%%%%%%
\section{Polytopal discontinuous Galerkin method for space discretization}
\label{sec:polydg}

In this section, we introduce the space discretization using the polydG method \cite{antonietti:2013,cangiani:2014,cangiani:2017}. First, let us introduce a polytopic mesh partition $\partition$ of the domain $\Omega$ made of disjoint polytopal elements $K$. For each element $K \in \partition$, we denote by $|K|$ its $d$-dimensional measure. Concerning the space mesh regularity, we refer to the assumptions in~\cite{cangiani:2017, cangiani:2014}. We define interface as the intersection of the $(d-1)$-dimensional facets of two neighbouring elements. Namely, the facets are line segments and (planar) triangles decomposition of the polygons, for $d=2$ and $d=3$, respectively. We denote the set of facets as $\faces = \facesinternal \cup \facesboundary$, where $\facesinternal$ is the set of internal facets and $\facesboundary$ the set of exterior ones. Finally, we split $\facesboundary = \facesD \cup \facesN$, where $\facesD$ and $\facesN$ are the boundary faces contained in $\Gamma_{\mathrm{D}}$ and $\Gamma_{\mathrm{N}}$, respectively. We recall that we assume that $\Gamma_{\mathrm{D}} \cap \Gamma_{\mathrm{N}} = \emptyset$. Using the standard notation of Lebesgue and Sobolev spaces for scalar and vector-valued functions, we introduce the discontinuous finite element spaces
$$V_{h}^{\mathrm{DG}} = \{v \in L^{2}(\Omega):\, v|_{K} \in \mathbb{P}_{\ell}(K)\;\forall K \in \partition, \quad \mathrm{and} \quad \Vh = [V_h^\mathrm{DG}]^{d},$$
where $\mathbb{P}_{\ell}(K)$ is the space of polynomials of total degree less than or equal to $\ell$. Then, we introduce the trace operators~\cite{arnold:2001} for a facet $F \in \facesinternal$ shared by the elements $K^{\pm}$. Let $\bm{n}_{K^{\pm}}$ by the unit normal vector on face $F$ pointing exterior to $K^{\pm}$, respectively. Then, for sufficiently regular scalar-valued functions $w$ we define the average operator $\ave{\cdot}$ such that $\ave{w} = \frac{1}{2}\rpth{w|_{K^{+}} + w|_{K^{-}}}$ and the jump operator $\jump{\cdot}$ such that $\jump{w} = w|_{K^{+}}\bm{n}_{K^{+}} + w|_{K^{-}} \bm{n}_{K^{-}}$. Analogously for a vector-valued function $\bm{v}$ the definitions are $ \ave{\bm{v}} = \frac{1}{2}\rpth{\bm{v}|_{K^{+}} + \bm{v}_{K^{-}}}$ and $\jump{\bm{v}} = \bm{v}|_{K^{+}} \cdot \bm{n}_{K^{+}} + \bm{v}|_{K^{-}} \cdot \bm{n}_{K^{-}}$. In the same way, on the face $F \in \facesD$ the average operators are defined as $\ave{w} = w|_{K}$ and $\ave{\bm{v}} = \bm{v}|_{K}$, while the jump operators as $\jump{w} = \rpth{w - g_{\mathrm{D}}}|_{K}\bm{n}_{K}$ and $\jump{\bm{v}} =\rpth{\bm{v} - \bm{g}_{\mathrm{D}}}|_{K} \cdot \bm{n}_{K}$, where $g_{\mathrm{D}}$ and $\bm{g}_{\mathrm{D}}$ are the corresponding Dirichlet boundary conditions. We introduce the following dG-norm
\begin{equation*}\label{eq:dG-norm}
    \|\bm{u}\|_{\mathrm{dG}}^{2} = \sum_{j=1}^{\ncomp} \rpth{\|\sigma_{j}\grad_{h}u_{j}\|_{\mathbf{L}^{2}(\Omega)}^{2} + \sum_{F \in \facesinternal \cup \facesD}\|\sqrt{\eta_{j}} \jump{u_{j}}\|^{2}_{\mathbf{L}^{2}(F)}}, \qquad \forall \bm{u} \in \Vh,
\end{equation*}
where $\grad_{h}$ is the elementwise gradient and $\eta_{j}:\facesinternal \cup \facesD \rightarrow \mathbb{R}_{+}$ is the penalization function defined as:
\begin{equation}\label{eq:penalty}
    \eta_{j} = \eta_{0}
    \begin{cases}
         \dfrac{\sigma_{j}\ell^{2}|F|}{\{|{K^\pm}|\}_\mathrm{H}}, & \mathrm{on}\; F \in \facesinternal, \\[6pt]
         \dfrac{\sigma_{j}\ell^{2}|F|}{|K|}, & \mathrm{on}\; F \in \facesD.
    \end{cases}
\end{equation}
In Equation~\eqref{eq:penalty}, we are considering the harmonic average operator $\{a^{\pm}\}_{\mathrm{H}} = \frac{2a^{+}a^{-}}{a^{+} + a^{-}}$ on $F\in\facesinternal$, where $a^{+}$ and $a^{-}$ are short notations for $a|_{K^{+}}$ and $a|_{K^{-}}$, respectively. Moreover, $\eta_{0}$ is a parameter to be chosen large enough to have stability~\cite{antonietti:2012}. Finally, we can define the following bilinear form $\mathscr{A}: \Vh \times \Vh \rightarrow \mathbb{R}$ as $\mathscr{A}(\bm{u},\bm{v}) = \sum\limits_{j=1}^{\ncomp} \mathscr{A}_{j}(u_{j},v_{j})$, where
\begin{equation*}
\begin{split}
    \mathscr{A}_{j}(u_{j},v_{j}) =& \int_{\Omega} \sigma_{j}\grad_{h}u_{j} \cdot \grad_{h}v_{j}\,\mathrm{d}\bm{x} + \sum_{F \in \facesinternal \cup \facesD}\int_{F} \eta_{j}\jump{u_{j}} \cdot \jump{v_{j}}\,\mathrm{d}\bm{\sigma} \\ 
    -& \sum_{F \in \facesinternal \cup \facesD}\int_{F} \rpth{\ave{\sigma_{j}\grad_{h}u_{j}} \cdot \jump{v_{j}} + \ave{\sigma_{j}\nabla_{h}v_{j}} \cdot \jump{u_{j}}}\,\mathrm{d}\bm{\sigma},\qquad \forall \bm{u},\bm{v} \in \Vh.
\end{split}
\end{equation*}
It is worth noting that we consider a diagonal diffusion tensor solely to simplify the presentation of the method. However, the method can be readily extended to non-diagonal diffusion tensors, provided that they are uniformly positive definite, so that the well-posedness analysis remains valid. The semi-discrete PolyDG symmetric interior penalty formulation reads: for each $t\in(0,T]$ find $\bm{u}_{h}(t) \in \Vh$ such that
\begin{subequations}\label{eq:DGFormulation}
\begin{alignat}{3}
    \int_{\Omega} \dot{\bm{u}}_{h}(t)\bm{v}_{h}\,\mathrm{d}\bm{x} + \mathscr{A}\rpth{\bm{u}_{h}(t),\bm{v}_{h}} &= \int_{\Omega} \mathbf{F}(\bm{u}_{h}(t),t) \cdot \bm{v}_{h}\,\mathrm{d}\bm{x} + \int_{\Gamma_\mathrm{N}} \bm{g}_\mathrm{N} \cdot \bm{v}_{h}\,\mathrm{d}\bm{\sigma}, && \qquad \forall \bm{v}_{h} \in \Vh, \\
    \bm{u}_{h}(0) &= \bm{u}_{h0}, && \qquad \mathrm{in}\,\Omega,
\end{alignat}
\end{subequations}
where $\bm{u}_{h0} \in \Vh$ is a suitable approximation of the initial condition $\bm{u}_{0}$ and we denote the time derivative with $\dot{\bullet} = \partial_{t}\bullet$. For more details about the derivation of the SIP-dG method for the laplacian term, we refer to~\cite{arnold:2001}. Moreover, we refer to~\cite{antonietti:2024, corti:2023}, for the derivation of the method for two particular cases with $\ncomp=1,2$.
\par
To derive the algebraic formulation of the system we introduce a suitable basis for each scalar equation $(\varphi_{jk})_{k=1}^{\dimDG}$ of the space $V_h^\mathrm{DG}$ such that $\dimDG = \mathrm{dim}(V_h^\mathrm{DG})$. Then, we construct a vectorial basis $(\boldsymbol{\Phi}_{k})_{k=1}^{\ncomp\dimDG}$ such that $\boldsymbol{\Phi}_{k} = [\varphi_{k},0,\dots,0]^{\top}$ for $k=1,\dots,\dimDG$, $\boldsymbol{\Phi}_{2k} = [0,\varphi_{k},0,\dots,0]^{\top}$ for $k=\dimDG+1,\dots,2\dimDG$, and so on. Then we rewrite the solution
\begin{equation*}
   \bm{u}_{h}(\bm{x},t) = \sum_{k=0}^{\ncomp\dimDG} U_{k}(t) \boldsymbol{\Phi}_{k}(\bm{x}),
\end{equation*}
where $\boldsymbol{U} \in \mathbb{R}^{\ncomp\dimDG}$ is the corresponding vector of the expansion coefficients written in terms of the chosen basis. Moreover, we define the matrices and the vectors for $i,j=1,...,\ncomp\dimDG$
\begin{alignat*}{3}
    [\mathbf{M}]_{ij} =& \int_{\Omega} \boldsymbol{\Phi}_{j} \cdot \boldsymbol{\Phi}_{i}\,\mathrm{d}\bm{x}, && \qquad \text{(Mass matrix),} \\
    [\mathbf{K}]_{ij} =& \mathscr{A}\rpth{\boldsymbol{\Phi}_{j}, \boldsymbol{\Phi}_{i}}, && \qquad \text{(Stiffness matrix)} \label{eq:stiffness}, \\
    [\boldsymbol{G}(\boldsymbol{U}(t))]_{j} =& \int_{\Omega} \mathbf{F}\rpth{\bm{u}_{h}(t),t} \cdot \boldsymbol{\Phi}_{j}\,\mathrm{d}\bm{x},
    && \qquad \text{(Nonlinear reaction vector),} \\
    [\boldsymbol{G}_\mathrm{E,I}(\boldsymbol{U}(t))]_{j} = & \int_{\Omega} \mathbf{F}_\mathrm{E,I}\rpth{\bm{u}_{h}(t),t} \cdot \boldsymbol{\Phi}_{j}\,\mathrm{d}\bm{x},
    && \qquad \text{(Nonlinear reaction vectors for SI-IMEX).}
\end{alignat*}

%%%%%%%%%%%%%%%%%%%%%%%%%%%%%%%%%%%%%%%
%%%%%%%%%%% Test cases: ODEs %%%%%%%%%%
%%%%%%%%%%%%%%%%%%%%%%%%%%%%%%%%%%%%%%%
\section{Numerical results: application to ODEs}
\label{sec:numerical_results_ODE}

The numerical methods and the IMEX-RK schemes presented in the previous sections are now validated through a series of test cases. While the ultimate goal of this work is the simulation of diffusion-reaction systems \eqref{eq:model}, the development of novel IMEX-RK schemes warrants a preliminary investigation on benchmark ODE problems. The tests of this section allow us to isolate and assess the fundamental properties of the proposed methods, including their accuracy and stability behaviour, and to compare them with existing approaches. More challenging diffusion-reaction applications are reported in Section~\ref{sec:numerical_results_PDE}. First, we consider a logistic ODE (or Verhulst model), which can be interpreted as a Fisher--Kolmogorov equation with vanishing diffusion. In the test case 2, we investigate the classical van der Pol oscillator~\cite{boscarino:2009_IMEX_relaxation, hairer:1996}, a well-known stiff ODE problem.
\par
For the reader’s convenience, Tables~\ref{tab:third_order_recap} and~\ref{tab:fourth_order_recap} summarize the third- and fourth-order IMEX-RK schemes employed in the numerical experiments, including both schemes available in the literature and those developed in the present work. The corresponding coefficients are reported in Appendix~\ref{app:IMEX_coeffs}. For each scheme, we indicate the number of stages and its main structural properties.

\begin{table}[h!]
    \centering
    \scriptsize
    \begin{tabular}{lllp{8.5cm}}
        \toprule
        \textbf{Scheme} & \textbf{Reference} & \textbf{Stages} & \textbf{Main features} \\
        \midrule
        \boscarinoduemilasedici & \cite{boscarino:2016} & 4 & Type I, \texttt{L}-stable, Strong stability preserving \\[6pt]
        \boscarinoduemilaventidue & \cite{boscarino:2022} & 4 & Type I, SA \\[3pt]
        \boscarinoduemilanove & \cite{boscarino:2009_third_order} & 5 & Type II, SA, \texttt{L}-stable, condition \eqref{eq:supplementary_algebraic_condition}, improved algebraic accuracy \\[6pt]
        \carpenterkennedyduemilatre & \cite{kennedy:2003} & 3 & Type II, SA, \texttt{L}-stable \\[3pt]
        \orderthreeLDsone & present work & 5 & Type II, SA, \texttt{L}-stable, condition \eqref{eq:supplementary_algebraic_condition}, improved algebraic accuracy, triplet $[c_3^{(1-)}, c_4^{(1-)}, \hat{a}_{54}^{(1)}]$ \\[3pt]
        \orderthreeLDstwo & present work & 5 & Type II, SA, \texttt{L}-stable, condition \eqref{eq:supplementary_algebraic_condition}, improved algebraic accuracy, triplet $[c_3^{(2)}, c_4^{(1-)}, \hat{a}_{54}^{(1)}]$ \\[3pt]
        \orderthreeLDp & present work & 5 & Type II, SA, \texttt{L}-stable, condition \eqref{eq:supplementary_algebraic_condition}, improved algebraic accuracy, triplet $[c_3^{(1-)}, c_4^{(1-)}, \hat{a}_{54}^{(2)}]$ \\
        \bottomrule
    \end{tabular}
    \caption{Third-order IMEX-RK schemes employed in the numerical simulations.}
    \label{tab:third_order_recap}
\end{table}

\begin{table}[h!]
    \centering
    \scriptsize
    \begin{tabular}{lllp{8.5cm}}
        \toprule
        \textbf{Scheme} & \textbf{Reference} & \textbf{Stages} & \textbf{Main features} \\
        \midrule
        \calvoduemilauno & \cite{calvo:2001} & 5 & Type ARS, SA, \texttt{L}-stable \\[3pt]
        \kennedycarpenterduemilatre & \cite{kennedy:2003} & 6 & Type II, SA, \texttt{L}-stable \\[3pt]
        \kennedycarpenterduemiladiciannove & \cite{kennedy:2019} & 6 & Type II, SA, \texttt{L}-stable, condition \eqref{eq:supplementary_algebraic_condition} \\[3pt]
        \orderfourLDpone & present work & 6 & Type II, SA, \texttt{L}-stable, condition \eqref{eq:supplementary_algebraic_condition}, improved algebraic accuracy, low dispersion implicit method, stage order 3 3rd stage, internal stability 4th-5th stage \\[3pt]
        \orderfourLDptwo & present work & 6 & Type II, SA, \texttt{L}-stable, condition \eqref{eq:supplementary_algebraic_condition}, improved algebraic accuracy, low dispersion implicit method, stage order 3 3rd-4th stage, internal stability 5th stage \\[3pt]
        \orderfourLDpthree & present work & 6 & Type II, SA, \texttt{L}-stable, condition \eqref{eq:supplementary_algebraic_condition}, improved algebraic accuracy, low dispersion implicit method, stage order 3 3rd-4th-5th stage \\
        \bottomrule
    \end{tabular}
    \caption{Fourth-order IMEX-RK schemes employed in the numerical simulations.}
    \label{tab:fourth_order_recap}
\end{table}

%%%%%%%%%%%% Pure reaction %%%%%%%%%%%%%%%
\subsection{Test case 1: Verhulst model}
\label{ssec:reaction}

As a first test case, we consider the Verhulst model. This is a single ODE, with a nonlinear quadratic term of the form
\begin{equation*}
    \odif{\Verhulstcomp}{t} = \alpha\,\Verhulstcomp\rpth{1 - \Verhulstcomp}, 
\end{equation*}
supplemented with a suitable initial condition $\Verhulstcomp(0) = \Verhulstcomp_{0}$. As we will see in Section~\ref{ssec:travelling_wave_FK}, this equation represents the limiting model of the Fisher--Kolmogorov equation in the absence of diffusion. The analytical solution is
\begin{equation*}
    \Verhulstcomp(t) = \frac{\Verhulstcomp_{0}e^{\alpha t}}{1 + \Verhulstcomp_{0}\rpth{e^{\alpha t} - 1}},
\end{equation*}
In the numerical experiments, we set $\alpha = 1$, $T = 1$, and $\Verhulstcomp_{0} = 0.2$. According to the notation set in Section~\ref{ssec:IMEX_RK_model}, the explicit and implicit components are given by $f_{\mathrm{E}} = -\alpha\,\Verhulstcomp^{2}$ and $f_{\mathrm{I}} = \alpha\,\Verhulstcomp$, respectively.
\par
\medskip
In Figure~\ref{fig:tc1_errors3}, we show that the proposed third-order IMEX--RK schemes yield the smallest errors among all the tested methods. In particular, the combination $c_{3}^{(1-)}$ in~\eqref{eq:c3_stage_order3_s5_p3} and $c_{4}^{(1-)}$~\eqref{eq:c4_stage_order3_s5_p3}, which ensures stage order 3 at the third and fourth internal stages of the implicit companion method, yields lower errors. Overall, we identify three schemes that provide superior accuracy: \orderthreeLDsone\ (Table~\ref{tab:ACO_order3_LDs_1}), characterized by the triplet $[c_{3}^{(1-)}, c_{4}^{(1-)}, \hat{a}_{54}^{(1)}]$, \orderthreeLDstwo\ (Table~\ref{tab:ACO_order3_LDs_2}), characterized by the triplet $[c_{3}^{(2)}, c_{4}^{(1-)}, \hat{a}_{54}^{(1)}]$, and \orderthreeLDp\ (Table~\ref{tab:ACO_order3_LDp}), characterized by the triplet $[c_{3}^{(1-)}, c_{4}^{(1-)}, \hat{a}_{54}^{(2)}]$. We recall that $c_{3}^{(2)}$~\eqref{eq:c3_stage_order2_exp} guarantees stage order 2 for both the third and fourth explicit stage, while $\hat{a}_{54}^{(1)}$~\eqref{eq:ahat54_order3_low_dissipation} maximizes the order of dissipation of the explicit scheme, and $\hat{a}_{54}^{(2)}$~\eqref{eq:ahat54_order3_low_dispersion} maximizes the order of dispersion of the explicit scheme. Among these, \orderthreeLDp\ produces the lowest error values. The IMEX-RK(4,3,3) scheme exhibits order reduction due to the stiffness of the problem and its failure to satisfy condition~\eqref{eq:supplementary_algebraic_condition}.

Regarding the fourth-order schemes, Figure~\ref{fig:tc1_errors4} shows that the proposed methods also perform very well. In particular, \orderfourLDptwo\ (Table~\ref{tab:ACO_order4_LDp_2}) attains the smallest errors, which are of the same order as those of the method by Kennedy and Carpenter~\cite{kennedy:2019}, while requiring one stage fewer. This scheme is associated with stage order 3 for the third and fourth intermediate stage, while internal stability is enforced for the fifth (penultimate) internal stage, i.e. $R^{[5]}(\infty) = 0$, where $R^{[l]}$ denotes the stability function of the $l$-th stage. As for the IMEX-RK(4,3,3) scheme, the ARK4(3)6L[2]SA scheme shows an order reduction, because it does not satisfy condition~\eqref{eq:supplementary_algebraic_condition}.

\begin{figure}[h!]
    \begin{subfigure}[b]{0.5\textwidth}
        \resizebox{\textwidth}{!}{\definecolor{myred}{rgb}{0.75000,0.000000,0.7500000}%
\definecolor{myblue}{rgb}{0.00000,0.750000,0.75000}%
\begin{tikzpicture}

\begin{axis}[%
width =4.50in,
height=3.50in,
at={(2.6in,1.099in)},
scale only axis,
xmode=log,
xmin = 3.9063e-04,
xmax= 1.0000e-01,
xminorticks=true,
xlabel = {$\Delta t$ [-]},
ylabel = {$|\Verhulstcomp(T)-\Verhulstcomp^{N_T})|$},
xticklabel={\pgfmathparse{exp(\tick)}\pgfmathprintnumber{\pgfmathresult}},
x tick label style={
/pgf/number format/.cd, fixed, fixed zerofill,
precision=3},
ymode=log,
ymin=1e-14,
ymax=1e-4,
yminorticks=true,
axis background/.style={fill=white},
title style={font=\bfseries},
xmajorgrids,
xminorgrids,
ymajorgrids,
yminorgrids,
legend style={at={(0.34,0.99)},legend cell align=left, draw=white!15!black, font=\footnotesize}
]

\addplot [color=blue, line width=2.0pt, mark=*]
  table[row sep=crcr]{%
   1.0000e-01   3.1179e-06 \\
   5.0000e-02   3.8774e-07 \\
   2.5000e-02   4.8345e-08 \\
   1.2500e-02   6.0356e-09 \\
   6.2500e-03   7.5398e-10 \\
   3.1250e-03   9.4219e-11 \\
   1.5625e-03   1.1776e-11 \\
   7.8125e-04   1.4709e-12 \\
   3.9063e-04   1.8269e-13 \\
};
\addlegendentry{\boscarinoduemilasedici}      
\addplot [color=cyan, line width=2.0pt, mark=*]
  table[row sep=crcr]{%
   1.0000e-01   3.0258e-06 \\
   5.0000e-02   1.3182e-06 \\
   2.5000e-02   3.9569e-07 \\
   1.2500e-02   1.0695e-07 \\
   6.2500e-03   2.7726e-08 \\
   3.1250e-03   7.0543e-09 \\
   1.5625e-03   1.7789e-09 \\
   7.8125e-04   4.4663e-10 \\
   3.9063e-04   1.1189e-10 \\
};
\addlegendentry{\boscarinoduemilaventidue}

\addplot [color=myblue, line width=2.0pt, mark=*]
  table[row sep=crcr]{%
   1.0000e-01   4.5768e-06 \\
   5.0000e-02   5.4855e-07 \\
   2.5000e-02   6.7175e-08 \\
   1.2500e-02   8.3121e-09 \\
   6.2500e-03   7.3951e-09 \\
   3.1250e-03   4.0194e-08 \\
   1.5625e-03   1.0044e-08 \\
   7.8125e-04   2.5104e-09 \\
   3.9063e-04   6.2752e-10 \\
};
\addlegendentry{\boscarinoduemilanove}

\addplot [color=green, line width=2.0pt, mark=*]
  table[row sep=crcr]{%
   1.0000e-01   4.7602e-06 \\
   5.0000e-02   5.8644e-07 \\
   2.5000e-02   7.2793e-08 \\
   1.2500e-02   9.0680e-09 \\
   6.2500e-03   1.1316e-09 \\
   3.1250e-03   1.4133e-10 \\
   1.5625e-03   1.7658e-11 \\
   7.8125e-04   2.2066e-12 \\
   3.9063e-04   2.7650e-13 \\
};
\addlegendentry{\carpenterkennedyduemilatre}

\addplot [color=red, line width=2.0pt, mark=square*]
  table[row sep=crcr]{%
   1.0000e-01   1.0725e-06 \\
   5.0000e-02   1.2444e-07 \\
   2.5000e-02   1.4969e-08 \\
   1.2500e-02   1.8350e-09 \\
   6.2500e-03   2.2714e-10 \\
   3.1250e-03   2.8252e-11 \\
   1.5625e-03   3.5234e-12 \\
   7.8125e-04   4.4020e-13 \\
   3.9063e-04   5.5178e-14 \\
};
\addlegendentry{\nostroordinetreprimo}

\addplot [color=myred, line width=2.0pt, mark=square*]
  table[row sep=crcr]{%
   1.0000e-01   6.9957e-07 \\
   5.0000e-02   7.8877e-08 \\
   2.5000e-02   9.3366e-09 \\
   1.2500e-02   1.1348e-09 \\
   6.2500e-03   1.3984e-10 \\
   3.1250e-03   1.7354e-11 \\
   1.5625e-03   2.1619e-12 \\
   7.8125e-04   2.6934e-13 \\
   3.9063e-04   3.4028e-14 \\
};
\addlegendentry{\nostroordinetresecondo}

\addplot [color=orange, line width=2.0pt, mark=square*]
  table[row sep=crcr]{%
   1.0000e-01   1.0353e-06 \\
   5.0000e-02   1.2900e-07 \\
   2.5000e-02   1.6103e-08 \\
   1.2500e-02   2.0117e-09 \\
   6.2500e-03   2.5139e-10 \\
   3.1250e-03   3.1419e-11 \\
   1.5625e-03   3.9274e-12 \\
   7.8125e-04   4.9177e-13 \\
   3.9063e-04   6.0507e-14 \\
};
\addlegendentry{\nostroordinetreterzo}

\addplot [color=black, dash dot, line width=2.0pt]
  table[row sep=crcr]{%
   1.0000e-01   1.0000e-07 \\
   3.0000e-04   2.7000e-15 \\
};
\node[right, align=left, text=black]
at (axis cs:0.006,1.7e-11) {$\sim\mathcal{O}(\Delta t^3)$};

\end{axis}
\end{tikzpicture}%
}
        \caption{Third-order methods}
        \label{fig:tc1_errors3}
    \end{subfigure}%
    \begin{subfigure}[b]{0.5\textwidth}
        \resizebox{\textwidth}{!}{\definecolor{mygreen}{rgb}{0.00000,0.75000,0.0000}%
\definecolor{mygreen2}{rgb}{0.75000,0.75000,0.0000}%
\definecolor{mycyan}{rgb}{0.00000,0.50000,1.00000}%
\definecolor{mypurple}{rgb}{0.50000,0.00000,0.75000}%

\begin{tikzpicture}

\begin{axis}[%
width =4.50in,
height=3.50in,
at={(2.6in,1.099in)},
scale only axis,
xmode=log,
xmin = 1.5625e-03,
xmax= 1.0000e-01,
xminorticks=true,
xlabel = {$\Delta t$ [-]},
ylabel = {$|\Verhulstcomp(T)-\Verhulstcomp^{N_T})|$},
xticklabel={\pgfmathparse{exp(\tick)}\pgfmathprintnumber{\pgfmathresult}},
x tick label style={
/pgf/number format/.cd, fixed, fixed zerofill,
precision=3},
ymode=log,
ymin=1e-17,
ymax=1e-6,
yminorticks=true,
axis background/.style={fill=white},
title style={font=\bfseries},
xmajorgrids,
xminorgrids,
ymajorgrids,
yminorgrids,
legend style={at={(0.34,0.99)},legend cell align=left, draw=white!15!black, font=\footnotesize}
]

\addplot [color=mycyan, line width=2.0pt, mark=*]
  table[row sep=crcr]{%
   1.0000e-01   7.6154e-08 \\
   5.0000e-02   4.4804e-09 \\
   2.5000e-02   2.7161e-10 \\
   1.2500e-02   1.6716e-11 \\
   6.2500e-03   1.0365e-12 \\
   3.1250e-03   6.4726e-14 \\
   1.5625e-03   3.7192e-15 \\
};
\addlegendentry{\calvoduemilauno}

\addplot [color=mygreen, line width=2.0pt, mark=*]
  table[row sep=crcr]{%
   1.0000e-01   6.5779e-09 \\
   5.0000e-02   3.3229e-10 \\
   2.5000e-02   2.2788e-11 \\
   1.2500e-02   2.3495e-11 \\
   6.2500e-03   1.2552e-11 \\
   3.1250e-03   6.3376e-12 \\
   1.5625e-03   3.1761e-12 \\
};
\addlegendentry{\kennedycarpenterduemilatre}

\addplot [color=mygreen2, line width=2.0pt, mark=*]
  table[row sep=crcr]{%
   1.0000e-01   1.6780e-09 \\
   5.0000e-02   1.1585e-10 \\
   2.5000e-02   7.5883e-12 \\
   1.2500e-02   4.8483e-13 \\ 
   6.2500e-03   3.0698e-14 \\
   3.1250e-03   2.2204e-15 \\
   1.5625e-03   4.4409e-16 \\
};
\addlegendentry{\kennedycarpenterduemiladiciannove}

\addplot [color=magenta, line width=2.0pt, mark=square*]
  table[row sep=crcr]{%
   1.0000e-01   5.1080e-08 \\
   5.0000e-02   3.0387e-09 \\
   2.5000e-02   1.8538e-10 \\
   1.2500e-02   1.1449e-11 \\
   6.2500e-03   7.1138e-13 \\
   3.1250e-03   4.3965e-14 \\
   1.5625e-03   2.3870e-15 \\
};
\addlegendentry{\nostroordinequattroprimo}

\addplot [color=purple, line width=2.0pt, mark=square*]
  table[row sep=crcr]{%
   1.0000e-01   9.6275e-09 \\
   5.0000e-02   3.9993e-10 \\
   2.5000e-02   1.8713e-11 \\
   1.2500e-02   9.7378e-13 \\
   6.2500e-03   5.4567e-14 \\
   3.1250e-03   3.1641e-15 \\
   1.5625e-03   1.1102e-16 \\
};
\addlegendentry{\nostroordinequattrosecondo}

\addplot [color=violet, line width=2.0pt, mark=square*]
  table[row sep=crcr]{%
   1.0000e-01   1.3749e-07 \\
   5.0000e-02   8.5096e-09 \\
   2.5000e-02   5.2943e-10 \\
   1.2500e-02   3.3022e-11 \\
   6.2500e-03   2.0633e-12 \\
   3.1250e-03   1.2818e-13 \\
   1.5625e-03   7.1054e-15 \\
};
\addlegendentry{\nostroordinequattroterzo}

\addplot [color=black, dashed, line width=2.0pt]
  table[row sep=crcr]{
   2.0000e-01   1.6000e-08 \\
   7.8125e-04   3.7253e-18 \\
};
\node[right, align=left, text=black]
at (axis cs:0.015,1.7e-13) {$\sim\mathcal{O}(\Delta t^4)$};

\end{axis}
\end{tikzpicture}%
}
        \caption{Fourth-order methods}
        \label{fig:tc1_errors4}
    \end{subfigure}%
    \caption{Test case 1, Verhulst model: convergence of different IMEX-RK schemes. The error is computed at final time for third-order (a) and fourth-order methods (b).}
    \label{fig:tc1_errors}
\end{figure}
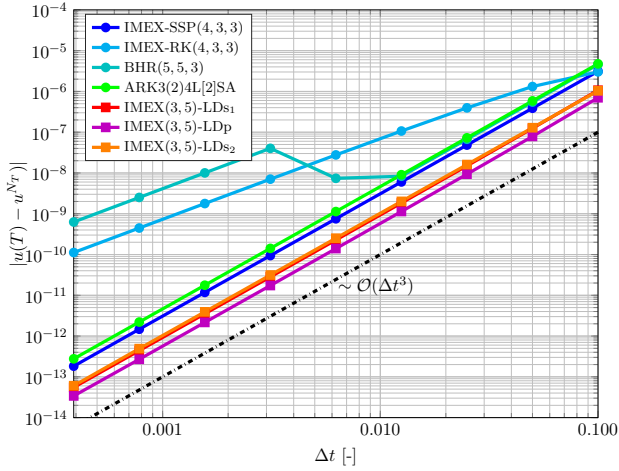
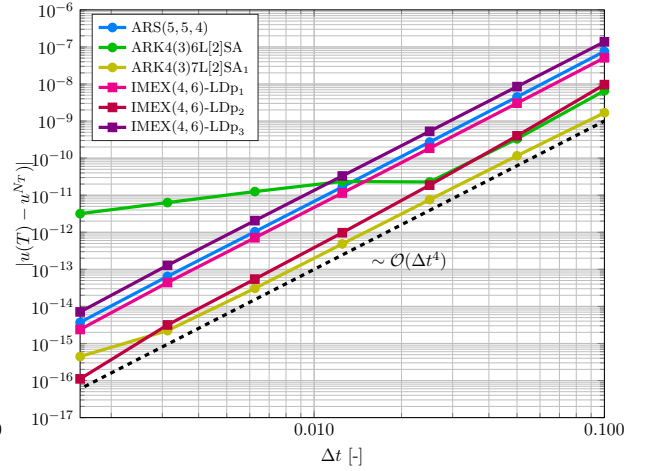

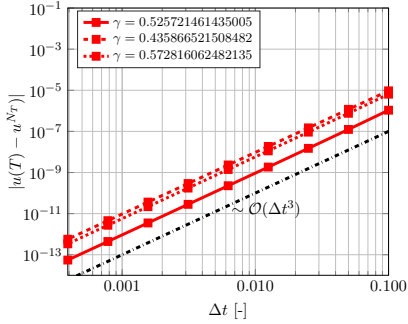
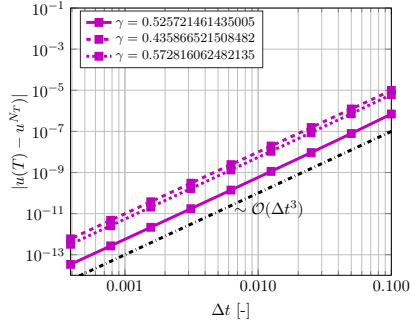
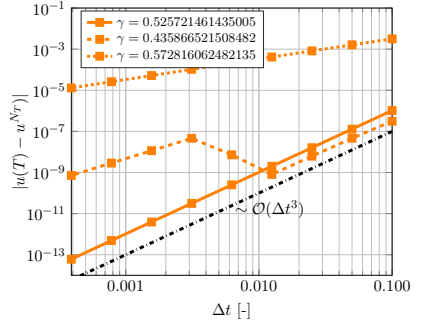
\begin{figure}[h!]
    \begin{subfigure}[b]{0.33\textwidth}
        \resizebox{\textwidth}{!}{\definecolor{myred}{rgb}{0.75000,0.000000,0.7500000}%
\definecolor{myblue}{rgb}{0.00000,0.750000,0.75000}%
\begin{tikzpicture}

\begin{axis}[%
width =3.00in,
height=2.50in,
at={(2.6in,1.099in)},
scale only axis,
xmode=log,
xmin = 3.9063e-04,
xmax= 1.0000e-01,
xminorticks=true,
xlabel = {$\Delta t$ [-]},
ylabel = {$|\Verhulstcomp(T)-\Verhulstcomp^{N_T})|$},
xticklabel={\pgfmathparse{exp(\tick)}\pgfmathprintnumber{\pgfmathresult}},
x tick label style={
/pgf/number format/.cd, fixed, fixed zerofill,
precision=3},
ymode=log,
ymin=1e-14,
ymax=1e-1,
yminorticks=true,
axis background/.style={fill=white},
title style={font=\bfseries},
xmajorgrids,
xminorgrids,
ymajorgrids,
yminorgrids,
legend style={at={(0.60,0.99)},legend cell align=left, draw=white!15!black, font=\footnotesize}
]

\addplot [color=red, line width=2.0pt, mark=square*]
  table[row sep=crcr]{%
   1.0000e-01   1.0725e-06 \\
   5.0000e-02   1.2444e-07 \\
   2.5000e-02   1.4969e-08 \\
   1.2500e-02   1.8350e-09 \\
   6.2500e-03   2.2714e-10 \\
   3.1250e-03   2.8252e-11 \\
   1.5625e-03   3.5234e-12 \\
   7.8125e-04   4.4020e-13 \\
   3.9063e-04   5.5178e-14 \\
};
\addlegendentry{$\gamma=0.525721461435005$}

\addplot [color=red, dashed, line width=2.0pt, mark=square*, mark options=solid]
  table[row sep=crcr]{%
   1.0000e-01   9.4025e-06 \\
   5.0000e-02   1.1823e-06 \\
   2.5000e-02   1.4834e-07 \\
   1.2500e-02   1.8580e-08 \\
   6.2500e-03   2.3249e-09 \\
   3.1250e-03   2.9076e-10 \\ 
   1.5625e-03   3.6355e-11 \\
   7.8125e-04   4.5457e-12 \\
   3.9063e-04   5.6821e-13 \\
};
\addlegendentry{$\gamma=0.435866521508482$}

\addplot [color=red, dotted, line width=2.0pt, mark=square*, mark options=solid]
  table[row sep=crcr]{%
   1.0000e-01   6.4587e-06 \\
   5.0000e-02   7.5806e-07 \\
   2.5000e-02   9.1851e-08 \\
   1.2500e-02   1.1305e-08 \\
   6.2500e-03   1.4022e-09 \\
   3.1250e-03   1.7460e-10 \\
   1.5625e-03   2.1782e-11 \\
   7.8125e-04   2.7207e-12 \\
   3.9063e-04   3.3989e-13 \\
};
\addlegendentry{$\gamma=0.572816062482135$}

\addplot [color=black, dash dot, line width=2.0pt]
  table[row sep=crcr]{%
   1.0000e-01   1.0000e-07 \\
   3.0000e-04   2.7000e-15 \\
};
\node[right, align=left, text=black]
at (axis cs:0.006,1.7e-11) {$\sim\mathcal{O}(\Delta t^3)$};

\end{axis}
\end{tikzpicture}%
}
        \caption{\orderthreeLDsone}
        \label{fig:tc1_errorsgamma1}
    \end{subfigure}%    
    \begin{subfigure}[b]{0.33\textwidth}
        \resizebox{\textwidth}{!}{\definecolor{myred}{rgb}{0.75000,0.000000,0.7500000}%
\definecolor{myblue}{rgb}{0.00000,0.750000,0.75000}%
\begin{tikzpicture}

\begin{axis}[%
width =3.00in,
height=2.50in,
at={(2.6in,1.099in)},
scale only axis,
xmode=log,
xmin = 3.9063e-04,
xmax= 1.0000e-01,
xminorticks=true,
xlabel = {$\Delta t$ [-]},
ylabel = {$|\Verhulstcomp(T)-\Verhulstcomp^{N_T})|$},
xticklabel={\pgfmathparse{exp(\tick)}\pgfmathprintnumber{\pgfmathresult}},
x tick label style={
/pgf/number format/.cd, fixed, fixed zerofill,
precision=3},
ymode=log,
ymin=1e-14,
ymax=1e-1,
yminorticks=true,
axis background/.style={fill=white},
title style={font=\bfseries},
xmajorgrids,
xminorgrids,
ymajorgrids,
yminorgrids,
legend style={at={(0.60,0.99)},legend cell align=left, draw=white!15!black, font=\footnotesize}
]

\addplot [color=myred, line width=2.0pt, mark=square*]
  table[row sep=crcr]{%
   1.0000e-01   6.9957e-07 \\
   5.0000e-02   7.8877e-08 \\
   2.5000e-02   9.3366e-09 \\
   1.2500e-02   1.1348e-09 \\
   6.2500e-03   1.3984e-10 \\
   3.1250e-03   1.7354e-11 \\
   1.5625e-03   2.1619e-12 \\
   7.8125e-04   2.6934e-13 \\
   3.9063e-04   3.4028e-14 \\
};
\addlegendentry{$\gamma=0.525721461435005$}

\addplot [color=myred, dashed, line width=2.0pt, mark=square*, mark options=solid]
  table[row sep=crcr]{%
   1.0000e-01   9.7959e-06 \\
   5.0000e-02   1.2308e-06 \\
   2.5000e-02   1.5436e-07 \\
   1.2500e-02   1.9330e-08 \\
   6.2500e-03   2.4186e-09 \\
   3.1250e-03   3.0247e-10 \\
   1.5625e-03   3.7819e-11 \\
   7.8125e-04   4.7277e-12 \\
   3.9063e-04   5.9253e-13 \\
};
\addlegendentry{$\gamma=0.435866521508482$}

\addplot [color=myred, dotted, line width=2.0pt, mark=square*, mark options=solid]
  table[row sep=crcr]{%
   1.0000e-01   6.1565e-06 \\
   5.0000e-02   7.2108e-07 \\
   2.5000e-02   8.7277e-08 \\
   1.2500e-02   1.0736e-08 \\
   6.2500e-03   1.3313e-09 \\
   3.1250e-03   1.6575e-10 \\
   1.5625e-03   2.0677e-11 \\
   7.8125e-04   2.5819e-12 \\
   3.9063e-04   3.2363e-13 \\
};
\addlegendentry{$\gamma=0.572816062482135$}

\addplot [color=black, dash dot, line width=2.0pt]
  table[row sep=crcr]{%
   1.0000e-01   1.0000e-07 \\
   3.0000e-04   2.7000e-15 \\
};
\node[right, align=left, text=black]
at (axis cs:0.006,1.7e-11) {$\sim\mathcal{O}(\Delta t^3)$};

\end{axis}
\end{tikzpicture}%
}
        \caption{\orderthreeLDp}
        \label{fig:tc1_errorsgamma2}
    \end{subfigure}%
    \begin{subfigure}[b]{0.33\textwidth}
        \resizebox{\textwidth}{!}{\definecolor{myred}{rgb}{0.75000,0.000000,0.7500000}%
\definecolor{myblue}{rgb}{0.00000,0.750000,0.75000}%
\begin{tikzpicture}

\begin{axis}[%
width =3.00in,
height=2.50in,
at={(2.6in,1.099in)},
scale only axis,
xmode=log,
xmin = 3.9063e-04,
xmax= 1.0000e-01,
xminorticks=true,
xlabel = {$\Delta t$ [-]},
ylabel = {$|\Verhulstcomp(T)-\Verhulstcomp^{N_T})|$},
xticklabel={\pgfmathparse{exp(\tick)}\pgfmathprintnumber{\pgfmathresult}},
x tick label style={
/pgf/number format/.cd, fixed, fixed zerofill,
precision=3},
ymode=log,
ymin=1e-14,
ymax=1e-1,
yminorticks=true,
axis background/.style={fill=white},
title style={font=\bfseries},
xmajorgrids,
xminorgrids,
ymajorgrids,
yminorgrids,
legend style={at={(0.60,0.99)},legend cell align=left, draw=white!15!black, font=\footnotesize}
]

\addplot [color=orange, line width=2.0pt, mark=square*]
  table[row sep=crcr]{%
   1.0000e-01   1.0353e-06 \\
   5.0000e-02   1.2900e-07 \\
   2.5000e-02   1.6103e-08 \\
   1.2500e-02   2.0117e-09 \\
   6.2500e-03   2.5139e-10 \\
   3.1250e-03   3.1419e-11 \\
   1.5625e-03   3.9274e-12 \\
   7.8125e-04   4.9177e-13 \\
   3.9063e-04   6.0507e-14 \\
};
\addlegendentry{$\gamma=0.525721461435005$}

\addplot [color=orange, dashed, line width=2.0pt, mark=square*, mark options=solid]
  table[row sep=crcr]{%
   1.0000e-01   3.1299e-07 \\
   5.0000e-02   4.6336e-08 \\
   2.5000e-02   6.2208e-09 \\
   1.2500e-02   8.0377e-10 \\
   6.2500e-03   7.2301e-09 \\
   3.1250e-03   4.6128e-08 \\
   1.5625e-03   1.1536e-08 \\
   7.8125e-04   2.8845e-09 \\
   3.9063e-04   7.2120e-10 \\
};
\addlegendentry{$\gamma=0.435866521508482$}

\addplot [color=orange, dotted, line width=2.0pt, mark=square*, mark options=solid]
  table[row sep=crcr]{%
   1.0000e-01   3.1706e-03 \\
   5.0000e-02   1.6268e-03 \\
   2.5000e-02   8.2362e-04 \\
   1.2500e-02   4.1435e-04 \\
   6.2500e-03   2.0781e-04 \\
   3.1250e-03   1.0406e-04 \\
   1.5625e-03   5.2070e-05 \\
   7.8125e-04   2.6045e-05 \\
   3.9063e-04   1.3025e-05 \\
};
\addlegendentry{$\gamma=0.572816062482135$}

\addplot [color=black, dash dot, line width=2.0pt]
  table[row sep=crcr]{%
   1.0000e-01   1.0000e-07 \\
   3.0000e-04   2.7000e-15 \\
};
\node[right, align=left, text=black]
at (axis cs:0.006,1.7e-11) {$\sim\mathcal{O}(\Delta t^3)$};

\end{axis}
\end{tikzpicture}%
}
        \caption{\orderthreeLDstwo}
        \label{fig:tc1_errorsgamma3}
    \end{subfigure}%
    \caption{Test case 1, Verhulst model: comparison of three different versions of our method with $\gamma \simeq 0.5257$, derived in this work, $\gamma \simeq 0.4358$, classical literature choice, and $\gamma \simeq 0.5728$, obtained by maximizing the order of dissipation of the implicit method.}
    \label{fig:tc1_errors_gamma}
\end{figure}

Additionally, we compare three versions of our method using three values of $\gamma$: $\gamma \simeq 0.5257$, derived here, the classical $\gamma \simeq 0.4358$~\cite{boscarino:2009_third_order, kennedy:2003}, and $\gamma \simeq 0.5728$, obtained by maximizing the order of dissipation (see Section~\ref{ssec:third_order_schemes}). Figure~\ref{fig:tc1_errors_gamma} shows that $\gamma \simeq 0.5257$ reduces errors by one order of magnitude for \orderthreeLDsone\ and \orderthreeLDp. For the \orderthreeLDstwo\ scheme, however, $\gamma \simeq 0.4358$ yields smaller errors at larger $\Delta t$, though they rise after $10^{-10}$, similarly to what happens for the BHR(5,5,3) scheme, see Figure~\ref{fig:tc1_errors}; this issue does not occur with $\gamma \simeq 0.5257$. These findings indicate that the choice of $\gamma$ may strongly affect the accuracy of the method, and that the selected value of $\gamma$ improves the performance of the proposed third-order methods.

%%%%%%%%%%%% van der Pol %%%%%%%%%%%%%%
\subsection{Van der Pol oscillator}
\label{ssec:van_der_pol}

The van der Pol’s equation~\cite{boscarino:2009_third_order, hairer:1996} is one of the most widely studied nonlinear equations in the stiff ODE literature. It can be recast as the following first-order ODE system
\begin{alignat*}{2}
    \odif{u_{1}}{t} &= u_{2} \\
    \odif{u_{2}}{t} &= \frac{1}{\varepsilon}\rpth{\rpth{1 - u_{1}^{2}}u_{2} - u_{1}},
\end{alignat*}
with $0 < \varepsilon \ll 1$. The IMEX-RK methods considered here treat the first equation explicitly and the second implicitly. This test problem is chosen to compare the accuracy of several types of IMEX-RK schemes and the new schemes in the presence of severe stiffness. In particular, the numerical experiments are performed with $\varepsilon = 10^{-6}$. Finally, following~\cite{boscarino:2009_third_order}, we consider $T = \SI{0.55139}{\second}$ and prescribe the initial conditions $u_{1}(0) = 2, u_{2}(0) = -\frac{2}{3} + \frac{10}{81}\varepsilon - \frac{292}{2187}\varepsilon^{2} - \frac{1814}{19683}\varepsilon^{3} + \mathcal{O}(\varepsilon^{4})$.

\begin{figure}[H]
    \begin{subfigure}[b]{0.5\textwidth}
        \resizebox{\textwidth}{!}{\definecolor{myred}{rgb}{0.75000,0.000000,0.7500000}%
\definecolor{myblue}{rgb}{0.00000,0.750000,0.75000}%
\begin{tikzpicture}

\begin{axis}[%
width =4.50in,
height=2.50in,
at={(2.6in,1.099in)},
scale only axis,
xmode=log,
xmin = 2.1539e-04,
xmax = 5.5139e-02,
xminorticks=true,
xlabel = {$\Delta t$ [-]},
ylabel = {$|u_1(T)-u_1^{N_T})|$},
xticklabel={\pgfmathparse{exp(\tick)}\pgfmathprintnumber{\pgfmathresult}},
x tick label style={
/pgf/number format/.cd, fixed, fixed zerofill,
precision=3},
ymode=log,
ymin=1e-14,
ymax=1e-1,
yminorticks=true,
axis background/.style={fill=white},
title style={font=\bfseries},
xmajorgrids,
xminorgrids,
ymajorgrids,
yminorgrids,
legend style={at={(0.34,0.99)},legend cell align=left, draw=white!15!black, font=\footnotesize}
]

\addplot [color=blue, line width=2.0pt, mark=*]
  table[row sep=crcr]{%
   5.5139e-02   5.0222e-06 \\
   2.7570e-02   6.4902e-07 \\
   1.3785e-02   9.1132e-08 \\
   6.8924e-03   1.6563e-08 \\
   3.4462e-03   4.6994e-09 \\
   1.7231e-03   1.9112e-09 \\
   8.6155e-04   9.0201e-10 \\
   4.3077e-04   4.4355e-10 \\
   2.1539e-04   2.1955e-10 \\
};
\addlegendentry{\boscarinoduemilasedici}
          
\addplot [color=cyan, line width=2.0pt, mark=*]
  table[row sep=crcr]{%
   5.5139e-02   6.7398e-05 \\
   2.7570e-02   1.5966e-05 \\
   1.3785e-02   3.8683e-06 \\
   6.8924e-03   9.5084e-07 \\
   3.4462e-03   2.3561e-07 \\
   1.7231e-03   5.8619e-08 \\
   8.6155e-04   1.4611e-08 \\
   4.3077e-04   3.6430e-09 \\   
   2.1539e-04   9.0719e-10 \\
};
\addlegendentry{\boscarinoduemilaventidue}

\addplot [color=myblue, line width=2.0pt, mark=*]
  table[row sep=crcr]{%
   5.5139e-02   9.8162e-05 \\
   2.7570e-02   1.1008e-05 \\
   1.3785e-02   1.3080e-06 \\
   6.8924e-03   1.5954e-07 \\
   3.4462e-03   1.9703e-08 \\
   1.7231e-03   2.4487e-09 \\
   8.6155e-04   3.0566e-10 \\
   4.3077e-04   3.8630e-11 \\
   2.1539e-04   5.3098e-12 \\
};
\addlegendentry{\boscarinoduemilanove}

\addplot [color=green, line width=2.0pt, mark=*]
  table[row sep=crcr]{%
   5.5139e-02   1.6417e-06 \\
   2.7570e-02   2.1930e-07 \\
   1.3785e-02   2.8039e-08 \\
   6.8924e-03   3.5436e-09 \\
   3.4462e-03   4.4744e-10 \\
   1.7231e-03   5.6473e-11 \\
   8.6155e-04   6.8205e-12 \\
   4.3077e-04   4.2721e-13 \\
   2.1539e-04   4.1345e-13 \\
};
\addlegendentry{\carpenterkennedyduemilatre}

\addplot [color=red, line width=2.0pt, mark=square*]
  table[row sep=crcr]{%
   5.5139e-02   2.1719e-05 \\
   2.7570e-02   2.6595e-06 \\
   1.3785e-02   3.2862e-07 \\
   6.8924e-03   4.0829e-08 \\
   3.4462e-03   5.0882e-09 \\
   1.7231e-03   6.3553e-10 \\
   8.6155e-04   7.9871e-11 \\
   4.3077e-04   1.0465e-11 \\
   2.1539e-04   1.7883e-12 \\
};
\addlegendentry{\nostroordinetreprimo}

\addplot [color=myred, line width=2.0pt, mark=square*]
  table[row sep=crcr]{%
   5.5139e-02   2.0632e-05 \\
   2.7570e-02   2.5279e-06 \\
   1.3785e-02   3.1242e-07 \\
   6.8924e-03   3.8818e-08 \\
   3.4462e-03   4.8377e-09 \\
   1.7231e-03   6.0426e-10 \\
   8.6155e-04   7.5962e-11 \\
   4.3077e-04   9.9780e-12 \\
   2.1539e-04   1.7271e-12 \\
};
\addlegendentry{\nostroordinetresecondo}

\addplot [color=orange, line width=2.0pt, mark=square*]
  table[row sep=crcr]{%
   5.5139e-02   7.7334e-06 \\
   2.7570e-02   9.5948e-07 \\
   1.3785e-02   1.1930e-07 \\
   6.8924e-03   1.4867e-08 \\
   3.4462e-03   1.8558e-09 \\
   1.7231e-03   2.3225e-10 \\
   8.6155e-04   2.9502e-11 \\
   4.3077e-04   4.1698e-12 \\
   2.1539e-04   1.0056e-12 \\
};
\addlegendentry{\nostroordinetreterzo}

\addplot [color=black, dash dot, line width=2.0pt]
  table[row sep=crcr]{%
   1.0000e-01   2.0000e-06 \\
   1.0000e-04   2.0000e-15 \\
};
\node[right, align=left, text=black]
at (axis cs:0.006,1.7e-10) {$\sim\mathcal{O}(\Delta t^3)$};

\end{axis}
\end{tikzpicture}%
}
        \caption{Third-order methods: differential variable}
        \label{fig:tc2_errors3_diff}
    \end{subfigure}%
    \begin{subfigure}[b]{0.5\textwidth}
        \resizebox{\textwidth}{!}{\definecolor{mygreen}{rgb}{0.00000,0.75000,0.0000}%
\definecolor{mygreen2}{rgb}{0.75000,0.75000,0.0000}%
\definecolor{mycyan}{rgb}{0.00000,0.50000,1.00000}%
\definecolor{mypurple}{rgb}{0.50000,0.00000,0.75000}%

\begin{tikzpicture}

\begin{axis}[%
width =4.50in,
height=2.50in,
at={(2.6in,1.099in)},
scale only axis,
xmode=log,
xmin = 1.0769e-03,
xmax = 2.7570e-01,
xminorticks=true,
xlabel = {$\Delta t$ [-]},
ylabel = {$|u_1(T)-u_1^{N_T})|$},
xticklabel={\pgfmathparse{exp(\tick)}\pgfmathprintnumber{\pgfmathresult}},
x tick label style={
/pgf/number format/.cd, fixed, fixed zerofill,
precision=3},
ymode=log,
ymin=1e-14,
ymax=1e-1,
yminorticks=true,
axis background/.style={fill=white},
title style={font=\bfseries},
xmajorgrids,
xminorgrids,
ymajorgrids,
yminorgrids,
legend style={at={(0.34,0.99)},legend cell align=left, draw=white!15!black, font=\footnotesize}
]

\addplot [color=mycyan, line width=2.0pt, mark=*]
  table[row sep=crcr]{%
   2.7570e-01   1.0671e-05 \\
   1.3785e-01   4.3493e-07 \\
   6.8924e-02   1.0012e-07 \\
   3.4462e-02   1.0298e-08 \\
   1.7231e-02   9.6831e-10 \\
   8.6155e-03   1.1500e-10 \\
   4.3077e-03   1.9573e-11 \\
   2.1539e-03   3.9078e-12 \\
   1.0769e-03   5.3046e-13 \\ 
};
\addlegendentry{\calvoduemilauno}

\addplot [color=mygreen, line width=2.0pt, mark=*]
  table[row sep=crcr]{%
   2.7570e-01   1.8759e-04 \\
   1.3785e-01   1.5448e-05 \\
   6.8924e-02   1.0700e-06 \\
   3.4462e-02   6.9857e-08 \\
   1.7231e-02   4.8575e-09 \\
   8.6155e-03   5.3276e-10 \\
   4.3077e-03   1.4745e-10 \\
   2.1539e-03   6.6045e-11 \\
   1.0769e-03   3.2285e-11 \\
};
\addlegendentry{\kennedycarpenterduemilatre}

\addplot [color=mygreen2, line width=2.0pt, mark=*]
  table[row sep=crcr]{%
   2.7570e-01   5.4899e-05 \\
   1.3785e-01   3.4273e-06 \\
   6.8924e-02   1.4946e-07 \\
   3.4462e-02   4.4940e-09 \\
   1.7231e-02   2.9767e-11 \\
   8.6155e-03   1.4358e-11 \\
   4.3077e-03   3.5356e-12 \\
   2.1539e-03   1.2097e-12 \\
   1.0769e-03   7.1165e-13 \\
};
\addlegendentry{\kennedycarpenterduemiladiciannove}

\addplot [color=magenta, line width=2.0pt, mark=square*]
  table[row sep=crcr]{%
   2.7570e-01   7.1659e-05 \\
   1.3785e-01   1.9112e-06 \\
   6.8924e-02   1.3226e-08 \\
   3.4462e-02   5.7920e-09 \\
   1.7231e-02   5.2317e-10 \\
   8.6155e-03   3.7139e-11 \\
   4.3077e-03   2.6550e-12 \\
   2.1539e-03   6.0330e-13 \\
   1.0769e-03   5.3135e-13 \\
};
\addlegendentry{\nostroordinequattroprimo}

\addplot [color=purple, line width=2.0pt, mark=square*]
  table[row sep=crcr]{%
   2.7570e-01   7.9054e-05 \\
   1.3785e-01   7.3548e-06 \\
   6.8924e-02   5.4782e-07 \\
   3.4462e-02   3.6880e-08 \\
   1.7231e-02   2.3726e-09 \\
   8.6155e-03   1.4731e-10 \\
   4.3077e-03   8.0520e-12 \\
   2.1539e-03   1.8607e-13 \\
   1.0769e-03   5.7088e-13 \\
};
\addlegendentry{\nostroordinequattrosecondo}

\addplot [color=violet, line width=2.0pt, mark=square*]
  table[row sep=crcr]{%
   2.7570e-01   5.5013e-04 \\
   1.3785e-01   5.7218e-05 \\
   6.8924e-02   4.5847e-06 \\
   3.4462e-02   3.2071e-07 \\
   1.7231e-02   2.1107e-08 \\
   8.6155e-03   1.3549e-09 \\
   4.3077e-03   8.6298e-11 \\
   2.1539e-03   5.2003e-12 \\
   1.0769e-03   1.1746e-13 \\
};
\addlegendentry{\nostroordinequattroterzo}

\addplot [color=black, dashed, line width=2.0pt]
  table[row sep=crcr]{
   1.0000e-00   1.0000e-04 \\
   1.0000e-03   1.0000e-16 \\
};
\node[right, align=left, text=black]
at (axis cs:0.015,1.7e-12) {$\sim\mathcal{O}(\Delta t^4)$};

\end{axis}
\end{tikzpicture}%
}
        \caption{Fourth-order methods: differential variable}
        \label{fig:tc2_errors4_diff}
    \end{subfigure} \\
    \begin{subfigure}[b]{0.5\textwidth}
        \resizebox{\textwidth}{!}{\definecolor{myred}{rgb}{0.75000,0.000000,0.7500000}%
\definecolor{myblue}{rgb}{0.00000,0.750000,0.75000}%
\begin{tikzpicture}

\begin{axis}[%
width =4.50in,
height=2.50in,
at={(2.6in,1.099in)},
scale only axis,
xmode=log,
xmin = 2.1539e-04,
xmax = 5.5139e-02,
xminorticks=true,
xlabel = {$\Delta t$ [-]},
ylabel = {$|u_1(T)-u_1^{N_T})|$},
xticklabel={\pgfmathparse{exp(\tick)}\pgfmathprintnumber{\pgfmathresult}},
x tick label style={
/pgf/number format/.cd, fixed, fixed zerofill,
precision=3},
ymode=log,
ymin=1e-13,
ymax=1e-1,
yminorticks=true,
axis background/.style={fill=white},
title style={font=\bfseries},
xmajorgrids,
xminorgrids,
ymajorgrids,
yminorgrids,
legend style={at={(0.74,0.99)},legend cell align=left, draw=white!15!black, font=\footnotesize}
]

\addplot [color=blue, line width=2.0pt, mark=*]
  table[row sep=crcr]{%
   5.5139e-02   2.0948e-02 \\
   2.7570e-02   1.2780e-02 \\
   1.3785e-02   7.0662e-03 \\
   6.8924e-03   3.7173e-03 \\
   3.4462e-03   1.9072e-03 \\
   1.7231e-03   9.6642e-04 \\
   8.6155e-04   4.8686e-04 \\
   4.3077e-04   2.4472e-04 \\
   2.1539e-04   1.2298e-04 \\
};
%\addlegendentry{\boscarinoduemilasedici}
          
\addplot [color=cyan, line width=2.0pt, mark=*]
  table[row sep=crcr]{%
   5.5139e-02   9.1226e-04 \\
   2.7570e-02   3.0225e-04 \\
   1.3785e-02   8.6817e-05 \\
   6.8924e-03   2.3248e-05 \\
   3.4462e-03   6.0072e-06 \\
   1.7231e-03   1.5228e-06 \\
   8.6155e-04   3.8130e-07 \\
   4.3077e-04   9.4392e-08 \\
   2.1539e-04   2.2988e-08 \\
};
%\addlegendentry{\boscarinoduemilaventidue}

\addplot [color=myblue, line width=2.0pt, mark=*]
  table[row sep=crcr]{%
   5.5139e-02   1.6128e-04 \\
   2.7570e-02   1.8164e-05 \\
   1.3785e-02   2.1555e-06 \\
   6.8924e-03   2.6263e-07 \\
   3.4462e-03   3.2457e-08 \\
   1.7231e-03   4.0471e-09 \\
   8.6155e-04   5.0903e-10 \\
   4.3077e-04   6.5308e-11 \\
   2.1539e-04   9.2957e-12 \\
};
%\addlegendentry{\boscarinoduemilanove}

\addplot [color=green, line width=2.0pt, mark=*]
  table[row sep=crcr]{%
   5.5139e-02   2.4054e-03 \\
   2.7570e-02   6.4672e-04 \\
   1.3785e-02   1.6787e-04 \\
   6.8924e-03   4.2770e-05 \\
   3.4462e-03   1.0790e-05 \\
   1.7231e-03   2.7074e-06 \\
   8.6155e-04   6.7685e-07 \\
   4.3077e-04   1.6860e-07 \\
   2.1539e-04   4.1772e-08 \\
};
%\addlegendentry{\carpenterkennedyduemilatre}

\addplot [color=red, line width=2.0pt, mark=square*]
  table[row sep=crcr]{%
   5.5139e-02   4.3035e-05 \\
   2.7570e-02   5.5822e-06 \\
   1.3785e-02   7.1319e-07 \\
   6.8924e-03   9.1200e-08 \\
   3.4462e-03   1.2055e-08 \\
   1.7231e-03   1.8081e-09 \\
   8.6155e-04   3.7395e-10 \\
   4.3077e-04   1.1987e-10 \\
   2.1539e-04   5.0792e-11 \\
};
%\addlegendentry{\nostroordinetreprimo}

\addplot [color=myred, line width=2.0pt, mark=square*]
  table[row sep=crcr]{%
   5.5139e-02   5.8422e-05 \\
   2.7570e-02   7.7089e-06 \\
   1.3785e-02   9.9390e-07 \\
   6.8924e-03   1.2729e-07 \\
   3.4462e-03   1.6630e-08 \\
   1.7231e-03   2.3839e-09 \\
   8.6155e-04   4.4618e-10 \\
   4.3077e-04   1.2890e-10 \\
   2.1539e-04   5.1855e-11 \\
};
%\addlegendentry{\nostroordinetresecondo}

\addplot [color=orange, line width=2.0pt, mark=square*]
  table[row sep=crcr]{%
   5.5139e-02   2.5484e-05 \\
   2.7570e-02   3.2706e-06 \\
   1.3785e-02   4.1384e-07 \\
   6.8924e-03   5.2117e-08 \\
   3.4462e-03   6.5633e-09 \\
   1.7231e-03   8.2971e-10 \\
   8.6155e-04   1.0571e-10 \\
   4.3077e-04   1.3636e-11 \\
   2.1539e-04   1.6389e-12 \\
};
%\addlegendentry{\nostroordinetreterzo}

\addplot [color=black, dash dot, line width=2.0pt]
  table[row sep=crcr]{%
   1.0000e-01   2.0000e-06 \\
   1.0000e-04   2.0000e-15 \\
};
\node[right, align=left, text=black]
at (axis cs:0.006,1.7e-10) {$\sim\mathcal{O}(\Delta t^3)$};

\end{axis}
\end{tikzpicture}%
}
        \caption{Third-order methods: algebraic variable}
        \label{fig:tc2_errors3_alg}
    \end{subfigure}%
    \begin{subfigure}[b]{0.5\textwidth}
        \resizebox{\textwidth}{!}{\definecolor{mygreen}{rgb}{0.00000,0.75000,0.0000}%
\definecolor{mygreen2}{rgb}{0.75000,0.75000,0.0000}%
\definecolor{mycyan}{rgb}{0.00000,0.50000,1.00000}%
\definecolor{mypurple}{rgb}{0.50000,0.00000,0.75000}%

\begin{tikzpicture}

\begin{axis}[%
width =4.50in,
height=2.50in,
at={(2.6in,1.099in)},
scale only axis,
xmode=log,
xmin = 1.0769e-03,
xmax = 2.7570e-01,
xminorticks=true,
xlabel = {$\Delta t$ [-]},
ylabel = {$|u_1(T)-u_1^{N_T})|$},
xticklabel={\pgfmathparse{exp(\tick)}\pgfmathprintnumber{\pgfmathresult}},
x tick label style={
/pgf/number format/.cd, fixed, fixed zerofill,
precision=3},
ymode=log,
ymin=1e-14,
ymax=1e-1,
yminorticks=true,
axis background/.style={fill=white},
title style={font=\bfseries},
xmajorgrids,
xminorgrids,
ymajorgrids,
yminorgrids,
legend style={at={(0.67,0.99)},legend cell align=left, draw=white!15!black, font=\footnotesize}
]

\addplot [color=mycyan, line width=2.0pt, mark=*]
  table[row sep=crcr]{%
   2.7570e-01   4.5895e-02 \\
   1.3785e-01   1.6042e-02 \\
   6.8924e-02   4.8553e-03 \\
   3.4462e-02   1.3440e-03 \\
   1.7231e-02   3.5413e-04 \\
   8.6155e-03   9.0906e-05 \\
   4.3077e-03   2.3021e-05 \\
   2.1539e-03   5.7869e-06 \\
   1.0769e-03   1.4478e-06 \\
};
%\addlegendentry{\calvoduemilauno}

\addplot [color=mygreen, line width=2.0pt, mark=*]
  table[row sep=crcr]{%
   2.7570e-01   5.8574e-04 \\
   1.3785e-01   4.0609e-05 \\
   6.8924e-02   1.1322e-06 \\
   3.4462e-02   8.2394e-07 \\
   1.7231e-02   1.6073e-07 \\
   8.6155e-03   2.4592e-08 \\
   4.3077e-03   3.5356e-09 \\
   2.1539e-03   5.5261e-10 \\
   1.0769e-03   1.1711e-10 \\
};
%\addlegendentry{\kennedycarpenterduemilatre}

\addplot [color=mygreen2, line width=2.0pt, mark=*]
  table[row sep=crcr]{%
   2.7570e-01   1.5444e-03 \\
   1.3785e-01   4.3477e-04 \\
   6.8924e-02   9.1179e-05 \\
   3.4462e-02   1.5292e-05 \\
   1.7231e-02   2.2382e-06 \\
   8.6155e-03   3.0405e-07 \\
   4.3077e-03   3.9798e-08 \\
   2.1539e-03   5.1592e-09 \\
   1.0769e-03   7.0401e-10 \\
};
%\addlegendentry{\kennedycarpenterduemiladiciannove}

\addplot [color=magenta, line width=2.0pt, mark=square*]
  table[row sep=crcr]{%
   2.7570e-01   2.9019e-03 \\
   1.3785e-01   3.4989e-04 \\
   6.8924e-02   3.3552e-05 \\
   3.4462e-02   2.7450e-06 \\
   1.7231e-02   2.2622e-07 \\
   8.6155e-03   3.0640e-08 \\
   4.3077e-03   9.8993e-09 \\
   2.1539e-03   4.6312e-09 \\
   1.0769e-03   2.2987e-09 \\
};
%\addlegendentry{\nostroordinequattroprimo}

\addplot [color=purple, line width=2.0pt, mark=square*]
  table[row sep=crcr]{%
   2.7570e-01   2.2156e-03 \\
   1.3785e-01   2.6876e-04 \\
   6.8924e-02   2.4763e-05 \\
   3.4462e-02   1.9879e-06 \\
   1.7231e-02   1.9572e-07 \\
   8.6155e-03   4.4193e-08 \\
   4.3077e-03   1.8958e-08 \\
   2.1539e-03   9.3664e-09 \\
   1.0769e-03   4.6908e-09 \\
};
%\addlegendentry{\nostroordinequattrosecondo}

\addplot [color=violet, line width=2.0pt, mark=square*]
  table[row sep=crcr]{%
   2.7570e-01   1.7499e-03 \\
   1.3785e-01   1.9074e-04 \\
   6.8924e-02   1.5595e-05 \\
   3.4462e-02   1.0817e-06 \\
   1.7231e-02   7.2810e-08 \\
   8.6155e-03   6.3005e-09 \\
   4.3077e-03   1.2482e-09 \\
   2.1539e-03   5.0121e-10 \\
   1.0769e-03   2.4252e-10 \\
};
%\addlegendentry{\nostroordinequattroterzo}

\addplot [color=black, dashed, line width=2.0pt]
  table[row sep=crcr]{
   1.0000e-00   1.0000e-04 \\
   1.0000e-03   1.0000e-16 \\
};
\node[right, align=left, text=black]
at (axis cs:0.015,1.7e-12) {$\sim\mathcal{O}(\Delta t^4)$};

\end{axis}
\end{tikzpicture}%
}
        \caption{Fourth-order methods: algebraic variable}
        \label{fig:tc2_errors4_alg}
    \end{subfigure} 
    \caption{Test case 2, van der Pol oscillator: convergence of different IMEX-RK schemes. The error is computed at final time on the differential variable for third-order methods (a) and fourth-order methods (b), and on the algebraic variable for third-order methods (c) and fourth-order methods (d).}
    \label{fig:tc2_errors}
\end{figure}
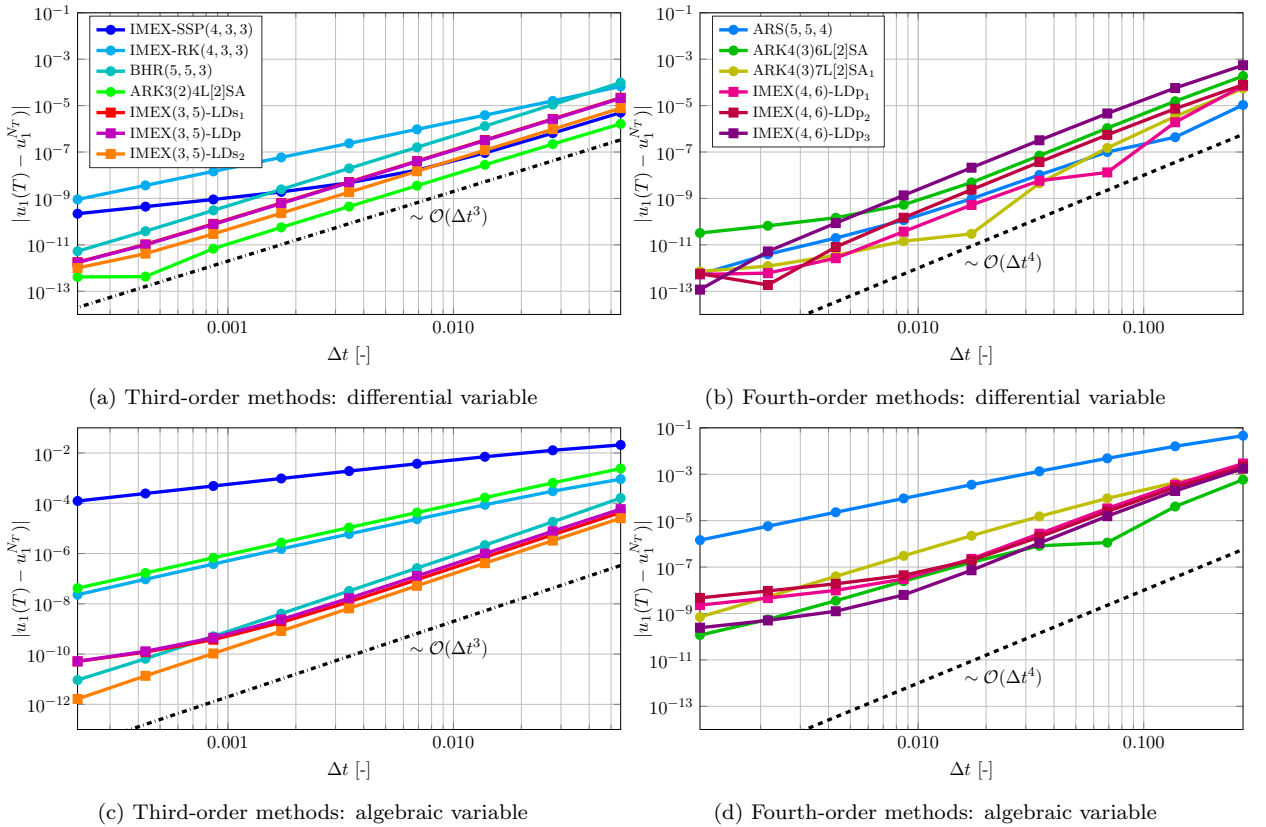

First, we compare the third-order IMEX-RK schemes on both the differential variable $u_{1}$ and the algebraic variable $u_{2}$ (see Figures~\ref{fig:tc2_errors3_diff}--\ref{fig:tc2_errors3_alg}). Some of the existing schemes from the literature exhibit noticeably larger global errors over the tested time-step range and show a more pronounced pre-asymptotic regime. In contrast, the Kennedy--Carpenter method from~\cite{kennedy:2003} and our new schemes display clean third-order convergence on both components. Moreover, on the algebraic variable our schemes systematically deliver smaller errors than the Kennedy--Carpenter method at comparable time steps, so that they not only retain their nominal order but also provide a more accurate resolution of the stiff dynamics of the van der Pol system.
\par
A similar behaviour is observed in the fourth-order comparison. The widely used ARS(5,5,4) method from~\cite{calvo:2001} shows a loss of order and relatively large errors on the algebraic variable (see Figure~\ref{fig:tc2_errors4_alg}), whereas the ARK4(3)6L[2]SA method from~\cite{kennedy:2003} suffers from a longer pre-asymptotic phase on the differential variable (see Figure~\ref{fig:tc2_errors4_diff}). In contrast, the three newly constructed fourth-order IMEX-RK schemes yield more robust error curves: they follow the expected fourth-order trend and only exhibit a loss of order in the algebraic component once the error in the differential variable falls below $10^{-12}$, which is likely due to residual inaccuracies of the nonlinear solver that are difficult to reduce further in this stiff regime. Overall, the performance of the proposed methods is comparable to that of the ARK4(3)7L[2]SA$_{1}$ scheme in~\cite{kennedy:2019}, while requiring one stage fewer.

%%%%%%%%%%%%%%%%%%%%%%%%%%%%%%%%%%%%%%%
%%%%%%%%%%% Test cases: PDEs %%%%%%%%%%
%%%%%%%%%%%%%%%%%%%%%%%%%%%%%%%%%%%%%%%
\section{Numerical results: application to PDEs}
\label{sec:numerical_results_PDE}

In this section, we present two applications of the IMEX methods to diffusion-reaction PDEs. First, we address the Fisher--Kolmogorov equation, a fundamental model in biological applications~\cite{corti:2024_IMEX, weickenmeier:2019}, which admits travelling wave solutions~\cite{antonietti:2026, corti:2024_IMEX, corti:2024_SP}. Finally, we consider a two-species diffusion-reaction system, the so-called Gray--Scott model~\cite{gray:1984, pearson:1993}, which exhibits complex patterns. The implementation of the polydG method described in Section~\ref{sec:polydg} is carried out within the \texttt{lymph} library framework~\cite{antonietti:2025}.

\subsection{Travelling-wave solution of FK equation}
\label{ssec:travelling_wave_FK}

In this section, we analyse the accuracy in simulating a travelling-wave solution of the Fisher--Kolmogorov (FK) problem. FK is a scalar problem ($n=1$) that can be rewritten as in Equation~\eqref{eq:model}, provided $F_{1}(u_{1}) = \alpha u_{1}(1-u_{1})$. Namely, considering $\FKcomp = u_{1}$ and $\sigma = \sigma_{1}$, the FK model reads
\begin{subequations}\label{eq:fisher_kpp}
\begin{alignat}{3}
    \pad{\FKcomp}{t} &= \dive\rpth{\sigma\grad \FKcomp} + \alpha\,\FKcomp\rpth{1 - \FKcomp} & & \quad \mathrm{in}\,\Omega \times (0,T), \\
    \sigma \grad \FKcomp \cdot \bm{n}_{\Omega} &= 0 & & \quad \mathrm{on}\,\partial\Omega \times (0,T), \\
    \FKcomp(\cdot, 0) &= \FKcomp_{0} & &\quad \mathrm{in}\,\Omega.
\end{alignat}
\end{subequations}
Model~\eqref{eq:fisher_kpp} represents a simplified formulation of the so-called heterodimer model~\cite{weickenmeier:2019}, which governs the dynamics of two different protein concentrations and is the classical mathematical description for neurodegenerative diseases. We refer, e.g., to~\cite{antonietti:2026, corti:2024_IMEX} for a more detailed discussion. In this section, we fix constant diffusion and reaction coefficients $\sigma$ and $\alpha$, and we consider travelling-wave solution of the form
\begin{equation*}
    \FKcomp(x,y,t) = \psi(x - vt),
\end{equation*}
where $v$ is wave speed depending on $\sigma$ and $\alpha$ defined by $v:= 5\sqrt{{\alpha\sigma}/{6}}$. Substituting $\FKcomp$ in the FK equation, we obtain an equivalent system of ODEs whose analytical solution is given by~\cite[\S7.2]{wen:2006}
\begin{equation*}
    \psi(x - vt) = \dfrac{1}{4}\rpth{1 + \tanh\rpth{8 - \sqrt{\dfrac{\alpha}{24\sigma}}(x - vt)}}^{2}.
\end{equation*}
This solution satisfies a homogeneous Neumann boundary condition at the limits $\xi\rightarrow \pm \infty$, which is equivalent to $x\rightarrow \pm \infty$ for each fixed value of $t \in (0,T]$. The homogeneous Neumann boundary condition is also respected in the $y$-direction, as the exact solution $c$ is independent of $y$. 
\par
In this simulation, we consider a rectangular space domain $\Omega = (0,3) \times (0,1)$, and the final time $T = 4$. Homogeneous Neumann boundary conditions are imposed not only at $y=0$ and $y=1$, but also at $x=0$ and $x=3$. Concerning the physical parameters, we fix $\sigma = 10^{-3}$ and $\alpha = 1$. The resulting velocity associated with these values is $v \simeq 6.45\times 10^{-2}$.
\par
The computational domain is discretized with a fixed polygonal mesh consisting of 800 elements generated by PolyMesher~\cite{talischi:2012}. The polynomial degree is set to $\ell = 8$ to emphasize the temporal errors associated with the IMEX time-integration scheme. We investigate temporal convergence by considering time steps $\Delta t = 0.8, 0.4, 0.2, 0.1$.
\par
We also consider, for the first time in this test case, the SI-IMEX-RK approach presented in Section~\ref{ssec:SI_IMEX_RK_method}. To this end, we rewrite~\eqref{eq:fisher_kpp} in the form of a partitioned system~\eqref{eq:partitioned_model} as follows
\begin{alignat}{3}
    \pad{\FKcomp_{\mathrm{E}}}{t} &= \dive\rpth{\sigma\grad \FKcomp_{\mathrm{I}}} + \alpha\,\FKcomp_{\mathrm{I}}\rpth{1 - \FKcomp_{\mathrm{E}}} & & \quad \mathrm{in}\,\Omega \times (0,T), \\[6pt]
    \pad{\FKcomp_{\mathrm{I}}}{t} &= \dive\rpth{\sigma\grad \FKcomp_{\mathrm{I}}} + \alpha\,\FKcomp_{\mathrm{I}}\rpth{1 - \FKcomp_{\mathrm{E}}} & & \quad \mathrm{in}\,\Omega \times (0,T),
\end{alignat}
coupled with appropriate initial and boundary conditions as in~\eqref{eq:fisher_kpp}. Following the discretization in Algorithm \ref{alg:si_imex_rk}, $\FKcomp_{\mathrm{E}}$ denotes the component treated explicitly, whereas $\FKcomp_{\mathrm{I}}$ represents the component treated implicitly.

\begin{figure}[h!]
    \begin{subfigure}[b]{0.33\textwidth}
        \resizebox{\textwidth}{!}{\definecolor{SIred}{rgb}{0.75000,0.000000,0.000000}%
\definecolor{myblue}{rgb}{0.00000,0.750000,0.75000}%
\begin{tikzpicture}

\begin{axis}[%
width =3.00in,
height=2.50in,
at={(2.6in,1.099in)},
scale only axis,
xmode=log,
xmin = 1.00e-1,
xmax = 8.00e-1,
xminorticks=true,
xlabel = {$\Delta t$ [-]},
ylabel = {Error},
xticklabel={\pgfmathparse{exp(\tick)}\pgfmathprintnumber{\pgfmathresult}},
x tick label style={
/pgf/number format/.cd, fixed, fixed zerofill,
precision=2},
ymode=log,
ymin=1e-7,
ymax=1e+0,
yminorticks=true,
axis background/.style={fill=white},
title style={font=\bfseries},
xmajorgrids,
xminorgrids,
ymajorgrids,
yminorgrids,
legend style={at={(0.67,0.99)},legend cell align=left, draw=white!15!black, font=\footnotesize}
]

\addplot [color=red, line width=2.0pt, mark options = solid, mark=square*]
  table[row sep=crcr]{%
  8.0000e-01    4.8516e-04 \\
  4.0000e-01    5.6860e-05 \\ 
  2.0000e-01    6.9063e-06 \\ 
  1.0000e-01    8.5122e-07 \\
};
\addlegendentry{IMEX: $\|\FKcomp(T)-\FKcomp_h^{N_T}\|_{L^2(\Omega)}$}

\addplot [color=red, line width=2.0pt, dashed, mark options = solid, mark=square*]
  table[row sep=crcr]{%
  8.0000e-01    3.3354e-03 \\
  4.0000e-01    4.6689e-04 \\ 
  2.0000e-01    6.4568e-05 \\ 
  1.0000e-01    8.5697e-06 \\
};
\addlegendentry{IMEX: $\|\FKcomp(T)-\FKcomp_h^{N_T}\|_{\mathrm{dG}}$}

\addplot [color=SIred, line width=2.0pt, mark options = solid, mark=*]
  table[row sep=crcr]{%
  8.0000e-01    5.9681e-03 \\
  4.0000e-01    2.7688e-04 \\ 
  2.0000e-01    2.6625e-05 \\ 
  1.0000e-01    3.0511e-06 \\
};
\addlegendentry{SI-IMEX: $\|\FKcomp(T)-\FKcomp_h^{N_T}\|_{L^2(\Omega)}$}

\addplot [color=SIred, line width=2.0pt, dashed, mark options = solid, mark=*]
  table[row sep=crcr]{%
  8.0000e-01    4.0037e-02 \\
  4.0000e-01    2.7507e-03 \\ 
  2.0000e-01    3.1008e-04 \\ 
  1.0000e-01    3.7365e-05 \\
};
\addlegendentry{SI-IMEX: $\|\FKcomp(T)-\FKcomp_h^{N_T}\|_{\mathrm{dG}}$}

\addplot [color=black, dash dot, line width=2.0pt]
  table[row sep=crcr]{%
   1.0000e-00   2.0000e-04 \\
   1.0000e-03   2.0000e-13 \\
};
\node[right, align=left, text=black]
at (axis cs:0.25,1.7e-6) {$\sim\mathcal{O}(\Delta t^3)$};

\end{axis}
\end{tikzpicture}%
}
        \caption{\orderthreeLDsone}
        \label{fig:tc3_errors_ordine_3_1}
    \end{subfigure}%    
    \begin{subfigure}[b]{0.33\textwidth}
        \resizebox{\textwidth}{!}{\definecolor{myred}{rgb}{0.75000,0.000000,0.750000}%
\definecolor{SImyred}{rgb}{0.5000,0.000000,0.50000}%
\definecolor{myblue}{rgb}{0.00000,0.750000,0.75000}%
\begin{tikzpicture}

\begin{axis}[%
width =3.00in,
height=2.50in,
at={(2.6in,1.099in)},
scale only axis,
xmode=log,
xmin = 1.00e-1,
xmax = 8.00e-1,
xminorticks=true,
xlabel = {$\Delta t$ [-]},
ylabel = {Error},
xticklabel={\pgfmathparse{exp(\tick)}\pgfmathprintnumber{\pgfmathresult}},
x tick label style={
/pgf/number format/.cd, fixed, fixed zerofill,
precision=2},
ymode=log,
ymin=1e-7,
ymax=1e+0,
yminorticks=true,
axis background/.style={fill=white},
title style={font=\bfseries},
xmajorgrids,
xminorgrids,
ymajorgrids,
yminorgrids,
legend style={at={(0.67,0.99)},legend cell align=left, draw=white!15!black, font=\footnotesize}
]

\addplot [color=myred, line width=2.0pt, mark options = solid, mark=square*]
  table[row sep=crcr]{%
  8.0000e-01    5.1998e-04 \\
  4.0000e-01    5.6529e-05 \\ 
  2.0000e-01    6.4817e-06 \\ 
  1.0000e-01    7.7068e-07 \\
};
\addlegendentry{IMEX: $\|\FKcomp(T)-\FKcomp_h^{N_T}\|_{L^2(\Omega)}$}

\addplot [color=myred, line width=2.0pt, dashed, mark options = solid, mark=square*]
  table[row sep=crcr]{%
  8.0000e-01    3.3126e-03 \\
  4.0000e-01    4.1843e-04 \\ 
  2.0000e-01    5.5436e-05 \\ 
  1.0000e-01    7.2356e-06 \\
};
\addlegendentry{IMEX: $\|\FKcomp(T)-\FKcomp_h^{N_T}\|_{\mathrm{dG}}$}

\addplot [color=SImyred, line width=2.0pt, mark options = solid, mark=*]
  table[row sep=crcr]{%
  8.0000e-01    6.0136e-03 \\
  4.0000e-01    2.8433e-04 \\ 
  2.0000e-01    2.7488e-05 \\ 
  1.0000e-01    3.1513e-06 \\
};
\addlegendentry{SI-IMEX: $\|\FKcomp(T)-\FKcomp_h^{N_T}\|_{L^2(\Omega)}$}

\addplot [color=SImyred, line width=2.0pt, dashed, mark options = solid, mark=*]
  table[row sep=crcr]{%
  8.0000e-01    3.9721e-02 \\
  4.0000e-01    2.7802e-03 \\ 
  2.0000e-01    3.1502e-04 \\ 
  1.0000e-01    3.8004e-05 \\
};
\addlegendentry{SI-IMEX: $\|\FKcomp(T)-\FKcomp_h^{N_T}\|_{\mathrm{dG}}$}

\addplot [color=black, dash dot, line width=2.0pt]
  table[row sep=crcr]{%
   1.0000e-00   2.0000e-04 \\
   1.0000e-03   2.0000e-13 \\
};
\node[right, align=left, text=black]
at (axis cs:0.25,1.7e-6) {$\sim\mathcal{O}(\Delta t^3)$};

\end{axis}
\end{tikzpicture}%
}
        \caption{\orderthreeLDp}
        \label{ffig:tc3_errors_ordine_3_2}
    \end{subfigure}%
    \begin{subfigure}[b]{0.33\textwidth}
        \resizebox{\textwidth}{!}{\definecolor{SIorange}{rgb}{0.75000,0.375000,0.00000}%
\begin{tikzpicture}

\begin{axis}[%
width =3.00in,
height=2.50in,
at={(2.6in,1.099in)},
scale only axis,
xmode=log,
xmin = 1.00e-1,
xmax = 8.00e-1,
xminorticks=true,
xlabel = {$\Delta t$ [-]},
ylabel = {Error},
xticklabel={\pgfmathparse{exp(\tick)}\pgfmathprintnumber{\pgfmathresult}},
x tick label style={
/pgf/number format/.cd, fixed, fixed zerofill,
precision=2},
ymode=log,
ymin=1e-7,
ymax=1e+0,
yminorticks=true,
axis background/.style={fill=white},
title style={font=\bfseries},
xmajorgrids,
xminorgrids,
ymajorgrids,
yminorgrids,
legend style={at={(0.67,0.99)},legend cell align=left, draw=white!15!black, font=\footnotesize}
]

\addplot [color=orange, line width=2.0pt, mark options = solid, mark=square*]
  table[row sep=crcr]{%
  8.0000e-01    3.5614e-04 \\
  4.0000e-01    3.2204e-05 \\ 
  2.0000e-01    3.3087e-06 \\ 
  1.0000e-01    3.7801e-07 \\
};
\addlegendentry{IMEX: $\|\FKcomp(T)-\FKcomp_h^{N_T}\|_{L^2(\Omega)}$}

\addplot [color=orange, line width=2.0pt, dashed, mark options = solid, mark=square*]
  table[row sep=crcr]{%
  8.0000e-01    3.2986e-03 \\
  4.0000e-01    3.2371e-04 \\ 
  2.0000e-01    3.7529e-05 \\ 
  1.0000e-01    4.7112e-06 \\
};
\addlegendentry{IMEX: $\|\FKcomp(T)-\FKcomp_h^{N_T}\|_{\mathrm{dG}}$}

\addplot [color=SIorange, line width=2.0pt, mark options = solid, mark=*]
  table[row sep=crcr]{%
  8.0000e-01    4.4148e-03 \\
  4.0000e-01    1.9020e-04 \\ 
  2.0000e-01    2.0824e-05 \\ 
  1.0000e-01    2.5464e-06 \\
};
\addlegendentry{SI-IMEX: $\|\FKcomp(T)-\FKcomp_h^{N_T}\|_{L^2(\Omega)}$}

\addplot [color=SIorange, line width=2.0pt, dashed, mark options = solid, mark=*]
  table[row sep=crcr]{%
  8.0000e-01    2.6404e-02 \\
  4.0000e-01    2.0313e-03 \\ 
  2.0000e-01    2.5334e-04 \\ 
  1.0000e-01    3.1559e-05 \\
};
\addlegendentry{SI-IMEX: $\|\FKcomp(T)-\FKcomp_h^{N_T}\|_{\mathrm{dG}}$}

\addplot [color=black, dash dot, line width=2.0pt]
  table[row sep=crcr]{%
   1.0000e-00   2.0000e-04 \\
   1.0000e-03   2.0000e-13 \\
};
\node[right, align=left, text=black]
at (axis cs:0.25,1.7e-6) {$\sim\mathcal{O}(\Delta t^3)$};

\end{axis}
\end{tikzpicture}%
}
        \caption{\orderthreeLDstwo}
        \label{fig:tc3_errors_ordine_3_3}
    \end{subfigure}%
    \caption{Test case 3, travelling-wave solution of FK equation: comparison of three different versions of our third-order method in the linearly implicit (IMEX) and semi-implicit (SI-IMEX) versions. The error is computed at final time both in $L^2(\Omega)$ and dG norms.}
    \label{fig:tc3_errors_ordine_3}
\end{figure}
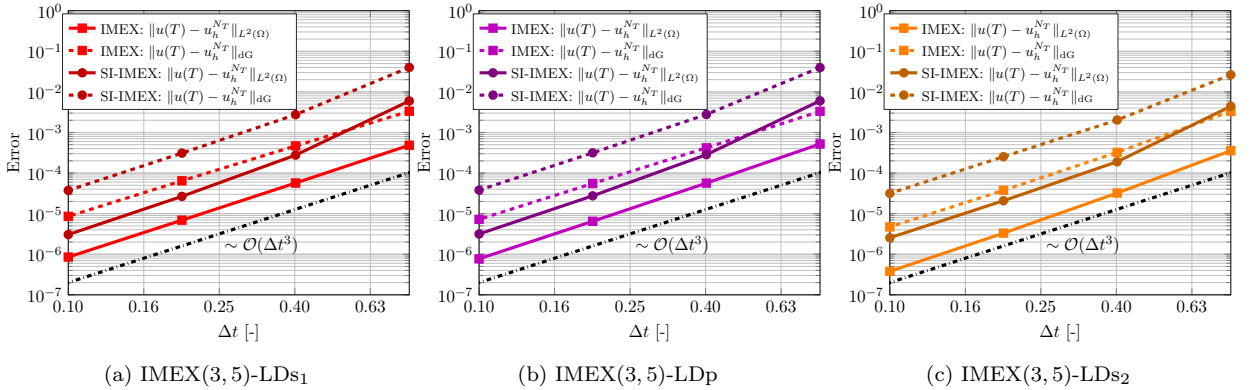

Figure~\ref{fig:tc3_errors_ordine_3} shows the error at the final time $T$ as a function of the time-step size $\Delta t$ for the three third-order IMEX-RK schemes proposed in this work. In each panel, we report the errors measured in the $L^{2}(\Omega)$ norm and in the dG norm, and we compare the standard linearly implicit IMEX formulation with its semi-implicit variant introduced in Section~\ref{ssec:SI_IMEX_RK_method}. The dashed reference line with slope three highlights the expected third-order convergence, which is clearly recovered by all schemes in the asymptotic regime.
\par
The curves also show that the SI-IMEX methods tend to exhibit slightly larger error constants than their IMEX counterparts, while preserving the same temporal convergence rate. Moreover, we observe that the \orderthreeLDstwo\ scheme, which attains stage order 2 for both third and fourth order stage of the explicit part, yields the smallest errors in magnitude, in particular for the linearly implicit formulation.

\begin{figure}[h!]
    \begin{subfigure}[b]{0.33\textwidth}
        \resizebox{\textwidth}{!}{\definecolor{SImagenta}{rgb}{0.75000,0.000000,0.750000}
\begin{tikzpicture}

\begin{axis}[%
width =3.00in,
height=2.50in,
at={(2.6in,1.099in)},
scale only axis,
xmode=log,
xmin = 1.00e-1,
xmax = 8.00e-1,
xminorticks=true,
xlabel = {$\Delta t$ [-]},
ylabel = {Error},
xticklabel={\pgfmathparse{exp(\tick)}\pgfmathprintnumber{\pgfmathresult}},
x tick label style={
/pgf/number format/.cd, fixed, fixed zerofill,
precision=2},
ymode=log,
ymin=1e-8,
ymax=1e-2,
yminorticks=true,
axis background/.style={fill=white},
title style={font=\bfseries},
xmajorgrids,
xminorgrids,
ymajorgrids,
yminorgrids,
legend style={at={(0.67,0.99)},legend cell align=left, draw=white!15!black, font=\footnotesize}
]

\addplot [color=magenta, line width=2.0pt, mark options = solid, mark=square*]
  table[row sep=crcr]{%
  8.0000e-01    5.3359e-04 \\
  4.0000e-01    2.4469e-05 \\ 
  2.0000e-01    1.0014e-06 \\ 
  1.0000e-01    4.3378e-08 \\
};
\addlegendentry{IMEX: $\|\FKcomp(T)-\FKcomp_h^{N_T}\|_{L^2(\Omega)}$}

\addplot [color=magenta, line width=2.0pt, dashed, mark options = solid, mark=square*]
  table[row sep=crcr]{%
  8.0000e-01    4.5996e-03 \\
  4.0000e-01    2.1445e-04 \\ 
  2.0000e-01    9.7196e-06 \\ 
  1.0000e-01    5.0037e-07 \\
};
\addlegendentry{IMEX: $\|\FKcomp(T)-\FKcomp_h^{N_T}\|_{\mathrm{dG}}$}

\addplot [color=SImagenta, line width=2.0pt, mark options = solid, mark=*]
  table[row sep=crcr]{%
  8.0000e-01    2.3795e-04 \\
  4.0000e-01    1.1115e-05 \\ 
  2.0000e-01    6.3692e-07 \\ 
  1.0000e-01    3.8610e-08 \\
};
\addlegendentry{SI-IMEX: $\|\FKcomp(T)-\FKcomp_h^{N_T}\|_{L^2(\Omega)}$}

\addplot [color=SImagenta, line width=2.0pt, dashed, mark options = solid, mark=*]
  table[row sep=crcr]{%
  8.0000e-01    1.9104e-03 \\
  4.0000e-01    8.8417e-05 \\ 
  2.0000e-01    5.2001e-06 \\ 
  1.0000e-01    3.2759e-07 \\
};
\addlegendentry{SI-IMEX: $\|\FKcomp(T)-\FKcomp_h^{N_T}\|_{\mathrm{dG}}$}

\addplot [color=black, dash dot, line width=2.0pt]
  table[row sep=crcr]{%
   1.0000e-00   2.0000e-04 \\
   1.0000e-03   2.0000e-16 \\
};
\node[right, align=left, text=black]
at (axis cs:0.25,4.0e-7) {$\sim\mathcal{O}(\Delta t^4)$};

\end{axis}
\end{tikzpicture}%
}
        \caption{\orderfourLDpone}
        \label{fig:tc3_errors_ordine_4_1}
    \end{subfigure}%    
    \begin{subfigure}[b]{0.33\textwidth}
        \resizebox{\textwidth}{!}{\definecolor{SIpurple}{rgb}{0.375000,0.000000,0.3750000}
\begin{tikzpicture}

\begin{axis}[%
width =3.00in,
height=2.50in,
at={(2.6in,1.099in)},
scale only axis,
xmode=log,
xmin = 1.00e-1,
xmax = 8.00e-1,
xminorticks=true,
xlabel = {$\Delta t$ [-]},
ylabel = {Error},
xticklabel={\pgfmathparse{exp(\tick)}\pgfmathprintnumber{\pgfmathresult}},
x tick label style={
/pgf/number format/.cd, fixed, fixed zerofill,
precision=2},
ymode=log,
ymin=1e-8,
ymax=1e-2,
yminorticks=true,
axis background/.style={fill=white},
title style={font=\bfseries},
xmajorgrids,
xminorgrids,
ymajorgrids,
yminorgrids,
legend style={at={(0.67,0.99)},legend cell align=left, draw=white!15!black, font=\footnotesize}
]

\addplot [color=purple, line width=2.0pt, mark options = solid, mark=square*]
  table[row sep=crcr]{%
  8.0000e-01    8.3308e-05 \\
  4.0000e-01    5.5580e-06 \\ 
  2.0000e-01    4.7485e-07 \\ 
  1.0000e-01    3.5730e-08 \\
};
\addlegendentry{IMEX: $\|\FKcomp(T)-\FKcomp_h^{N_T}\|_{L^2(\Omega)}$}

\addplot [color=purple, line width=2.0pt, dashed, mark options = solid, mark=square*]
  table[row sep=crcr]{%
  8.0000e-01    1.4554e-03 \\
  4.0000e-01    9.6083e-05 \\ 
  2.0000e-01    6.9339e-06 \\ 
  1.0000e-01    4.8319e-07 \\
};
\addlegendentry{IMEX: $\|\FKcomp(T)-\FKcomp_h^{N_T}\|_{\mathrm{dG}}$}

\addplot [color=SIpurple, line width=2.0pt, mark options = solid, mark=*]
  table[row sep=crcr]{%
  8.0000e-01    2.3380e-04 \\
  4.0000e-01    8.5050e-06 \\ 
  2.0000e-01    4.0784e-07 \\ 
  1.0000e-01    2.2332e-08 \\
};
\addlegendentry{SI-IMEX: $\|\FKcomp(T)-\FKcomp_h^{N_T}\|_{L^2(\Omega)}$}

\addplot [color=SIpurple, line width=2.0pt, dashed, mark options = solid, mark=*]
  table[row sep=crcr]{%
  8.0000e-01    1.7698e-03 \\
  4.0000e-01    6.8898e-05 \\ 
  2.0000e-01    3.4606e-06 \\ 
  1.0000e-01    2.0654e-07 \\
};
\addlegendentry{SI-IMEX: $\|\FKcomp(T)-\FKcomp_h^{N_T}\|_{\mathrm{dG}}$}

\addplot [color=black, dash dot, line width=2.0pt]
  table[row sep=crcr]{%
   1.0000e-00   2.0000e-04 \\
   1.0000e-03   2.0000e-16 \\
};
\node[right, align=left, text=black]
at (axis cs:0.25,4.0e-7) {$\sim\mathcal{O}(\Delta t^4)$};

\end{axis}
\end{tikzpicture}%
}
        \caption{\orderfourLDptwo}
        \label{ffig:tc3_errors_ordine_4_2}
    \end{subfigure}%
    \begin{subfigure}[b]{0.33\textwidth}
        \resizebox{\textwidth}{!}{\definecolor{SIviolet}{rgb}{0.25000,0.000000,0.750000}
\begin{tikzpicture}

\begin{axis}[%
width =3.00in,
height=2.50in,
at={(2.6in,1.099in)},
scale only axis,
xmode=log,
xmin = 1.00e-1,
xmax = 8.00e-1,
xminorticks=true,
xlabel = {$\Delta t$ [-]},
ylabel = {Error},
xticklabel={\pgfmathparse{exp(\tick)}\pgfmathprintnumber{\pgfmathresult}},
x tick label style={
/pgf/number format/.cd, fixed, fixed zerofill,
precision=2},
ymode=log,
ymin=1e-8,
ymax=1e-2,
yminorticks=true,
axis background/.style={fill=white},
title style={font=\bfseries},
xmajorgrids,
xminorgrids,
ymajorgrids,
yminorgrids,
legend style={at={(0.67,0.99)},legend cell align=left, draw=white!15!black, font=\footnotesize}
]

\addplot [color=violet, line width=2.0pt, mark options = solid, mark=square*]
  table[row sep=crcr]{%
  8.0000e-01    1.8424e-04 \\
  4.0000e-01    1.0898e-05 \\ 
  2.0000e-01    6.8683e-07 \\ 
  1.0000e-01    4.3477e-08 \\
};
\addlegendentry{IMEX: $\|\FKcomp(T)-\FKcomp_h^{N_T}\|_{L^2(\Omega)}$}

\addplot [color=violet, line width=2.0pt, dashed, mark options = solid, mark=square*]
  table[row sep=crcr]{%
  8.0000e-01    2.7466e-03 \\
  4.0000e-01    1.6253e-04 \\ 
  2.0000e-01    1.0184e-05 \\ 
  1.0000e-01    6.4597e-07 \\
};
\addlegendentry{IMEX: $\|\FKcomp(T)-\FKcomp_h^{N_T}\|_{\mathrm{dG}}$}

\addplot [color=SIviolet, line width=2.0pt, mark options = solid, mark=*]
  table[row sep=crcr]{%
  8.0000e-01    2.4006e-04 \\
  4.0000e-01    9.9400e-06 \\ 
  2.0000e-01    5.2460e-07 \\ 
  1.0000e-01    3.0534e-08 \\
};
\addlegendentry{SI-IMEX: $\|\FKcomp(T)-\FKcomp_h^{N_T}\|_{L^2(\Omega)}$}
 
\addplot [color=SIviolet, line width=2.0pt, dashed, mark options = solid, mark=*]
  table[row sep=crcr]{%
  8.0000e-01    2.2220e-03 \\
  4.0000e-01    1.2206e-04 \\ 
  2.0000e-01    7.2262e-06 \\ 
  1.0000e-01    4.4634e-07 \\
};
\addlegendentry{SI-IMEX: $\|\FKcomp(T)-\FKcomp_h^{N_T}\|_{\mathrm{dG}}$}

\addplot [color=black, dash dot, line width=2.0pt]
  table[row sep=crcr]{%
   1.0000e-00   2.0000e-04 \\
   1.0000e-03   2.0000e-16 \\
};
\node[right, align=left, text=black]
at (axis cs:0.25,4.0e-7) {$\sim\mathcal{O}(\Delta t^4)$};

\end{axis}
\end{tikzpicture}%
}
        \caption{\orderfourLDpthree}
        \label{fig:tc3_errors_ordine_4_3}
    \end{subfigure}%
    \caption{Test case 3, travelling-wave solution of FK equation: comparison of three different versions of our fourth-order method in the linearly implicit (IMEX) and semi-implicit (SI-IMEX) versions. The error is computed at final time both in $L^2(\Omega)$ and dG norms.}
    \label{fig:tc3_errors_ordine_4}
\end{figure}
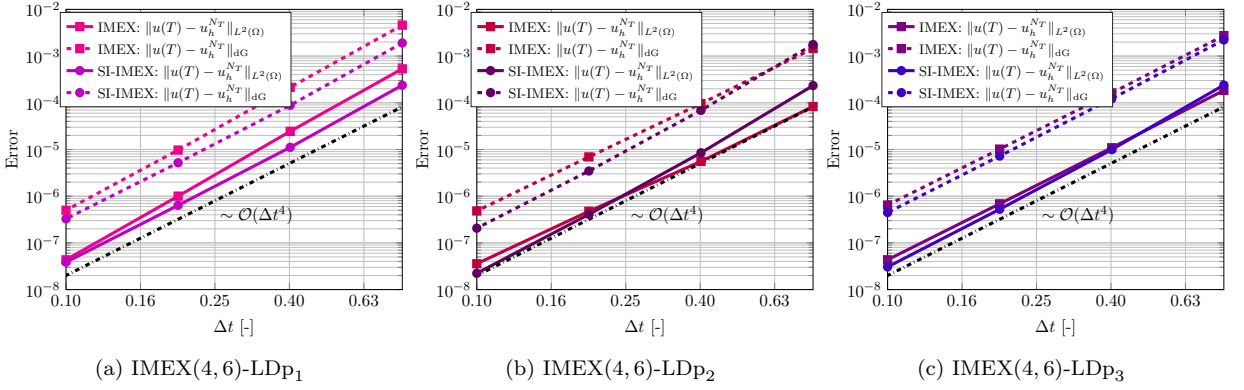

Figure~\ref{fig:tc3_errors_ordine_4} shows the error at the final time $T$ as a function of the time-step size $\Delta t$ for the three fourth-order IMEX-RK schemes presented in this work. The expected fourth-order convergence is clearly observed in the asymptotic regime for all schemes. The results indicate that all schemes yield errors of comparable magnitude. In addition, the plots show that, in this case, the SI-IMEX variants produce slightly smaller errors than their IMEX counterparts.
\par
\begin{figure}[t]
    \begin{subfigure}[b]{0.5\textwidth}
        \resizebox{\textwidth}{!}{\definecolor{myred}{rgb}{0.75000,0.000000,0.7500000}%
\definecolor{myblue}{rgb}{0.00000,0.750000,0.75000}%
\begin{tikzpicture}

\begin{axis}[%
width =3.50in,
height=2.50in,
at={(2.6in,1.099in)},
scale only axis,
xmode=log,
xmin = 1.0000e-01,
xmax = 8.0000e-01,
xminorticks=true,
xlabel = {$\Delta t$ [-]},
ylabel = {$\|u(T)-u_h^{N_T})\|_{L^2(\Omega)}$},
xticklabel={\pgfmathparse{exp(\tick)}\pgfmathprintnumber{\pgfmathresult}},
x tick label style={
/pgf/number format/.cd, fixed, fixed zerofill,
precision=3},
ymode=log,
ymin = 1.0000e-07,
ymax = 1.0000e-02,
yminorticks=true,
axis background/.style={fill=white},
title style={font=\bfseries},
xmajorgrids,
xminorgrids,
ymajorgrids,
yminorgrids,
legend style={at={(0.40,0.99)},legend cell align=left, draw=white!15!black, font=\footnotesize}
]

\addplot [color=blue, line width=2.0pt, mark=*]
  table[row sep=crcr]{%
  8.0000e-01    2.5397e-03 \\
  4.0000e-01    4.2402e-04 \\ 
  2.0000e-01    6.1732e-05 \\ 
  1.0000e-01    8.3359e-06 \\
};
\addlegendentry{\boscarinoduemilasedici}     

\addplot [color=cyan, line width=2.0pt, mark=*]
  table[row sep=crcr]{%
  8.0000e-01    2.0719e-03 \\
  4.0000e-01    4.0203e-04 \\ 
  2.0000e-01    1.4732e-04 \\ 
  1.0000e-01    4.5863e-05 \\
};
\addlegendentry{\boscarinoduemilaventidue}

\addplot [color=myblue, line width=2.0pt, mark=*]
  table[row sep=crcr]{%
  8.0000e-01    2.6685e-03 \\
  4.0000e-01    3.5539e-04 \\ 
  2.0000e-01    4.6554e-05 \\ 
  1.0000e-01    6.0038e-06 \\
};
\addlegendentry{\boscarinoduemilanove}

\addplot [color=green, line width=2.0pt, mark=*]
  table[row sep=crcr]{%
  8.0000e-01    8.5408e-04 \\
  4.0000e-01    1.2558e-04 \\ 
  2.0000e-01    1.7062e-05 \\ 
  1.0000e-01    2.2216e-06 \\
};
\addlegendentry{\carpenterkennedyduemilatre}

\addplot [color=red, line width=2.0pt, mark options = solid, mark=square*]
  table[row sep=crcr]{%
  8.0000e-01    4.8516e-04 \\
  4.0000e-01    5.6860e-05 \\ 
  2.0000e-01    6.9063e-06 \\ 
  1.0000e-01    8.5122e-07 \\
};
\addlegendentry{\nostroordinetreprimo}

\addplot [color=myred, line width=2.0pt, mark options = solid, mark=square*]
  table[row sep=crcr]{%
  8.0000e-01    5.1998e-04 \\
  4.0000e-01    5.6529e-05 \\ 
  2.0000e-01    6.4817e-06 \\ 
  1.0000e-01    7.7068e-07 \\
};
\addlegendentry{\nostroordinetresecondo}

\addplot [color=orange, line width=2.0pt, mark options = solid, mark=square*]
  table[row sep=crcr]{%
  8.0000e-01    3.5614e-04 \\
  4.0000e-01    3.2204e-05 \\ 
  2.0000e-01    3.3087e-06 \\ 
  1.0000e-01    3.7801e-07 \\
};
\addlegendentry{\nostroordinetreterzo}

\addplot [color=black, dash dot, line width=2.0pt]
  table[row sep=crcr]{%
   1.0000e-00   1.0000e-04 \\
   3.0000e-04   2.7000e-15 \\
};
\node[right, align=left, text=black]
at (axis cs:0.3,1.7e-6) {$\sim\mathcal{O}(\Delta t^3)$};

\end{axis}
\end{tikzpicture}%
}
        \caption{Third-order methods}
        \label{fig:tc1_errors3}
    \end{subfigure}%
    \begin{subfigure}[b]{0.5\textwidth}
        \resizebox{\textwidth}{!}{\definecolor{mygreen}{rgb}{0.00000,0.75000,0.0000}%
\definecolor{mygreen2}{rgb}{0.75000,0.75000,0.0000}%
\definecolor{mycyan}{rgb}{0.00000,0.50000,1.00000}%
\definecolor{mypurple}{rgb}{0.50000,0.00000,0.75000}%

\begin{tikzpicture}

\begin{axis}[%
width =3.50in,
height=2.50in,
at={(2.6in,1.099in)},
scale only axis,
xmode=log,
xmin = 1.00e-1,
xmax = 8.00e-1,
xminorticks=true,
xlabel = {$\Delta t$ [-]},
ylabel = {$|\Verhulstcomp(T)-\Verhulstcomp^{N_T})|$},
xticklabel={\pgfmathparse{exp(\tick)}\pgfmathprintnumber{\pgfmathresult}},
x tick label style={
/pgf/number format/.cd, fixed, fixed zerofill,
precision=3},
ymode=log,
ymin=5e-9,
ymax=1e-3,
yminorticks=true,
axis background/.style={fill=white},
title style={font=\bfseries},
xmajorgrids,
xminorgrids,
ymajorgrids,
yminorgrids,
legend style={at={(0.40,0.99)},legend cell align=left, draw=white!15!black, font=\footnotesize}
]

\addplot [color=mycyan, line width=2.0pt, mark=*]
  table[row sep=crcr]{%
  8.0000e-01    4.4159e-04 \\
  4.0000e-01    3.5194e-05 \\ 
  2.0000e-01    2.4929e-06 \\ 
  1.0000e-01    1.6590e-07 \\
};
\addlegendentry{\calvoduemilauno}

\addplot [color=mygreen, line width=2.0pt, mark=*]
  table[row sep=crcr]{%
  8.0000e-01    3.3772e-05 \\
  4.0000e-01    2.1173e-06 \\ 
  2.0000e-01    1.3686e-07 \\ 
  1.0000e-01    9.8789e-09 \\
};
\addlegendentry{\kennedycarpenterduemilatre}

\addplot [color=mygreen2, line width=2.0pt, mark=*]
  table[row sep=crcr]{%
  8.0000e-01    5.8789e-05 \\
  4.0000e-01    4.1389e-06 \\ 
  2.0000e-01    2.7333e-07 \\ 
  1.0000e-01    1.7532e-08 \\
};
\addlegendentry{\kennedycarpenterduemiladiciannove}

\addplot [color=magenta, line width=2.0pt, mark options = solid, mark=square*]
  table[row sep=crcr]{%
  8.0000e-01    5.3359e-04 \\
  4.0000e-01    2.4469e-05 \\ 
  2.0000e-01    1.0014e-06 \\ 
  1.0000e-01    4.3378e-08 \\
};
\addlegendentry{\nostroordinequattroprimo}

\addplot [color=purple, line width=2.0pt, mark options = solid, mark=square*]
  table[row sep=crcr]{%
  8.0000e-01    8.3308e-05 \\
  4.0000e-01    5.5580e-06 \\ 
  2.0000e-01    4.7485e-07 \\ 
  1.0000e-01    3.5730e-08 \\
};
\addlegendentry{\nostroordinequattrosecondo}

\addplot [color=violet, line width=2.0pt, mark options = solid, mark=square*]
  table[row sep=crcr]{%
  8.0000e-01    1.8424e-04 \\
  4.0000e-01    1.0898e-05 \\ 
  2.0000e-01    6.8683e-07 \\ 
  1.0000e-01    4.3477e-08 \\
};
\addlegendentry{\nostroordinequattroterzo}

\addplot [color=black, dashed, line width=2.0pt]
  table[row sep=crcr]{
   1.0000e-00   5.0000e-05 \\
   1.0000e-04   5.0000e-21 \\
};
\node[right, align=left, text=black]
at (axis cs:0.3,1.7e-7) {$\sim\mathcal{O}(\Delta t^4)$};

\end{axis}
\end{tikzpicture}%
}
        \caption{Fourth-order methods}
        \label{fig:tc1_errors4}
    \end{subfigure}%
    \caption{Test case 3, travelling-wave solution of FK equation: convergence of different IMEX-RK schemes. The $L^{2}(\Omega)$ error is computed at final time for third-order (a) and fourth-order methods (b).}
    \label{fig:tc1_errors}
\end{figure}
The observed behaviour may depend on the specific structure of the nonlinearity of the FK equation. The reaction term $\alpha\,\FKcomp\rpth{1 - \FKcomp}$ is globally bounded and exhibits moderate nonlinear character, so that the discrepancy between the duplicated variables $\FKcomp_\mathrm{E}^{(n,l)}$ and $\FKcomp_\mathrm{I}^{(n,l)}$ at intermediate stages remains well-controlled. At fourth order, the increased accuracy of the stage approximations results in error constants that are comparable to, or even smaller than, those of the linearly implicit IMEX counterpart. At third order, however, the coarser stage representation seems to yield a more pronounced disagreement between $\FKcomp_{\mathrm{E}}^{(n,l)}$ and $\FKcomp_{\mathrm{I}}^{(n,l)}$, ultimately resulting in systematically larger error magnitudes despite the two formulations sharing the same asymptotic convergence rate.
\par
Finally, we present a comparison with other IMEX schemes, in analogy to what we performed with the Verhulst equation in Section~\ref{ssec:reaction}. For third-order schemes, consistently with the results obtained for the Verhulst equation, the errors of our scheme are at least one order of magnitude lower than those of the existing schemes (Figure~\ref{fig:tc1_errors}, left). For fourth-order schemes, the results are again consistent with those obtained for the Verhulst model. Both ARK4(3)6L[2]SA and ARK4(3)7L[2]SA$_{1}$ perform slightly better than the fourth-order schemes proposed in this work (Figure~\ref{fig:tc1_errors}, right), likely because they were designed to minimize the leading error constant. The best-performing scheme is ARK4(3)6L[2]SA, followed by ARK4(3)7L[2]SA$_{1}$, which, we recall, requires one more stage than our schemes and is therefore more expensive per time step. Moreover, as shown in Section~\ref{ssec:van_der_pol}, ARK4(3)6L[2]SA exhibits a longer pre-asymptotic phase on the differential variable for the van der Pol equation. Our schemes, which satisfy the additional conditions~\eqref{eq:index_1_DAE_conditions_order_4}, may therefore be preferable for very stiff problems with an algebraic variable.

\subsection{Gray--Scott model}
\label{ssec:Gray-Scott}

As a final test, we consider the Gray--Scott model~\cite{gray:1984, pearson:1993}, which is widely used to study the formation of complex patterns arising from simple chemical reactions. The resulting reaction-diffusion system~\cite{pearson:1993} can be cast in the form~\eqref{eq:model}, by setting $\bm{u} = (\GSfirstcomp,\GSsecondcomp)^{\top}$,$\sigma_{1} = \sigma_{\GSfirstcomp}, \sigma_{2} = \sigma_{\GSsecondcomp}$, and
$$\mathbf{F}(\bm{u}, t) = \rpth{-\GSfirstcomp\,\GSsecondcomp^{2} + f\rpth{1 - \GSfirstcomp}, \GSfirstcomp\,\GSsecondcomp^{2} - \rpth{f + k}\GSsecondcomp}^{\top},$$
where $f$ and $k$ are characteristic (dimensionless) constant parameters of the the chemical reactions (see~\cite{pearson:1993} for further details). Hence, we obtain
\begin{subequations}
\label{eq:Gray_scott}
\begin{alignat}{3}
    \pad{\GSfirstcomp}{t} &\,= \dive\rpth{\sigma_{\GSfirstcomp}\grad \GSfirstcomp} - \GSfirstcomp\,\GSsecondcomp^{2} + f\rpth{1 - \GSfirstcomp} && \quad \mathrm{in}\,\Omega \times (0,T], \\[3pt]
    \pad{\GSsecondcomp}{t} &\,= \dive\rpth{\sigma_{\GSsecondcomp}\grad \GSsecondcomp} + \GSfirstcomp\,\GSsecondcomp^{2} - \rpth{f + k}\GSsecondcomp && \quad \mathrm{in}\,\Omega \times (0,T], \\[3pt]
    \sigma_{\GSfirstcomp}\grad \GSfirstcomp \cdot \bm{n}_{\Omega} = \sigma_{\GSsecondcomp}\grad \GSsecondcomp \cdot \bm{n}_{\Omega} &\, = 0 && \quad \mathrm{on}\,\Omega \times (0,T], \\[6pt]
    (\GSfirstcomp,\GSsecondcomp)(\cdot,0) & \,= (\GSfirstcomp[,0],\GSsecondcomp[,0]) && \quad \mathrm{in}\,\Omega.
\end{alignat}
\end{subequations}
The numerical simulation of the Gray--Scott model requires particular attention in the design of the time integration scheme. The diffusion forces fully explicit methods to adopt time steps $\Delta t = \mathcal{O}(h^{2})$. Since the accurate resolution of localised patterns such as spots and travelling waves typically demands fine spatial meshes, this restriction renders long-time simulations computationally very expensive. High-order IMEX Runge--Kutta methods overcome this limitation by treating the diffusion operator implicitly and the nonlinear reaction term explicitly, allowing time steps substantially larger than the explicit stability bound, while preserving both the qualitative and quantitative structure of the Gray--Scott patterns. These patterns thus provide a particularly demanding benchmark for the schemes proposed in this work.
\par
Due to the non-linearity, the term $\GSfirstcomp\,\GSsecondcomp^{2}$ is treated explicitly, whereas the remaining terms are discretized implicitly. However, thanks to the linearity with respect to $\GSfirstcomp$, the SI-IMEX-RK approach presented in Section~\ref{ssec:SI_IMEX_RK_method} allows for an implicit discretization of $\GSfirstcomp$, i.e. of the linear factor in the nonlinear term $\GSfirstcomp\,\GSsecondcomp^{2}$. In order to recast the system in the form~\eqref{eq:partitioned_model}, we introduce $(\GSfirstcomp[,\mathrm{E}], \GSsecondcomp[,\mathrm{E}])^{\top}$ and $(\GSfirstcomp[,\mathrm{I}], \GSsecondcomp[,\mathrm{I}])^{\top}$, and define
$$\mathbf{F}(\bm{u}_{\mathrm{E}}, \bm{u}_{\mathrm{I}}) = \rpth{-\GSfirstcomp[,\mathrm{I}]\,\GSsecondcomp[,\mathrm{E}]^{2} + f\rpth{1 - \GSfirstcomp[,\mathrm{I}]}, \GSfirstcomp[,\mathrm{I}]\,\GSsecondcomp[,\mathrm{E}]^{2} - (f + k)\GSsecondcomp[,\mathrm{I}]}^{\top}.$$
Hence, the resulting partitioned system writes
\begin{subequations}\label{eq:Gray_scott_SI}
\begin{alignat}{3}
    \pad{\GSfirstcomp[,\mathrm{E}]}{t} &\,= \dive\rpth{\sigma_{\GSfirstcomp}\grad \GSfirstcomp[,\mathrm{I}]} - \GSfirstcomp[,\mathrm{I}]\,\GSsecondcomp[,\mathrm{E}]^{2} + f\rpth{1 - \GSfirstcomp[,\mathrm{I}]} && \quad \mathrm{in}\,\Omega \times (0,T], \label{eq:Gray_scott_SI_1} \\[3pt]
    \pad{\GSfirstcomp[,\mathrm{I}]}{t} &\,= \dive\rpth{\sigma_{\GSfirstcomp}\grad \GSfirstcomp[,\mathrm{I}]} - \GSfirstcomp[,\mathrm{I}]\,\GSsecondcomp[,\mathrm{E}]^{2} + f\rpth{1 - \GSfirstcomp[,\mathrm{I}]} && \quad \mathrm{in}\,\Omega \times (0,T], \\[6pt]
    \pad{\GSsecondcomp[,\mathrm{E}]}{t} &\,= \dive\rpth{\sigma_{\GSsecondcomp}\grad \GSsecondcomp[,\mathrm{I}]} + \GSfirstcomp[,\mathrm{I}]\,\GSsecondcomp[,\mathrm{E}]^{2} - \rpth{f + k}\GSsecondcomp[,\mathrm{I}] && \quad \mathrm{in}\,\Omega \times (0,T], \\[3pt]
    \pad{\GSsecondcomp[,\mathrm{I}]}{t} &\,= \dive\rpth{\sigma_{\GSsecondcomp}\grad \GSsecondcomp[,\mathrm{I}]} + \GSfirstcomp[,\mathrm{I}]\,\GSsecondcomp[,\mathrm{E}]^{2} - \rpth{f + k}\GSsecondcomp[,\mathrm{I}] && \quad \mathrm{in}\,\Omega \times (0,T],
\end{alignat}
\end{subequations}
where Equation~\eqref{eq:Gray_scott_SI_1} can be neglected in practical computations due to the absence of $\GSfirstcomp[,\mathrm{E}]$ in the others.
\begin{figure}[t]
    \includegraphics[width=\textwidth]{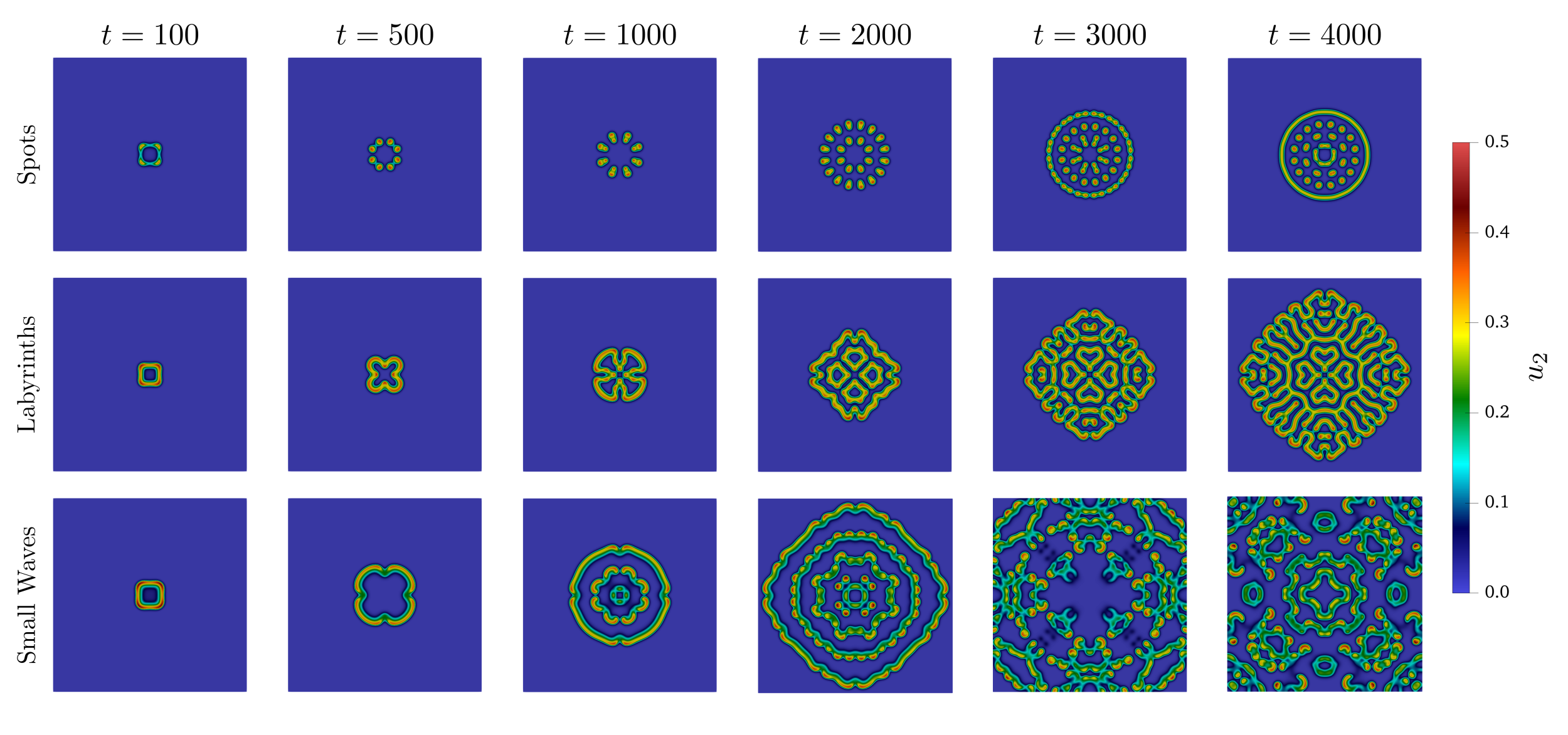}
    \caption{Test case 4, Gray--Scott model: numerical solutions at different times for three different values of the reaction parameters and with different patterns.}
    \label{fig:tc4_solutions}
\end{figure}
\begin{figure}[t]
    \begin{subfigure}[b]{0.5\textwidth}
        \resizebox{\textwidth}{!}{\definecolor{SIorange}{rgb}{0.75000,0.375000,0.00000}%
\begin{tikzpicture}
\begin{axis}[%
width=3.875in,
height=2.0in,
at={(2.6in,1.099in)},
scale only axis,
xmin=0.0,
xmax=0.7,
xlabel = {$x$},
ylabel = {$\GSsecondcomp$},
x tick label style={/pgf/number format/.cd, fixed, fixed zerofill, precision=1},
ymin=-0.02,
ymax=0.4,
ytick = {0, 0.1, 0.2, 0.3, 0.4},
y tick label style={/pgf/number format/.cd, fixed, fixed zerofill, precision=1},
axis background/.style={fill=white},
title style={font=\bfseries},
xmajorgrids,
xminorgrids,
ymajorgrids,
yminorgrids,
legend style={at={(0.28,0.25)},legend cell align=left, draw=white!15!black}
]
              
\addplot [color=orange, line width=2.0pt]
  table[row sep=crcr]{%
0.000000 -0.000000 \\ 
0.007071 -0.000003 \\ 
0.014142 -0.000004 \\ 
0.021213 -0.000004 \\ 
0.028284 -0.000006 \\ 
0.035355 -0.000012 \\ 
0.042426 -0.000023 \\ 
0.049497 0.000065 \\ 
0.056569 0.000107 \\ 
0.063640 -0.000010 \\ 
0.070711 -0.000169 \\ 
0.077782 -0.000274 \\ 
0.084853 -0.000290 \\ 
0.091924 0.001082 \\ 
0.098995 0.002000 \\ 
0.106066 0.002033 \\ 
0.113137 0.003525 \\ 
0.120208 0.009797 \\ 
0.127279 0.023137 \\ 
0.134350 0.043928 \\ 
0.141421 0.095114 \\ 
0.148492 0.171360 \\ 
0.155564 0.250272 \\ 
0.162635 0.308199 \\ 
0.169706 0.356254 \\ 
0.176777 0.365610 \\ 
0.183848 0.360302 \\ 
0.190919 0.357295 \\ 
0.197990 0.355124 \\ 
0.205061 0.349577 \\ 
0.212132 0.336294 \\ 
0.219203 0.312367 \\ 
0.226274 0.268446 \\ 
0.233345 0.203242 \\ 
0.240416 0.134461 \\ 
0.247487 0.078118 \\ 
0.254558 0.042302 \\ 
0.261629 0.019407 \\ 
0.268701 0.010458 \\ 
0.275772 0.004997 \\ 
0.282843 0.002749 \\ 
0.289914 0.001885 \\ 
0.296985 0.001797 \\ 
0.304056 0.001094 \\ 
0.311127 0.000043 \\ 
0.318198 -0.000132 \\ 
0.325269 -0.000111 \\ 
0.332340 -0.000001 \\ 
0.339411 0.000081 \\ 
0.346482 0.000091 \\ 
0.353553 -0.000098 \\ 
0.360624 0.000092 \\ 
0.367696 0.000082 \\ 
0.374767 -0.000001 \\ 
0.381838 -0.000112 \\ 
0.388909 -0.000132 \\ 
0.395980 0.000041 \\ 
0.403051 0.001095 \\ 
0.410122 0.001798 \\ 
0.417193 0.001884 \\ 
0.424264 0.002744 \\ 
0.431335 0.004988 \\ 
0.438406 0.010445 \\ 
0.445477 0.019408 \\ 
0.452548 0.042298 \\ 
0.459619 0.078063 \\ 
0.466691 0.134319 \\ 
0.473762 0.203032 \\ 
0.480833 0.268206 \\ 
0.487904 0.312190 \\ 
0.494975 0.336145 \\ 
0.502046 0.349437 \\ 
0.509117 0.355002 \\ 
0.516188 0.357181 \\ 
0.523259 0.360217 \\ 
0.530330 0.367830 \\ 
0.537401 0.356253 \\ 
0.544472 0.308229 \\ 
0.551543 0.250334 \\ 
0.558614 0.171440 \\ 
0.565685 0.095199 \\ 
0.572756 0.043983 \\ 
0.579827 0.023166 \\ 
0.586899 0.009813 \\ 
0.593970 0.003533 \\ 
0.601041 0.002036 \\ 
0.608112 0.002001 \\ 
0.615183 0.001083 \\ 
0.622254 -0.000290 \\ 
0.629325 -0.000274 \\ 
0.636396 -0.000170 \\ 
0.643467 -0.000010 \\ 
0.650538 0.000108 \\ 
0.657609 0.000066 \\ 
0.664680 -0.000023 \\ 
0.671751 -0.000012 \\ 
0.678823 -0.000006 \\ 
0.685894 -0.000004 \\ 
0.692965 -0.000003 \\ 
0.700036 -0.000003 \\ 
0.707107 -0.000001 \\ 
};
\addlegendentry{IMEX}

\addplot [color=SIorange, dashed, line width=2.0pt]
  table[row sep=crcr]{%
0.000000 -0.000000 \\ 
0.007071 -0.000003 \\ 
0.014142 -0.000004 \\ 
0.021213 -0.000004 \\ 
0.028284 -0.000006 \\ 
0.035355 -0.000012 \\ 
0.042426 -0.000023 \\ 
0.049497 0.000068 \\ 
0.056569 0.000112 \\ 
0.063640 -0.000008 \\ 
0.070711 -0.000175 \\ 
0.077782 -0.000289 \\ 
0.084853 -0.000320 \\ 
0.091924 0.001058 \\ 
0.098995 0.001935 \\ 
0.106066 0.002009 \\ 
0.113137 0.003796 \\ 
0.120208 0.010861 \\ 
0.127279 0.025580 \\ 
0.134350 0.048155 \\ 
0.141421 0.103578 \\ 
0.148492 0.181509 \\ 
0.155564 0.259766 \\ 
0.162635 0.315801 \\ 
0.169706 0.360235 \\ 
0.176777 0.364720 \\ 
0.183848 0.358708 \\ 
0.190919 0.355629 \\ 
0.197990 0.353622 \\ 
0.205061 0.348016 \\ 
0.212132 0.333775 \\ 
0.219203 0.307431 \\ 
0.226274 0.261387 \\ 
0.233345 0.195253 \\ 
0.240416 0.127367 \\ 
0.247487 0.072902 \\ 
0.254558 0.038894 \\ 
0.261629 0.017801 \\ 
0.268701 0.009737 \\ 
0.275772 0.004643 \\ 
0.282843 0.002574 \\ 
0.289914 0.001812 \\ 
0.296985 0.001755 \\ 
0.304056 0.001067 \\ 
0.311127 0.000042 \\ 
0.318198 -0.000127 \\ 
0.325269 -0.000108 \\ 
0.332340 -0.000001 \\ 
0.339411 0.000078 \\ 
0.346482 0.000088 \\ 
0.353553 -0.000095 \\ 
0.360624 0.000090 \\ 
0.367696 0.000080 \\ 
0.374767 -0.000001 \\ 
0.381838 -0.000109 \\ 
0.388909 -0.000129 \\ 
0.395980 0.000041 \\ 
0.403051 0.001069 \\ 
0.410122 0.001757 \\ 
0.417193 0.001812 \\ 
0.424264 0.002573 \\ 
0.431335 0.004642 \\ 
0.438406 0.009737 \\ 
0.445477 0.017830 \\ 
0.452548 0.038932 \\ 
0.459619 0.072892 \\ 
0.466691 0.127264 \\ 
0.473762 0.195054 \\ 
0.480833 0.261134 \\ 
0.487904 0.307214 \\ 
0.494975 0.333606 \\ 
0.502046 0.347906 \\ 
0.509117 0.353556 \\ 
0.516188 0.355605 \\ 
0.523259 0.358719 \\ 
0.530330 0.367250 \\ 
0.537401 0.360105 \\ 
0.544472 0.315574 \\ 
0.551543 0.259513 \\ 
0.558614 0.181259 \\ 
0.565685 0.103392 \\ 
0.572756 0.048081 \\ 
0.579827 0.025534 \\ 
0.586899 0.010844 \\ 
0.593970 0.003797 \\ 
0.601041 0.002015 \\ 
0.608112 0.001941 \\ 
0.615183 0.001061 \\ 
0.622254 -0.000319 \\ 
0.629325 -0.000288 \\ 
0.636396 -0.000175 \\ 
0.643467 -0.000008 \\ 
0.650538 0.000113 \\ 
0.657609 0.000068 \\ 
0.664680 -0.000023 \\ 
0.671751 -0.000012 \\ 
0.678823 -0.000006 \\ 
0.685894 -0.000004 \\ 
0.692965 -0.000004 \\ 
0.700036 -0.000003 \\ 
0.707107 -0.000001 \\ 
};
\addlegendentry{SI-IMEX}

\end{axis}
\end{tikzpicture}%
}
        \caption{\orderthreeLDstwo}
        \label{fig:tc4_comp_3}
    \end{subfigure}%    
    \begin{subfigure}[b]{0.5\textwidth}
        \resizebox{\textwidth}{!}{\definecolor{SIviolet}{rgb}{0.25000,0.000000,0.750000}
\begin{tikzpicture}
\begin{axis}[%
width=3.875in,
height=2.0in,
at={(2.6in,1.099in)},
scale only axis,
xmin=0.0,
xmax=0.7,
xlabel = {$x$},
ylabel = {$\GSsecondcomp$},
x tick label style={/pgf/number format/.cd, fixed, fixed zerofill, precision=1},
ymin=-0.02,
ymax=0.4,
ytick = {0, 0.1, 0.2, 0.3, 0.4},
y tick label style={/pgf/number format/.cd, fixed, fixed zerofill, precision=1},
axis background/.style={fill=white},
title style={font=\bfseries},
xmajorgrids,
xminorgrids,
ymajorgrids,
yminorgrids,
legend style={at={(0.28,0.25)},legend cell align=left, draw=white!15!black}
]
              
\addplot [color=violet, line width=2.0pt]
  table[row sep=crcr]{%
0.000000 -0.000000 \\ 
0.007071 -0.000003 \\ 
0.014142 -0.000004 \\ 
0.021213 -0.000004 \\ 
0.028284 -0.000006 \\ 
0.035355 -0.000012 \\ 
0.042426 -0.000023 \\ 
0.049497 0.000066 \\ 
0.056568 0.000109 \\ 
0.063640 -0.000010 \\ 
0.070711 -0.000171 \\ 
0.077782 -0.000278 \\ 
0.084853 -0.000299 \\ 
0.091924 0.001074 \\ 
0.098995 0.001980 \\ 
0.106066 0.002023 \\ 
0.113137 0.003595 \\ 
0.120208 0.010079 \\ 
0.127279 0.023787 \\ 
0.134350 0.045053 \\ 
0.141421 0.097432 \\ 
0.148492 0.174170 \\ 
0.155563 0.252934 \\ 
0.162635 0.310356 \\ 
0.169706 0.357398 \\ 
0.176777 0.365370 \\ 
0.183848 0.359844 \\ 
0.190919 0.356832 \\ 
0.197990 0.354733 \\ 
0.205061 0.349207 \\ 
0.212132 0.335685 \\ 
0.219203 0.311084 \\ 
0.226274 0.266496 \\ 
0.233345 0.200952 \\ 
0.240416 0.132377 \\ 
0.247487 0.076552 \\ 
0.254558 0.041255 \\ 
0.261629 0.018899 \\ 
0.268701 0.010229 \\ 
0.275772 0.004882 \\ 
0.282843 0.002692 \\ 
0.289914 0.001862 \\ 
0.296985 0.001785 \\ 
0.304056 0.001087 \\ 
0.311127 0.000042 \\ 
0.318198 -0.000131 \\ 
0.325269 -0.000111 \\ 
0.332340 -0.000001 \\ 
0.339411 0.000081 \\ 
0.346482 0.000091 \\ 
0.353553 -0.000097 \\ 
0.360625 0.000091 \\ 
0.367696 0.000081 \\ 
0.374767 -0.000001 \\ 
0.381838 -0.000111 \\ 
0.388909 -0.000132 \\ 
0.395980 0.000041 \\ 
0.403051 0.001088 \\ 
0.410122 0.001786 \\ 
0.417193 0.001861 \\ 
0.424264 0.002689 \\ 
0.431335 0.004878 \\ 
0.438406 0.010225 \\ 
0.445477 0.018921 \\ 
0.452548 0.041297 \\ 
0.459620 0.076566 \\ 
0.466691 0.132326 \\ 
0.473762 0.200842 \\ 
0.480833 0.266350 \\ 
0.487904 0.310974 \\ 
0.494975 0.335574 \\ 
0.502046 0.349097 \\ 
0.509117 0.354631 \\ 
0.516188 0.356745 \\ 
0.523259 0.359788 \\ 
0.530330 0.367690 \\ 
0.537401 0.357322 \\ 
0.544472 0.310242 \\ 
0.551543 0.252826 \\ 
0.558614 0.174090 \\ 
0.565685 0.097382 \\ 
0.572757 0.045042 \\ 
0.579828 0.023779 \\ 
0.586899 0.010079 \\ 
0.593970 0.003599 \\ 
0.601041 0.002027 \\ 
0.608112 0.001982 \\ 
0.615183 0.001076 \\ 
0.622254 -0.000298 \\ 
0.629325 -0.000278 \\ 
0.636396 -0.000171 \\ 
0.643467 -0.000009 \\ 
0.650538 0.000109 \\ 
0.657609 0.000066 \\ 
0.664680 -0.000023 \\ 
0.671751 -0.000012 \\ 
0.678822 -0.000006 \\ 
0.685894 -0.000004 \\ 
0.692965 -0.000004 \\ 
0.700036 -0.000003 \\ 
0.707107 -0.000001 \\
};
\addlegendentry{IMEX}

\addplot [color=SIviolet, dashed, line width=2.0pt]
  table[row sep=crcr]{%
0.000000 -0.000000 \\ 
0.007071 -0.000000 \\ 
0.014142 -0.000000 \\ 
0.021213 -0.000002 \\ 
0.028284 -0.000006 \\ 
0.035355 -0.000012 \\ 
0.042426 -0.000020 \\ 
0.049497 0.000101 \\ 
0.056569 0.000171 \\ 
0.063640 0.000027 \\ 
0.070711 -0.000190 \\ 
0.077782 -0.000355 \\ 
0.084853 -0.000447 \\ 
0.091924 0.000824 \\ 
0.098995 0.001338 \\ 
0.106066 0.002039 \\ 
0.113137 0.006581 \\ 
0.120208 0.020021 \\ 
0.127279 0.044957 \\ 
0.134350 0.079332 \\ 
0.141421 0.155441 \\ 
0.148492 0.235955 \\ 
0.155564 0.303606 \\ 
0.162635 0.345233 \\ 
0.169706 0.369724 \\ 
0.176777 0.355880 \\ 
0.183848 0.349530 \\ 
0.190919 0.347866 \\ 
0.197990 0.346947 \\ 
0.205061 0.340360 \\ 
0.212132 0.319596 \\ 
0.219203 0.279142 \\ 
0.226274 0.223354 \\ 
0.233345 0.154172 \\ 
0.240416 0.092208 \\ 
0.247487 0.047923 \\ 
0.254558 0.023118 \\ 
0.261629 0.010446 \\ 
0.268701 0.006168 \\ 
0.275772 0.002841 \\ 
0.282843 0.001636 \\ 
0.289914 0.001368 \\ 
0.296985 0.001479 \\ 
0.304056 0.000897 \\ 
0.311127 0.000020 \\ 
0.318198 -0.000118 \\ 
0.325269 -0.000097 \\ 
0.332340 -0.000004 \\ 
0.339411 0.000065 \\ 
0.346482 0.000075 \\ 
0.353553 -0.000079 \\ 
0.360624 0.000072 \\ 
0.367696 0.000063 \\ 
0.374767 -0.000005 \\ 
0.381838 -0.000096 \\ 
0.388909 -0.000115 \\ 
0.395980 0.000020 \\ 
0.403051 0.000903 \\ 
0.410122 0.001489 \\ 
0.417193 0.001379 \\ 
0.424264 0.001656 \\ 
0.431335 0.002880 \\ 
0.438406 0.006249 \\ 
0.445477 0.010624 \\ 
0.452548 0.023548 \\ 
0.459619 0.048646 \\ 
0.466691 0.093292 \\ 
0.473762 0.155550 \\ 
0.480833 0.224811 \\ 
0.487904 0.280485 \\ 
0.494975 0.320396 \\ 
0.502046 0.340851 \\ 
0.509117 0.347316 \\ 
0.516188 0.348233 \\ 
0.523259 0.349951 \\ 
0.530330 0.359510 \\ 
0.537401 0.369357 \\ 
0.544472 0.343795 \\ 
0.551543 0.301302 \\ 
0.558614 0.232904 \\ 
0.565685 0.152343 \\ 
0.572756 0.077202 \\ 
0.579827 0.043582 \\ 
0.586899 0.019336 \\ 
0.593970 0.006344 \\ 
0.601041 0.002011 \\ 
0.608112 0.001368 \\ 
0.615183 0.000835 \\ 
0.622254 -0.000443 \\ 
0.629325 -0.000353 \\ 
0.636396 -0.000190 \\ 
0.643467 0.000023 \\ 
0.650538 0.000166 \\ 
0.657609 0.000098 \\ 
0.664680 -0.000020 \\ 
0.671751 -0.000012 \\ 
0.678823 -0.000006 \\ 
0.685894 -0.000002 \\ 
0.692965 -0.000000 \\ 
0.700036 -0.000001 \\ 
0.707107 -0.000003 \\ 
};
\addlegendentry{SI-IMEX}

\end{axis}
\end{tikzpicture}%
}
        \caption{\orderthreeLDstwo}
        \label{fig:tc4_comp_4}
    \end{subfigure}
    \caption{Test case 4, Gray--Scott model: profile of $\GSsecondcomp$ along the diagonal of the computational domain at $T = 200$ with $\Delta t = 0.5$. Comparison between the IMEX and SI-IMEX schemes for the third-order (a) and fourth-order (b).}
    \label{fig:tc4_comp}
\end{figure}
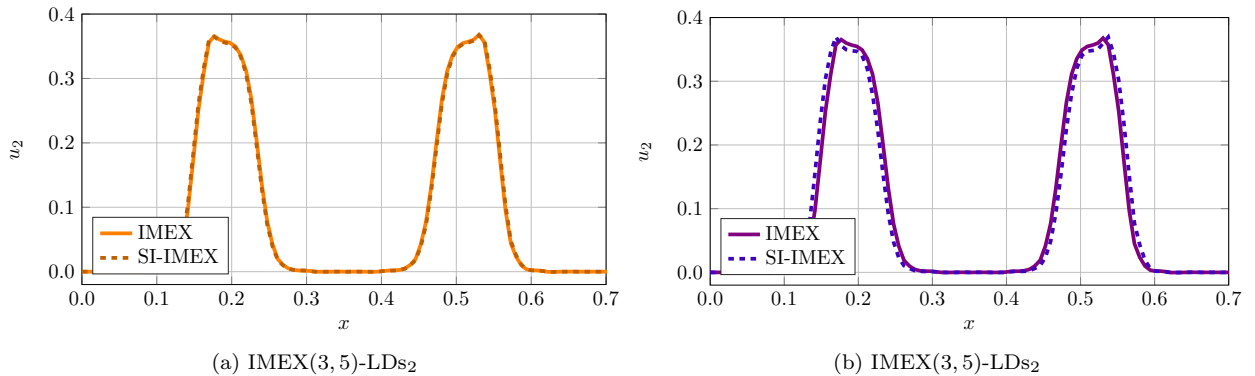
\par
For the numerical test we consider a two–dimensional square domain $\Omega = (0,2.5)^{2}$ with diffusion coefficients $\sigma_{\GSfirstcomp} = 2\times 10^{-5}$ and $\sigma_{\GSsecondcomp} = 10^{-5}$. For the spatial discretisation we consider a polynomial degree $\ell = 5$ on a Cartesian mesh of $25\,600$ elements. The code is implemented in \texttt{lymph}~\cite{antonietti:2025}. As initial data we consider $\GSfirstcomp = 1$ throughout the domain and a localised square seed for the second component, with $\GSsecondcomp = 0.5$ in $\Omega_{\mathrm{seed}} = \{(x,y) \in \mathbb{R}^{2}: 1.1562 \leq x \leq 1.3438 \wedge 1.1562 \leq y \leq 1.3438\}$ and $\GSsecondcomp = 0$ elsewhere. To illustrate the ability of the scheme to reproduce qualitatively different dynamics until the maximum time $T = 5\,000$, we consider three choices of the reaction parameters~\cite{pearson:1993}:
\begin{itemize}
    \item $(f,k) = (0.030, 0.062)$, which produces a pattern of localised spots;
    \item $(f,k) = (0.037, 0.060)$, which generates labyrinthine structures;
    \item $(f,k) = (0.018, 0.051)$, which gives rise to small--amplitude wave patterns.
\end{itemize}
\begin{figure}[t]
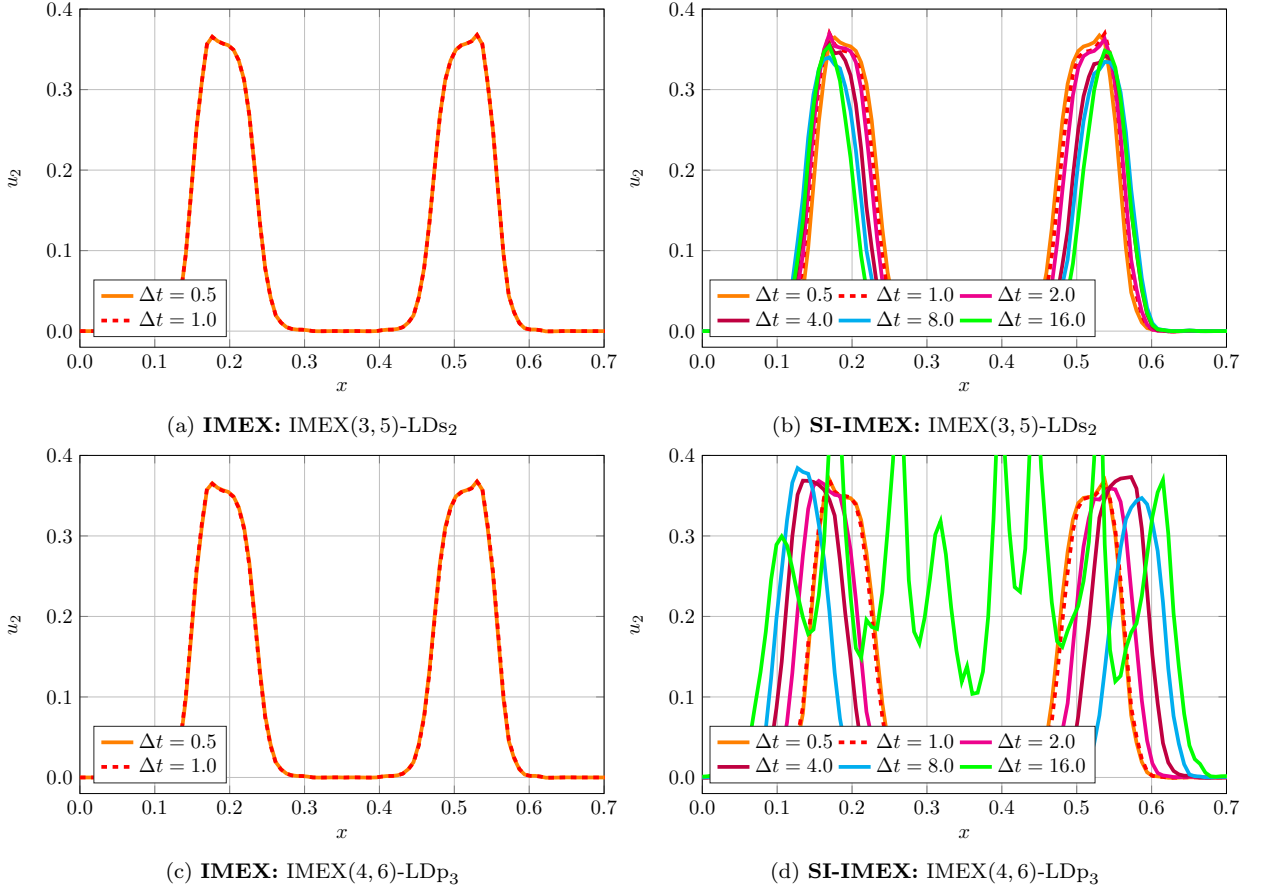

    \begin{subfigure}[b]{0.5\textwidth}
        \resizebox{\textwidth}{!}{\definecolor{mycolor2}{rgb}{0.00000,1.00000,1.00000}%
\begin{tikzpicture}
\begin{axis}[%
width=3.875in,
height=2.5in,
at={(2.6in,1.099in)},
scale only axis,
xmin=0.0,
xmax=0.7,
xlabel = {$x$},
ylabel = {$u_2$},
x tick label style={/pgf/number format/.cd, fixed, fixed zerofill, precision=1},
ymin=-0.02,
ymax=0.4,
ytick = {0, 0.1, 0.2, 0.3, 0.4},
y tick label style={/pgf/number format/.cd, fixed, fixed zerofill, precision=1},
axis background/.style={fill=white},
title style={font=\bfseries},
xmajorgrids,
xminorgrids,
ymajorgrids,
yminorgrids,
legend style={at={(0.28,0.20)},legend cell align=left, draw=white!15!black}
]
              
\addplot [color=orange, line width=2.0pt]
  table[row sep=crcr]{%
0.000000 -0.000000 \\ 
0.007071 -0.000003 \\ 
0.014142 -0.000004 \\ 
0.021213 -0.000004 \\ 
0.028284 -0.000006 \\ 
0.035355 -0.000012 \\ 
0.042426 -0.000023 \\ 
0.049497 0.000065 \\ 
0.056569 0.000107 \\ 
0.063640 -0.000010 \\ 
0.070711 -0.000169 \\ 
0.077782 -0.000274 \\ 
0.084853 -0.000290 \\ 
0.091924 0.001082 \\ 
0.098995 0.002000 \\ 
0.106066 0.002033 \\ 
0.113137 0.003525 \\ 
0.120208 0.009797 \\ 
0.127279 0.023137 \\ 
0.134350 0.043928 \\ 
0.141421 0.095114 \\ 
0.148492 0.171360 \\ 
0.155564 0.250272 \\ 
0.162635 0.308199 \\ 
0.169706 0.356254 \\ 
0.176777 0.365610 \\ 
0.183848 0.360302 \\ 
0.190919 0.357295 \\ 
0.197990 0.355124 \\ 
0.205061 0.349577 \\ 
0.212132 0.336294 \\ 
0.219203 0.312367 \\ 
0.226274 0.268446 \\ 
0.233345 0.203242 \\ 
0.240416 0.134461 \\ 
0.247487 0.078118 \\ 
0.254558 0.042302 \\ 
0.261629 0.019407 \\ 
0.268701 0.010458 \\ 
0.275772 0.004997 \\ 
0.282843 0.002749 \\ 
0.289914 0.001885 \\ 
0.296985 0.001797 \\ 
0.304056 0.001094 \\ 
0.311127 0.000043 \\ 
0.318198 -0.000132 \\ 
0.325269 -0.000111 \\ 
0.332340 -0.000001 \\ 
0.339411 0.000081 \\ 
0.346482 0.000091 \\ 
0.353553 -0.000098 \\ 
0.360624 0.000092 \\ 
0.367696 0.000082 \\ 
0.374767 -0.000001 \\ 
0.381838 -0.000112 \\ 
0.388909 -0.000132 \\ 
0.395980 0.000041 \\ 
0.403051 0.001095 \\ 
0.410122 0.001798 \\ 
0.417193 0.001884 \\ 
0.424264 0.002744 \\ 
0.431335 0.004988 \\ 
0.438406 0.010445 \\ 
0.445477 0.019408 \\ 
0.452548 0.042298 \\ 
0.459619 0.078063 \\ 
0.466691 0.134319 \\ 
0.473762 0.203032 \\ 
0.480833 0.268206 \\ 
0.487904 0.312190 \\ 
0.494975 0.336145 \\ 
0.502046 0.349437 \\ 
0.509117 0.355002 \\ 
0.516188 0.357181 \\ 
0.523259 0.360217 \\ 
0.530330 0.367830 \\ 
0.537401 0.356253 \\ 
0.544472 0.308229 \\ 
0.551543 0.250334 \\ 
0.558614 0.171440 \\ 
0.565685 0.095199 \\ 
0.572756 0.043983 \\ 
0.579827 0.023166 \\ 
0.586899 0.009813 \\ 
0.593970 0.003533 \\ 
0.601041 0.002036 \\ 
0.608112 0.002001 \\ 
0.615183 0.001083 \\ 
0.622254 -0.000290 \\ 
0.629325 -0.000274 \\ 
0.636396 -0.000170 \\ 
0.643467 -0.000010 \\ 
0.650538 0.000108 \\ 
0.657609 0.000066 \\ 
0.664680 -0.000023 \\ 
0.671751 -0.000012 \\ 
0.678823 -0.000006 \\ 
0.685894 -0.000004 \\ 
0.692965 -0.000003 \\ 
0.700036 -0.000003 \\ 
0.707107 -0.000001 \\ 
};
\addlegendentry{$\Delta t=0.5$}

\addplot [color=red, line width=2.0pt, dashed]
  table[row sep=crcr]{%
0.000000 -0.000000 \\ 
0.007071 -0.000003 \\ 
0.014142 -0.000004 \\ 
0.021213 -0.000004 \\ 
0.028284 -0.000006 \\ 
0.035355 -0.000012 \\ 
0.042426 -0.000023 \\ 
0.049497 0.000066 \\ 
0.056569 0.000108 \\ 
0.063640 -0.000010 \\ 
0.070711 -0.000170 \\ 
0.077782 -0.000275 \\ 
0.084853 -0.000293 \\ 
0.091924 0.001080 \\ 
0.098995 0.001994 \\ 
0.106066 0.002030 \\ 
0.113137 0.003546 \\ 
0.120208 0.009879 \\ 
0.127279 0.023328 \\ 
0.134350 0.044259 \\ 
0.141421 0.095799 \\ 
0.148492 0.172200 \\ 
0.155564 0.251072 \\ 
0.162635 0.308846 \\ 
0.169706 0.356608 \\ 
0.176777 0.365550 \\ 
0.183848 0.360174 \\ 
0.190919 0.357158 \\ 
0.197990 0.355009 \\ 
0.205061 0.349467 \\ 
0.212132 0.336110 \\ 
0.219203 0.311988 \\ 
0.226274 0.267861 \\ 
0.233345 0.202562 \\ 
0.240416 0.133837 \\ 
0.247487 0.077649 \\ 
0.254558 0.041989 \\ 
0.261629 0.019255 \\ 
0.268701 0.010389 \\ 
0.275772 0.004963 \\ 
0.282843 0.002732 \\ 
0.289914 0.001879 \\ 
0.296985 0.001793 \\ 
0.304056 0.001092 \\ 
0.311127 0.000042 \\ 
0.318198 -0.000131 \\ 
0.325269 -0.000111 \\ 
0.332340 -0.000001 \\ 
0.339411 0.000081 \\ 
0.346482 0.000091 \\ 
0.353553 -0.000098 \\ 
0.360624 0.000092 \\ 
0.367696 0.000081 \\ 
0.374767 -0.000001 \\ 
0.381838 -0.000112 \\ 
0.388909 -0.000132 \\ 
0.395980 0.000041 \\ 
0.403051 0.001093 \\ 
0.410122 0.001795 \\ 
0.417193 0.001877 \\ 
0.424264 0.002727 \\ 
0.431335 0.004954 \\ 
0.438406 0.010376 \\ 
0.445477 0.019257 \\ 
0.452548 0.041985 \\ 
0.459619 0.077595 \\ 
0.466691 0.133695 \\ 
0.473762 0.202342 \\ 
0.480833 0.267623 \\ 
0.487904 0.311804 \\ 
0.494975 0.335961 \\ 
0.502046 0.349337 \\ 
0.509117 0.354887 \\ 
0.516188 0.357051 \\ 
0.523259 0.360099 \\ 
0.530330 0.367800 \\ 
0.537401 0.356599 \\ 
0.544472 0.308870 \\ 
0.551543 0.251118 \\ 
0.558614 0.172270 \\ 
0.565685 0.095879 \\ 
0.572756 0.044311 \\ 
0.579827 0.023356 \\ 
0.586899 0.009895 \\ 
0.593970 0.003553 \\ 
0.601041 0.002033 \\ 
0.608112 0.001996 \\ 
0.615183 0.001081 \\ 
0.622254 -0.000293 \\ 
0.629325 -0.000275 \\ 
0.636396 -0.000170 \\ 
0.643467 -0.000010 \\ 
0.650538 0.000108 \\ 
0.657609 0.000066 \\ 
0.664680 -0.000023 \\ 
0.671751 -0.000012 \\ 
0.678823 -0.000006 \\ 
0.685894 -0.000004 \\ 
0.692965 -0.000004 \\ 
0.700036 -0.000003 \\ 
0.707107 -0.000001 \\ 
};
\addlegendentry{$\Delta t=1.0$}

\end{axis}
\end{tikzpicture}%
}
        \caption{\textbf{IMEX:} \orderthreeLDstwo}
        \label{fig:tc4_imex_3}
    \end{subfigure}%    
    \begin{subfigure}[b]{0.5\textwidth}
        \resizebox{\textwidth}{!}{\input{TC4_SI_IMEX_3.tex}}
        \caption{\textbf{SI-IMEX:} \orderthreeLDstwo}
        \label{fig:tc4_si_imex_3}
    \end{subfigure}\\
    \begin{subfigure}[b]{0.5\textwidth}
        \resizebox{\textwidth}{!}{\definecolor{mycolor2}{rgb}{0.00000,1.00000,1.00000}%
\begin{tikzpicture}
\begin{axis}[%
width=3.875in,
height=2.5in,
at={(2.6in,1.099in)},
scale only axis,
xmin=0.0,
xmax=0.7,
xlabel = {$x$},
ylabel = {$u_2$},
x tick label style={/pgf/number format/.cd, fixed, fixed zerofill, precision=1},
ymin=-0.02,
ymax=0.4,
ytick = {0, 0.1, 0.2, 0.3, 0.4},
y tick label style={/pgf/number format/.cd, fixed, fixed zerofill, precision=1},
axis background/.style={fill=white},
title style={font=\bfseries},
xmajorgrids,
xminorgrids,
ymajorgrids,
yminorgrids,
legend style={at={(0.28,0.20)},legend cell align=left, draw=white!15!black}
]
              
\addplot [color=orange, line width=2.0pt]
  table[row sep=crcr]{%
0.000000 -0.000000 \\ 
0.007071 -0.000003 \\ 
0.014142 -0.000004 \\ 
0.021213 -0.000004 \\ 
0.028284 -0.000006 \\ 
0.035355 -0.000012 \\ 
0.042426 -0.000023 \\ 
0.049497 0.000066 \\ 
0.056568 0.000109 \\ 
0.063640 -0.000010 \\ 
0.070711 -0.000171 \\ 
0.077782 -0.000278 \\ 
0.084853 -0.000299 \\ 
0.091924 0.001074 \\ 
0.098995 0.001980 \\ 
0.106066 0.002023 \\ 
0.113137 0.003595 \\ 
0.120208 0.010079 \\ 
0.127279 0.023787 \\ 
0.134350 0.045053 \\ 
0.141421 0.097432 \\ 
0.148492 0.174170 \\ 
0.155563 0.252934 \\ 
0.162635 0.310356 \\ 
0.169706 0.357398 \\ 
0.176777 0.365370 \\ 
0.183848 0.359844 \\ 
0.190919 0.356832 \\ 
0.197990 0.354733 \\ 
0.205061 0.349207 \\ 
0.212132 0.335685 \\ 
0.219203 0.311084 \\ 
0.226274 0.266496 \\ 
0.233345 0.200952 \\ 
0.240416 0.132377 \\ 
0.247487 0.076552 \\ 
0.254558 0.041255 \\ 
0.261629 0.018899 \\ 
0.268701 0.010229 \\ 
0.275772 0.004882 \\ 
0.282843 0.002692 \\ 
0.289914 0.001862 \\ 
0.296985 0.001785 \\ 
0.304056 0.001087 \\ 
0.311127 0.000042 \\ 
0.318198 -0.000131 \\ 
0.325269 -0.000111 \\ 
0.332340 -0.000001 \\ 
0.339411 0.000081 \\ 
0.346482 0.000091 \\ 
0.353553 -0.000097 \\ 
0.360625 0.000091 \\ 
0.367696 0.000081 \\ 
0.374767 -0.000001 \\ 
0.381838 -0.000111 \\ 
0.388909 -0.000132 \\ 
0.395980 0.000041 \\ 
0.403051 0.001088 \\ 
0.410122 0.001786 \\ 
0.417193 0.001861 \\ 
0.424264 0.002689 \\ 
0.431335 0.004878 \\ 
0.438406 0.010225 \\ 
0.445477 0.018921 \\ 
0.452548 0.041297 \\ 
0.459620 0.076566 \\ 
0.466691 0.132326 \\ 
0.473762 0.200842 \\ 
0.480833 0.266350 \\ 
0.487904 0.310974 \\ 
0.494975 0.335574 \\ 
0.502046 0.349097 \\ 
0.509117 0.354631 \\ 
0.516188 0.356745 \\ 
0.523259 0.359788 \\ 
0.530330 0.367690 \\ 
0.537401 0.357322 \\ 
0.544472 0.310242 \\ 
0.551543 0.252826 \\ 
0.558614 0.174090 \\ 
0.565685 0.097382 \\ 
0.572757 0.045042 \\ 
0.579828 0.023779 \\ 
0.586899 0.010079 \\ 
0.593970 0.003599 \\ 
0.601041 0.002027 \\ 
0.608112 0.001982 \\ 
0.615183 0.001076 \\ 
0.622254 -0.000298 \\ 
0.629325 -0.000278 \\ 
0.636396 -0.000171 \\ 
0.643467 -0.000009 \\ 
0.650538 0.000109 \\ 
0.657609 0.000066 \\ 
0.664680 -0.000023 \\ 
0.671751 -0.000012 \\ 
0.678822 -0.000006 \\ 
0.685894 -0.000004 \\ 
0.692965 -0.000004 \\ 
0.700036 -0.000003 \\ 
0.707107 -0.000001 \\ 
};
\addlegendentry{$\Delta t=0.5$}

\addplot [color=red, line width=2.0pt, dashed]
  table[row sep=crcr]{%
0.000000 -0.000000 \\ 
0.007071 -0.000003 \\ 
0.014142 -0.000004 \\ 
0.021213 -0.000004 \\ 
0.028284 -0.000006 \\ 
0.035355 -0.000012 \\ 
0.042426 -0.000023 \\ 
0.049497 0.000067 \\ 
0.056568 0.000111 \\ 
0.063640 -0.000009 \\ 
0.070711 -0.000174 \\ 
0.077782 -0.000285 \\ 
0.084853 -0.000312 \\ 
0.091924 0.001061 \\ 
0.098995 0.001945 \\ 
0.106066 0.002008 \\ 
0.113137 0.003716 \\ 
0.120208 0.010562 \\ 
0.127279 0.024893 \\ 
0.134350 0.046957 \\ 
0.141421 0.101268 \\ 
0.148492 0.178769 \\ 
0.155563 0.257215 \\ 
0.162635 0.313764 \\ 
0.169706 0.359161 \\ 
0.176777 0.364950 \\ 
0.183848 0.359117 \\ 
0.190919 0.356099 \\ 
0.197990 0.354115 \\ 
0.205061 0.348596 \\ 
0.212132 0.334634 \\ 
0.219203 0.308879 \\ 
0.226274 0.263194 \\ 
0.233345 0.197122 \\ 
0.240416 0.128933 \\ 
0.247487 0.073984 \\ 
0.254558 0.039556 \\ 
0.261629 0.018082 \\ 
0.268701 0.009857 \\ 
0.275772 0.004697 \\ 
0.282843 0.002599 \\ 
0.289914 0.001823 \\ 
0.296985 0.001764 \\ 
0.304056 0.001073 \\ 
0.311127 0.000041 \\ 
0.318198 -0.000130 \\ 
0.325269 -0.000110 \\ 
0.332340 -0.000001 \\ 
0.339411 0.000079 \\ 
0.346482 0.000089 \\ 
0.353553 -0.000096 \\ 
0.360625 0.000089 \\ 
0.367696 0.000079 \\ 
0.374767 -0.000002 \\ 
0.381838 -0.000110 \\ 
0.388909 -0.000130 \\ 
0.395980 0.000040 \\ 
0.403051 0.001075 \\ 
0.410122 0.001766 \\ 
0.417193 0.001824 \\ 
0.424264 0.002599 \\ 
0.431335 0.004699 \\ 
0.438406 0.009866 \\ 
0.445477 0.018132 \\ 
0.452548 0.039659 \\ 
0.459620 0.074096 \\ 
0.466691 0.129021 \\ 
0.473762 0.197163 \\ 
0.480833 0.263190 \\ 
0.487904 0.308862 \\ 
0.494975 0.334573 \\ 
0.502046 0.348517 \\ 
0.509117 0.354043 \\ 
0.516188 0.356055 \\ 
0.523259 0.359099 \\ 
0.530330 0.367400 \\ 
0.537401 0.358961 \\ 
0.544472 0.313437 \\ 
0.551543 0.256842 \\ 
0.558614 0.178409 \\ 
0.565685 0.100994 \\ 
0.572757 0.046830 \\ 
0.579828 0.024819 \\ 
0.586899 0.010534 \\ 
0.593970 0.003713 \\ 
0.601041 0.002012 \\ 
0.608112 0.001950 \\ 
0.615183 0.001063 \\ 
0.622254 -0.000310 \\ 
0.629325 -0.000284 \\ 
0.636396 -0.000173 \\ 
0.643467 -0.000009 \\ 
0.650538 0.000111 \\ 
0.657609 0.000067 \\ 
0.664680 -0.000023 \\ 
0.671751 -0.000012 \\ 
0.678822 -0.000006 \\ 
0.685894 -0.000004 \\ 
0.692965 -0.000004 \\ 
0.700036 -0.000003 \\ 
0.707107 -0.000001 \\  
};
\addlegendentry{$\Delta t=1.0$}

\end{axis}
\end{tikzpicture}%
}
        \caption{\textbf{IMEX:} \orderfourLDpthree}
        \label{fig:tc4_imex_4}
    \end{subfigure}%   
    \begin{subfigure}[b]{0.5\textwidth}
        \resizebox{\textwidth}{!}{\input{TC4_SI_IMEX_4.tex}}
        \caption{\textbf{SI-IMEX:} \orderfourLDpthree}
        \label{fig:tc4_si_imex_4}
    \end{subfigure}
    \caption{Test case 4, Gray--Scott model: profile of $u_2$ along the diagonal of the computational domain at $T = 200$ for several values of the time step $\Delta t$, computed with the third-order schemes IMEX (a) and SI-IMEX (b), and the fourth-order schemes IMEX (c) and SI-IMEX (d).}
    \label{fig:tc4_shapes}
\end{figure}
\par
In Figure~\ref{fig:tc4_solutions}, we report the results obtained using the fourth--order \orderfourLDpthree\ scheme with a time step $\Delta t = 0.25$. The three parameter sets lead to the expected, qualitatively different patterns. In the first row (spots), a single initial bump destabilises into a ring of spots and then undergoes successive replications. In the first row (spots), a single initial bump destabilises into a ring of spots and then undergoes successive replications. In the second row (labyrinths), the initial annulus progressively folds and branches into a connected labyrinthine network. In the third row (small waves), the solution generates outward--propagating concentric waves that subsequently fragment into smaller wavelets. At later times, these wavelets interact with one another and with the boundary, giving rise to an increasingly irregular pattern.
\par
We next compare the IMEX and SI-IMEX formulations and examine their sensitivity to the time step. We perform the same simulation up to $T = 200$ with parameters $(f,k) = (0.030,0.062)$. In Figure~\ref{fig:tc4_comp}, we report the profile of $u_{2}$ along the diagonal connecting $(1,1)$ to $(1.5, 1.5)$ in the computational domain, for both the third-order scheme \orderthreeLDstwo\ and the fourth-order scheme \orderfourLDpthree\ for $\Delta t = 0.5$. We can see that for the third order, the IMEX and the SI-IMEX methods provide equivalent results. A small difference can be detected comparing the fourth-order methods.
\par
\begin{figure}[t]
    \centering
    \includegraphics[width=\textwidth]{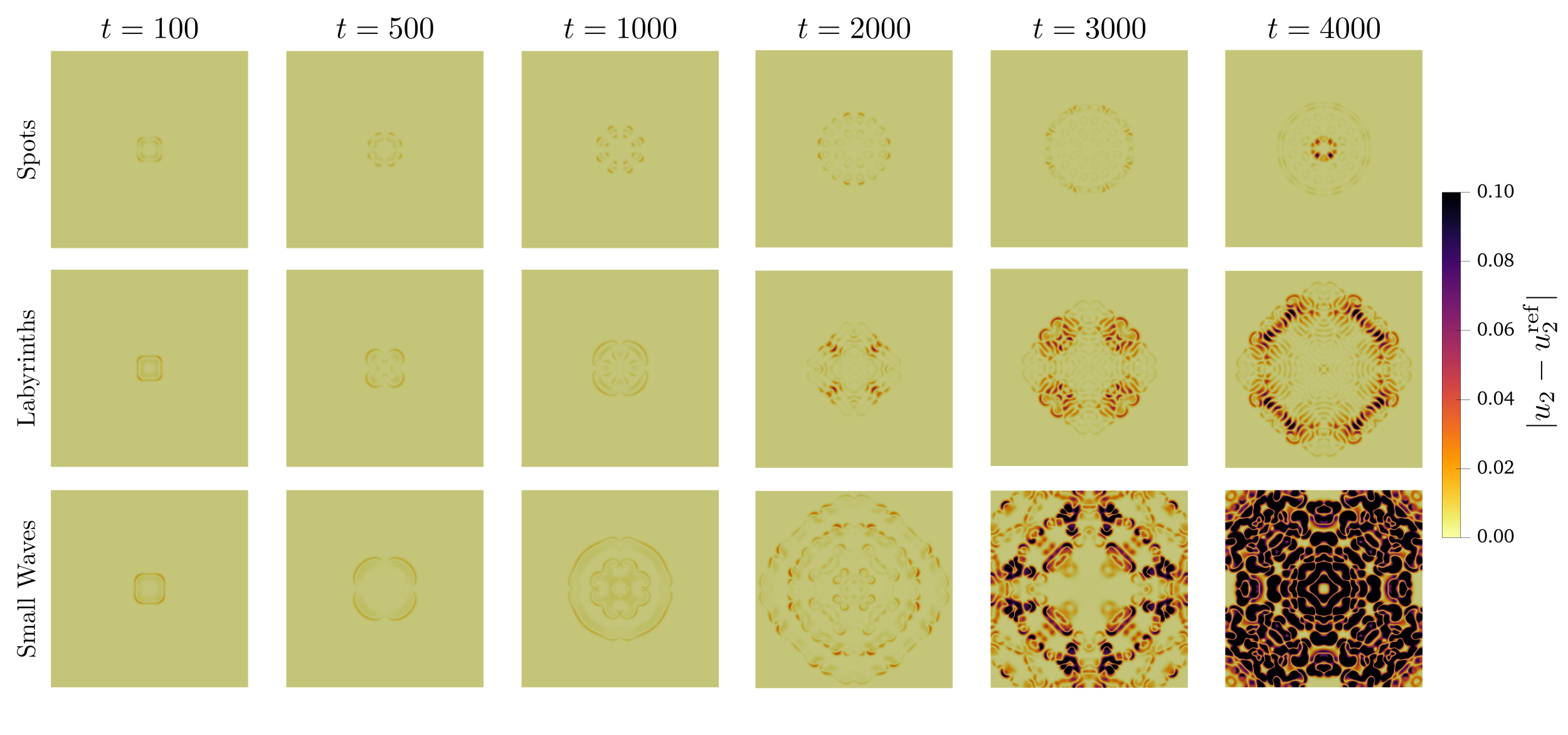}
    \caption{Pointwise absolute difference $|u_{2} - u_{2}^{\mathrm{ref}}|$ the IMEX-DG solution and the reference solution at different times (columns), for the spots, labyrinths and small waves cases (rows). The same colour scale is used in all panels.}
    \label{fig:GrayScottErrors}
\end{figure}
In Figure~\ref{fig:tc4_shapes}, we report the profile of $u_{2}$ for several values of the time step, namely $\Delta t \in \{0.5, 1.0, 2.0, 4.0, 8.0, 16.0\}$. For the standard IMEX formulations, the curves corresponding to $\Delta t = 0.5$ and $\Delta t = 1.0$ are virtually indistinguishable, confirming that both the propagation speed and the shape of the travelling fronts are accurately captured at small time steps; however, the method becomes unstable for larger values of $\Delta t$ and fails to produce a solution. In contrast, the SI-IMEX variants remain stable across the entire range of time steps considered, at the cost of a progressive deterioration of the wave profile: as $\Delta t$ increases, phase shifts and shape distortions become visible, and are particularly pronounced for the fourth-order scheme.

Finally, we compare the solutions for the three parameter sets with an independently computed reference solution. The reference is obtained in \texttt{FEniCSx} \cite{alnaes_fenics_2015} using conforming $\mathbb{Q}_{2}$ elements on a Cartesian mesh of $409\,600$ quadrilateral cells and fully implicit BDF2 time integration with $\Delta t_{\mathrm{ref}} = 0.025$. A backward Euler step is used to initialise the time integration. Figure~\ref{fig:GrayScottErrors} shows the pointwise absolute difference $|u_{2} - u_{2}^{\mathrm{ref}}|$ at selected times up to $T = 4\,000$, using the same colour scale for all panels. Up to $t = 2\,000$, the differences are mostly localised near the fronts and spots. At later times, they remain localised in the spot case but grow along the labyrinthine branches and become widespread in the wave case. The strongest discrepancies occur during the stage of multiple front interactions and contact with the boundary, although the maps alone do not establish whether these interactions are responsible for their growth. Overall, the method reproduces the expected patterns arising from Turing instability, while larger discrepancies with the reference solution develop during the more complex stages of wave interaction.

%%%%%%%%% Conclusions %%%%%%%%%%%%%%%%%
\section{Conclusions}
\label{sec:conclu}

In this work, we have presented IMplicit–EXplicit (IMEX) Runge--Kutta schemes coupled with a polytopal discontinuous Galerkin discretization for possibly degenerate diffusion–reaction problems. Since such models may admit travelling wave solutions, as in the Fisher--Kolmogorov equation, we have designed new third- and fourth-order IMEX–RK schemes tailored to improve the numerical description of wave propagation, with reduced artificial dissipation and dispersion. In addition, inspired by the partitioned formulation introduced in~\cite{boscarino:2016}, we have proposed a new class of semi-implicit IMEX Runge--Kutta (SI-IMEX-RK) methods. This approach avoids a rigid separation between linear and nonlinear reaction components, while enabling the selective treatment of the dominant (possibly linear) sources of stiffness. This strategy yields efficient and flexible schemes that can be naturally extended to more complex nonlinear dependencies, including those arising in the diffusion operator. The capabilities of the proposed methods have been assessed through a hierarchy of test cases of increasing complexity, targeting regimes characterized by strong stiffness and nonlinear dynamics. The experiments demonstrate that the considered schemes deliver accurate and stable solutions across a wide range of scenarios, including biologically relevant models.
\par
Possible future directions of this work are the rigorous theoretical analysis of the SI-IMEX-RK framework, establishing entropy stability properties at the discrete level and the extension of dispersion and dissipation analysis to IMEX methods, defining and optimizing these properties directly for the combined IMEX schemes on a suitable model problem, rather than separately on the implicit and explicit tableaux. A further promising direction is to couple the proposed time integration strategies with entropy-stable spatial discretisations based on summation-by-parts (SBP) operators~\cite{gassner:2013, gassner:2021}. Extending these SBP constructions to general polygonal meshes and combining them with IMEX time stepping is a challenging but natural step towards fully structure-preserving numerical methods.

%%%%%%%%%%%%%%%%%%%%%%% Appendix %%%%%%%%%%%%%%%%%%%%%%%%%%
\appendix

%%%%%%%% Used tableaux %%%%%%%%%%%%
\section{Coefficients of employed IMEX-RK schemes}
\label{app:IMEX_coeffs}

We report here, for the reader’s convenience, all the details on the IMEX–RK schemes employed in the numerical simulations, whose coefficients are approximated to the eighth decimal place. The interested reader may contact the authors to obtain the full-precision values.

\subsection{Third-order schemes}

Regarding the \emph{third-order} schemes, we list here the schemes considered.

\begin{table}[H]
    \small
    \begin{tabular}{c|cccc}
        $0$ & $0$ & $0$ & $0$ & $0$ \\[3pt]
        $0.43586652$ & %1508$ & 
        $\;\;\:0.43586652$ & %1508$ & 
        $0$ & $0$ & $0$ \\[3pt]
        $0.71793326$ & %0754$ & 
        $\;\;\:0.43586652$ & %1508$ & 
        $0.28206673$ & %9245$ & 
        $0$ & $0$ \\[3pt]
        $1$ & 
        $-0.73353408$ & %2748750$ & 
        $2.15052738$ & %1100$ & 
        $-0.41699329$ & %8352$ & 
        $0$ \\[3pt]
        \hline \\[-6pt]
        & $0$ & 
        $1.20849664$ & %9176$ & 
        $-0.64436317$ & %0684$ & 
        $0.43586652$ %1508$
    \end{tabular}
    \vskip 0.3cm
    \begin{tabular}{c|cccc}
        $0.43586652$ & %1508$ & 
        $0.43586652$ & %1508$ & 
        $0$ & $0$ & $0$ \\[3pt]
        $0.43586652$ & %1508$ & 
        $0$ & 
        $0.43586652$ & %1508$ 
        $0$ & $0$ \\[3pt]
        $0.71793326$ & %0754$ & 
        $0$ & 
        $0.28206673$ & %9245$ & 
        $\;\;\:0.43586652$ & %1508$ & 
        $0$ \\[3pt]
        $1$ & $0$ & 
        $1.20849664$ & %9176$ & 
        $-0.64436317$ & %0684$ & 
        $0.43586652$ %1508$ 
        \\[3pt]
        \hline \\[-6pt]
        & $0$ & 
        $1.20849664$ & %9176$ & 
        $-0.64436317$ & %0684$ & 
        $0.43586652$ %1508$
    \end{tabular}
    \caption{\textbf{Butcher tableaux of IMEX-RK(4,3,3) scheme employed in~\cite{boscarino:2022}.} Explicit method (top) and implicit method (bottom).}
    \label{tab:rk3_butch_type_I}
\end{table}

\begin{table}[H]
    \small
	\begin{tabular}{c|cccc}
		$0$ & $0$ & $0$ & $0$ & $0$ \\[3pt]
        $0.87173304$ & 
        $0.87173304$ & $0$ & $0$ & $0$ \\ [3pt]
        $0.60000000$ & 
        $0.52758901$ & 
        $0.07241098$ & $0$ & $0$ \\[3pt]
        $1$ & 
        $0.39909600$ & 
        $-0.43755765$ & 
        $1.03846164$ & $0$ \\[3pt]
        \hline \\[-6pt]
        & 
        $0.18764102$ & 
        $-0.59529747$ &
        $0.97178992$ & 
        $0.43586652$
	\end{tabular}
    \vskip 0.3cm
    \begin{tabular}{c|cccc}
		$0$ & $0$ & $0$ & $0$ & $0$ \\[3pt]
        $0.87173304$ & 
        $0.43586652$ & 
        $\;\;\:0.43586652$ & $0$ & $0$ \\[3pt]
		$0.60000000$ & 
        $0.25764824$ &
        $-0.09351476$ & 
        $0.43586652$ & $0$ \\[3pt]
        $1$ & 
        $0.18764102$ & 
        $-0.59529747$ &
        $0.97178992$ & 
        $0.43586652$ \\[3pt]
		\hline \\[-6pt]
		& 
        $0.18764102$ & 
        $-0.59529747$ &
        $0.97178992$ & 
        $0.43586652$
	\end{tabular}
    \caption{\textbf{Butcher tableaux of ARK3(2)4L[2]SA scheme employed in~\cite{kennedy:2003}.} Explicit method (top) and implicit method (bottom).}
    \label{tab:rk3_butch_type_II}
\end{table}

\begin{table}[h!]
    \begin{minipage}{0.5\textwidth}
	\begin{center}
		\begin{tabular}{>{\centering\arraybackslash}p{7mm}|>{\centering\arraybackslash}p{7mm}>{\centering\arraybackslash}p{7mm}>{\centering\arraybackslash}p{7mm}>{\centering\arraybackslash}p{7mm}}
			$0$ & $0$ & $0$ & $0$ & $0$ \\[2mm]
            $0$ & $0$ & $0$ & $0$ & $0$ \\[2mm]
            $1$ & $0$ & $1$ & $0$ & $0$ \\[2mm]
            $\frac{1}{2}$ & $0$ & $\frac{1}{4}$ & $\frac{1}{4}$ & $0$ \\ [2mm]
            \hline \\ [-0.35cm]
            & $0$ & $\frac{1}{6}$ & $\frac{1}{6}$ & $\frac{2}{3}$
		\end{tabular}
	\end{center}
    \end{minipage}
    \begin{minipage}{0.5\textwidth}
	\begin{center}
		\begin{tabular}{>{\centering\arraybackslash}p{7mm}|>{\centering\arraybackslash}p{7mm}cc>{\centering\arraybackslash}p{7mm}}
			$\alpha$ & $\alpha$ & $0$ & $0$ & $0$ \\[2mm]
            $0$ & $-\alpha$ & $\alpha$ & $0$ & $0$ \\[2mm]
			$1$ & $0$ & $1 - \alpha$ & $\alpha$ & $0$ \\[2mm]
            $\frac{1}{2}$ & $\frac{\alpha}{4}$ & $\eta$ & $\frac{1}{2} - \frac{5}{4}\alpha - \eta$ & $\alpha$ \\[2mm]
			\hline \\[-0.35cm]
			& $0$ & $\frac{1}{6}$ & $\frac{1}{6}$ & $\frac{2}{3}$
		\end{tabular}
	\end{center}
    \end{minipage}
    \caption{\textbf{Butcher tableaux of IMEX-SSP(4,3,3) scheme employed in~\cite{boscarino:2016}.} Method with  $\alpha = 0.24169426$ and $\eta = 0.12915286$. Explicit method (left) and implicit method (right).}
    \label{tab:rk3_butch_type_I_SSP}
\end{table}

\begin{table}[H]
    \small
	\begin{tabular}{c|ccccc}
		$0$ & $0$ & $0$ & $0$ & $0$ & $0$ \\[3pt]
        $0.87173304$ & $\;\;\:0.87173304$ & $0$ & $0$ & $0$ & $0$ \\[3pt]
        $0.87173304$ & $\;\;\:0.43586652$ & $\;\;\:0.43586652$ & $0$ & $0$ & $0$ \\[3pt]
        $0.65905793$ & $-0.80099845$ & $0$ & $3.14121102$ & $0$ & $0$ \\[3pt]
        $1$ & $\;\;\:0.35675320$ & $-0.19733989$ & $0.88194884$ & $-0.04136215$ & $0$ \\[3pt]
		\hline\\[-6pt]
		& $\;\;\:0.41289804$ & $0$ & $0.19733989$ & $-0.04610445$ & $0.43586652$ \\[3pt]
	\end{tabular}
    \vskip 0.3cm
	\begin{tabular}{c|ccccccc}
        $0$ & $0$ & $0$ & $0$ & $0$ & $0$ \\[3pt]
        $0.87173304$ & $\;\;\:0.43586652$ & $0.43586652$ & $0$ & $0$ & $0$ \\[3pt]
        $0.87173304$ & $\;\;\:0.43586652$ & $0$ & $0.43586652$ & $0$ & $0$ \\[3pt]
        $0.65905793$ & $-0.06675868$ & $0$ & $1.97110474$ & $\;\;\:0.43586652$ & $0$ \\[3pt]
        $1$ & $\;\;\:0.41289804$ & $0$ & $0.19733989$ & $-0.04610445$ & $0.43586652$ \\[3pt]
		\hline\\[-6pt]
		& $\;\;\:0.41289804$ & $0$ & $0.19733989$ & $-0.04610445$ & $0.43586652$
	\end{tabular}
    \caption{\textbf{Butcher tableaux of BHR(5,5,3) scheme employed in~\cite{boscarino:2009_third_order}.} Explicit method (top) and implicit method (bottom).}
    \label{tab:rk3_butch_type_II_BHR}
\end{table}

\begin{table}[H]
\small
	\begin{tabular}{c|ccccc}
		$0$ & $0$ & $0$ & $0$ & $0$ & $0$ \\[3pt]
        $1.05144292$ & %2870010$ & 
        $1.05144292$ & %2870010$ & 
        $0$ & $0$ & $0$ & $0$ \\[3pt]
        $0.66658810$ & %2470224$ & 
        $0.45528814$ & %0917395$ & 
        $0.21129996$ & %1552829$ & 
        $0$ & $0$ & $0$ \\[3pt]
        $1.47683871$ & %8193715$ & 
        $0.06936473$ & %4498722$ & 
        $0$ & 
        $1.40747398$ & %3694993$ & 
        $0$ & $0$ \\[3pt]
        $1$ & 
        $0.31801663$ & %8728572$ & 
        $0.02902514$ & %3044335$ & 
        $0.61071506$ & %7260961$ & 
        $0.04224315$ & %0966131$ & 
        $0$ \\[3pt]
        \hline \\[-6pt]
		& 
        $0.34225438$ & %8132523$ & 
        $0.04450272$ & %5380489$ & 
        $0.24901994$ & %8096975$ & 
        $-0.16149852$ & %3044992$ & 
        $0.52572146$ %1435005$
	\end{tabular}
    \vskip 0.3cm
	\begin{tabular}{c|ccccc}
		$0$ & $0$ & $0$ & $0$ & $0$ & $0$ \\[3pt]
        $1.05144292$ & %2870010$ & 
        $0.52572146$ & %1435005$ & 
        $0.52572146$ & %1435005$ & 
        $0$ & $0$ & $0$ \\[3pt]
		$0.66658810$ & %2470224$ & 
        $0.26286073$ & %0717503$ & 
        $-0.12199408$ & %9682283$ & 
        $0.52572146$ & %1435005$ & 
        $0$ &$0$ \\[3pt]
        $1.47683871$ & %8193715$ & 
        $0.07706201$ & %0265066$ & 
        $-0.69770293$ & %4810540$ & 
        $1.57175818$ & %1304184$ & 
        $0.52572146$ & %1435005$ & 
        $0$ \\[3pt]
        $1$ & 
        $0.34225438$ & %8132523$ & 
        $0.04450272$ & %5380489$ & 
        $0.24901994$ & %8096975$ & 
        $-0.16149852$ & %3044992$ & 
        $0.52572146$ %1435005$
        \\[3pt]
		\hline \\[-6pt]
		& 
        $0.34225438$ & %8132523$ & 
        $0.04450272$ & %5380489$ & 
        $0.24901994$ & %8096975$ & 
        $-0.16149852$ & %3044992$ & 
        $0.52572146$ %1435005$
	\end{tabular}
    \caption{\textbf{Butcher tableaux of \orderthreeLDsone\ scheme.} Explicit method (top) and implicit method (bottom).}
    \label{tab:ACO_order3_LDs_1}
\end{table}

\begin{table}[H]
\small
    \begin{tabular}{c|ccccc}
		$0$ & $0$ & $0$ & $0$ & $0$ & $0$ \\[3pt]
        $1.05144292$ & %2870010$ & 
        $1.05144292$ & %2870010$ & 
        $0$ & $0$ & $0$ & $0$ \\[3pt]
        $0.66658810$ & %2470224$ & 
        $0.45528814$ & %0917395$ & 
        $0.21129996$ & %1552829$ & 
        $0$ & $0$ & $0$ \\[3pt]
        $1.47683871$ & %8193715$ & 
        $0.06936473$ & %4498722$ & 
        $0$ & 
        $1.40747398$ & %3694993$ & 
        $0$ & $0$ \\[3pt]
        $1$ & 
        $0.33449382$ & %2441987$ & 
        $0.03977721$ & %5228850$ & 
        $0.57503718$ & %1169806$ & 
        $\;\;\:0.05069178$ & %1159357$ & 
        $0$ \\[3pt]
        \hline \\[-6pt]
        & 
        $0.34225438$ & %8132523$ & 
        $0.04450272$ & %5380489$ & 
        $0.24901994$ & %8096975$ & 
        $-0.16149852$ & %3044992$ & 
        $0.52572146$ %1435005$
	\end{tabular}
    \vskip 0.3cm
    \begin{tabular}{c|ccccc}
		$0$ & $0$ & $0$ & $0$ & $0$ & $0$ \\[3pt]
        $1.05144292$ & %2870010$ & 
        $0.52572146$ & %1435005$ & 
        $\;\;\:0.52572146$ & %1435005$ & 
        $0$ & $0$ & $0$ \\[3pt]
		$0.66658810$ & %2470224$ & 
        $0.26286073$ & %0717503$ & 
        $-0.12199408$ & %9682283$ & 
        $0.52572146$ & %1435005$ & 
        $0$ &$0$ \\[3pt]
        $1.47683871$ & %8193715$ & 
        $0.07706201$ & %0265066$ & 
        $-0.69770293$ & %4810540$ & 
        $1.57175818$ & %1304184$ & 
        $\;\;\:0.52572146$ & %1435005$ & 
        $0$ \\[3pt]
        $1$ & 
        $0.34225438$ & %8132523$ & 
        $\;\;\:0.04450272$ & %5380489$ & 
        $0.24901994$ & %8096975$ & 
        $-0.16149852$ & %3044992$ & 
        $0.52572146$ %1435005$
        \\[3pt]
		\hline \\[-6pt]
		&         
        $0.34225438$ & %8132523$ & 
        $\;\;\:0.04450272$ & %5380489$ & 
        $0.24901994$ & %8096975$ & 
        $-0.16149852$ & %3044992$ & 
        $0.52572146$ %1435005$
	\end{tabular}
    \caption{\textbf{Butcher tableaux of \orderthreeLDp\ scheme.} Explicit method (top) and implicit method (bottom).}
    \label{tab:ACO_order3_LDp}
\end{table}

\begin{table}[H]
\small
    \begin{tabular}{c|ccccc}
	    $0$ & $0$ & $0$ & $0$ & $0$ & $0$ \\[3pt]
        $1.05144292$ & %2870010$ & 
        $\;\;\:1.05144292$ & %2870010$ & 
        $0$ & $0$ & $0$ & $0$ \\[3pt]
        $0.60111815$ & %2198527$ & 
        $\;\;\:0.42928617$ & %4901912$ & 
        $\;\;\:0.17183197$ & %7296614$ & 
        $0$ & $0$ & $0$ \\[3pt]
        $1.26219537$ & %2957851$ & 
        $-0.06294940$ & %4062670$ & 
        $0$ & 
        $1.32514477$ & %7020521$ & 
        $0$ & $0$ \\[3pt]
        $1$ & 
        $\;\;\:0.22291660$ & %3586566$ & 
        $-0.00797809$ & %1493323$ & 
        $0.72988821$ & %0533753$ & 
        $\;\;\:0.05517327$ & %7373004$ & 
        $0$ \\[3pt]
        \hline \\[-6pt]
        & 
        $\;\;\:0.28387977$ & %0698957$ & 
        $0$ & 
        $0.40243695$ & %7284376$ & 
        $-0.21203818$ & %9418338$ & 
        $0.52572146$ %1435005$
	\end{tabular}
    \vskip 0.3cm
    \begin{tabular}{c|ccccc}
		$0$ & $0$ & $0$ & $0$ & $0$ & $0$ \\[3pt]
        $1.05144292$ & %2870010$ & 
        $0.52572146$ & %1435005$ & 
        $\;\;\:0.52572146$ & %1435005$ & 
        $0$ & $0$ & $0$ \\[3pt]
		$0.60111815$ & %2198527$ & 
        $0.20412379$ & %0.204123789566170$ & 
        $-0.12872709$ & %8802649$ & 
        $0.52572146$ & %1435005$ & 
        $0$ &$0$ \\[3pt]
        $1.26219537$ & %2957851$ & 
        $0.12407344$ & %4260833$ & 
        $-0.52211130$ & %4837317$ & 
        $1.13451177$ & %2099330$ & 
        $\;\;\:0.52572146$ & %1435005$ & 
        $0$ \\[3pt]
        $1$ & 
        $0.28387977$ & %0698957$ 
        $0$ & 
        $0.40243695$ & %7284376$ & 
        $-0.21203818$ & %9418338$ & 
        $0.52572146$ %1435005$ 
        \\[3pt]
		\hline \\[-6pt]
		& 
        $0.28387977$ & %0698957$ & 
        $0$ & 
        $0.40243695$ & %7284376$ & 
        $-0.21203818$ & %9418338$ & 
        $0.52572146$ %1435005$
	\end{tabular}
    \caption{\textbf{Butcher tableaux of \orderthreeLDstwo\ scheme.} Explicit method (top) and implicit method (bottom).}
    \label{tab:ACO_order3_LDs_2}
\end{table}

\color{black}
\subsection{Fourth-order schemes}

Regarding the \emph{fourth-order} schemes, we list here the schemes considered.

\begin{table}[H]
    \begin{minipage}{0.5\textwidth}
	\begin{center}
		\begin{tabular}{>{\centering\arraybackslash}p{7mm}|>{\centering\arraybackslash}p{7mm}>{\centering\arraybackslash}p{7mm}>{\centering\arraybackslash}p{7mm}>{\centering\arraybackslash}p{7mm}>{\centering\arraybackslash}p{7mm}>{\centering\arraybackslash}p{7mm}}
			$0$ & $0$ & $0$ & $0$ & $0$ & $0$ & $0$ \\[2mm]
            $\frac{1}{4}$ & $\frac{1}{4}$ & $0$ & $0$ & $0$ & $0$ & $0$ \\[2mm]
            $\frac{3}{4}$ & $-\frac{1}{4}$ & $1$ & $0$ & $0$ & $0$ & $0$ \\[2mm]
            $\frac{11}{20}$ & $-\frac{13}{100}$ & $\frac{43}{75}$ & $\frac{8}{75}$ & $0$ & $0$ & $0$ \\[2mm]
            $\frac{1}{2}$ & $-\frac{6}{85}$ & $\frac{42}{85}$ & $\frac{179}{1360}$ & $-\frac{15}{272}$ & $0$ & $0$ \\[2mm]
            $1$ & $0$ & $\frac{79}{24}$ & $-\frac{5}{8}$ & $\frac{25}{2}$ & $-\frac{85}{6}$ & $0$ \\[2mm]
            \hline \\[-0.35cm]
            & $0$ & $\frac{25}{24}$ & $-\frac{49}{48}$ & $\frac{125}{16}$ & $-\frac{85}{12}$ & $\frac{1}{4}$
		\end{tabular}
	\end{center}
    \end{minipage}
    \begin{minipage}{0.5\textwidth}
	\begin{center}
		\begin{tabular}{>{\centering\arraybackslash}p{7mm}|>{\centering\arraybackslash}p{7mm}>{\centering\arraybackslash}p{7mm}>{\centering\arraybackslash}p{7mm}>{\centering\arraybackslash}p{7mm}>{\centering\arraybackslash}p{7mm}>{\centering\arraybackslash}p{7mm}}
			$0$ & $0$ & $0$ & $0$ & $0$ & $0$ & $0$ \\[2mm]
            $\frac{1}{4}$ & $0$ & $\frac{1}{4}$ & $0$ & $0$ & $0$ & $0$ \\[2mm]
            $\frac{3}{4}$ & $0$ & $\frac{1}{2}$ & $\frac{1}{4}$ & $0$ & $0$ & $0$ \\[2mm]
            $\frac{11}{20}$ & $0$ & $\frac{17}{50}$ & $-\frac{1}{25}$ & $\frac{1}{4}$ & $0$ & $0$ \\[2mm]
            $\frac{1}{2}$ & $0$ & $\frac{371}{1360}$ & $-\frac{137}{2720}$ & $\frac{15}{544}$ & $\frac{1}{4}$ & $0$ \\[2mm]
            $1$ & $0$ & $\frac{25}{24}$ & $-\frac{49}{48}$ & $\frac{125}{16}$ & $-\frac{85}{12}$ & $\frac{1}{4}$ \\[2mm]
            \hline \\[-0.35cm]
            & $0$ & $\frac{25}{24}$ & $-\frac{49}{48}$ & $\frac{125}{16}$ & $-\frac{85}{12}$ & $\frac{1}{4}$
		\end{tabular}
	\end{center}
    \end{minipage}
    \caption{\textbf{Butcher tableaux of IMEX-ARS(5,5,4) scheme employed in~\cite{boscarino:2016}.} Explicit method (left) and implicit method (right).}
    \label{tab:rk4_butch_type_ARS}
\end{table}

\begin{table}[H]
    \small
	\begin{tabular}{c|cccccc}
		$0$ & $0$ & $0$ & $0$ & $0$ & $0$ & $0$ \\[3pt]
        $0.50000000$ & 
        $\;\;\:0.50000000$ & $0$ & $0$ & $0$ & $0$ & $0$ \\[3pt]
        $0.33200000$ & 
        $\;\;\:0.22177600$ & 
        $\;\;\:0.11022400$ & $0$ & $0$ & $0$ & $0$ \\[3pt]
        $0.62000000$ & 
        $-0.04884659$ & 
        $-0.17772065$ &
        $0.84656724$ & $0$ & $0$ & $0$ \\[3pt]
        $0.85000000$ & 
        $-0.15541685$ & 
        $-0.35670501$ & 
        $1.05872588$ & 
        $0.30339599$ & $0$ & $0$ \\[3pt]
        $1$ & 
        $\;\;\:0.20142435$ & 
        $\;\;\:0.00874206$ & 
        $0.15993996$ & 
        $0.40382906$ & 
        $\;\;\:0.22606457$ & $0$ \\[3pt]
		\hline \\[-6pt]
        & 
        $\;\;\:0.15791629$ & 
        $0$ & 
        $0.18675894$ & 
        $0.68056529$ & 
        $-0.27524053$ & 
        $0.25000000$
	\end{tabular}
	\vskip 0.3cm
	\begin{tabular}{c|cccccc}
		$0$ & $0$ & $0$ & $0$ & $0$ & $0$ & $0$ \\[3pt]
        $0.50000000$ & 
        $0.25000000$ & 
        $\;\;\:0.25000000$ & $0$ & $0$ & $0$ & $0$ \\[3pt]
        $0.33200000$ & 
        $0.13777600$ & 
        $-0.05577600$ & 
        $0.25000000$ & $0$ & $0$ & $0$ \\[3pt]
        $0.62000000$ & 
        $0.14463687$ & 
        $-0.22393190$ & 
        $0.44929504$ & 
        $0.25000000$ & $0$ & $0$ \\[3pt]
        $0.85000000$ & 
        $0.09825878$ & 
        $-0.59154424$ & 
        $0.81012105$ & 
        $0.28316440$ & 
        $\;\;\:0.25000000$ & $0$ \\[3pt]
        $1$ & 
        $0.15791629$ & $0$ & 
        $0.18675894$ & 
        $0.68056529$ & 
        $-0.27524053$ & 
        $0.25000000$ \\[3pt]
        \hline \\[-6pt]
        &        
        $\;\;\:0.15791629$ & 
        $0$ & 
        $0.18675894$ & 
        $0.68056529$ & 
        $-0.27524053$ & 
        $0.25000000$
    \end{tabular}    
	\caption{\textbf{Butcher tableaux of ARK4(3)6L[2]SA scheme employed in~\cite{kennedy:2003}.} Explicit method (top) and implicit method (bottom).}
	\label{tab:rk4_butch_type_II_2003}
\end{table}

\begin{table}[H]
    \small
    \begin{tabular}{c|ccccccc}
        $0$ & $0$ & $0$ & $0$ & $0$ & $0$ & $0$ & $0$ \\[3pt]
        $2\gamma$ & $2\gamma$ & $0$ & $0$ & $0$ & $0$ & $0$ & $0$ \\[3pt]
        $0.42165537$ & $0.06175000$ & $\;\;\:0.35990537$ & $0$ & $0$ & $0$ & $0$ & $0$ \\[3pt]
        $0.33500000$ & $0.05301658$ & $\;\;\:0.35949264$ & $-0.07750923$ & $0$ & $0$ & $0$ & $0$ \\[3pt]
        $0.07500000$ & $0.05841715$ & $-0.16313824$ & $-0.19732090$ & $\;\;\:0.37704199$ & $0$ & $0$ & $0$ \\[3pt]
        $0.70000000$ & $0.53853032$ & $-0.45497746$ & $\;\;\:1.25629056$ & $-0.47828452$ & $-0.16155888$ & $0$ & $0$ \\[3pt]
        $1$ & $0.23221715$ & $\;\;\:0.23221715$ & $-6.80999437$ & $\;\;\:7.36185855$ & $-1.37487908$ & $1.35858058$ & $0$ \\[3pt]
        \hline\\[-6pt]
        & $0$ &	$0$ & $\;\;\:0.51611072$ & $-0.14606356$ & $\;\;\:0.23473048$ & $0.27172234$ & $\gamma$
    \end{tabular}
	\vskip 0.3cm
    \begin{tabular}{c|ccccccc}
        $0$ & $0$ & $0$ & $0$ & $0$ & $0$ & $0$ & $0$ \\[3pt]
        $2\gamma$ & $\gamma$ & $\gamma$ & $0$ & $0$ & $0$ & $0$ & $0$ \\[3pt]
        $0.42165537$ & $\;\;\:0.14907768$ & $\;\;\:0.14907768$ & $\gamma$ & $0$ & $0$ & $0$ & $0$ \\[3pt]
        $0.33500000$ & $\;\;\:0.12483442$ & $\;\;\:0.12483442$ & $-0.03816885$ & 
        $\gamma$ & $0$ & $0$ & $0$ \\[3pt]
        $0.07500000$ & $-0.07303194$ & $-0.07303194$ & $-0.24343569$ & $\;\;\:0.34099956$ & $\gamma$ & $0$ & $0$ \\[3pt]
        $0.70000000$ & $-0.15296500$ &	$-0.15296500$ & $\;\;\:0.07220562$ & $\;\;\:0.40430630$ & $0.40591807$ & $\gamma$ & $0$ \\[3pt]
        $1$ & $0$ &	$0$ & $\;\;\:0.51611072$ & $-0.14606356$ & $0.23473048$ & $0.27172234$ & $\gamma$ \\[3pt]
        \hline\\[-6pt]
        & $0$ &	$0$ & $\;\;\:0.51611072$ & $-0.14606356$ & $0.23473048$ & $0.27172234$ & $\gamma$
    \end{tabular}
	\caption{\textbf{Butcher tableaux of ARK4(3)7L[2]SA$_{1}$ scheme employed in~\cite{kennedy:2019}.} Method with $\gamma = 0.1235$. Explicit method (top) and implicit method (bottom).}
	\label{tab:rk4_butch_type_II_2019}
\end{table}

\begin{table}[H]
    \small
	\begin{tabular}{c|ccccccc}
        $0$ & $0$ & $0$ & $0$ & $0$ & $0$ & $0$ \\[3pt]
        $0.55610768$ & %2272905$ & 
        $\;\;\:0.55610768$ & %227290464986317237023839$ 
        $0$ & $0$ & $0$ & $0$ & $0$ \\[3pt]
        $1.31576490$ & %3497378$ & 
        $-0.33180581$ & %032565871254392355403328$ & 
        $\;\;\:1.64757071$ & %382303625045407770459630$ & 
        $0$ & $0$ & $0$ & $0$ \\[3pt]
        $0.08673938$ & %4684563$ & 
        $\;\;\:0.13811029$ & %769046660532597640572413$ & 
        $-0.00270761$ & %127732173009126628066634$ & 
        $-0.04866330$ & %301581983009531538519497$ & 
        $0$ & $0$ & $0$ \\[3pt]
        $0.43558989$ & %2310426$ & 
        $\;\;\:1.03496518$ & %884090638146231919917500$ & 
        $\;\;\:0.54237034$ & %776306547993572581119410$ & 
        $-0.07939527$ & %807662554836611104295002$ & 
        $-1.06235036$ & %621692072916679280153180$ & 
        $0$ & $0$ \\[3pt]
        $1$ & 
        $-1.25469717$ & %090818147632558505554870$ & 
        $-0.59939455$ & %496817102968818520335416$ & 
        $\;\;\:0.17680815$ & %355709059002670038356215$ & 
        $\;\;\:1.62106030$ & %222482609603456971823010$ & 
        $\;\;\:1.05622327$ & %009443581995250015711060$ & 
        $0$ \\[3pt]
        \hline\\[-6pt]
        & $\;\;\:0.09317511$ & %019158532269929086458649$ & 
        $0$ & 
        $-0.03410445$ & %733053458230064096824696$ &
        $0.06284142$ & %612304550131324344021923$ & 
        $0.60003407$ & %987945143335652047832205$ & 
        $0.27805384$ %113645232493158618511920$
	\end{tabular}
	\vskip 0.3cm
	\begin{tabular}{c|cccccc}
        $0$ & $0$ & $0$ & $0$ & $0$ & $0$ & $0$ \\[3pt]
        $0.55610768$ & %2272905$ & 
        $\;\;\:0.27805384$ & %113645232493158618511920$ & 
        $\;\;\:0.27805384$ & %113645232493158618511920$ & 
        $0$ & $0$ & $0$ & $0$ \\[3pt]
        $1.31576490$ & %3497378$ & 
        $\;\;\:0.13902692$ & %056822616246579309255960$ & 
        $\;\;\:0.89868414$ & %179269905051277487288426$ & 
        $\;\;\:0.27805384$ & %113645232493158618511920$ & 
        $0$ & $0$ & $0$ \\[3pt]
        $0.08673938$ & %4684563$ & 
        $-0.12975100$ & %972489478661258367819389$ & 
        $-0.07983424$ & %916612824842428625843588$ & 
        $\;\;\:0.01827080$ & %243913366342683406624498$ & 
        $0.27805384$ & %113645232493158618511920$ & 
        $0$ & $0$ \\[3pt] 
        $0.43558989$ & %2310426$ & 
        $-0.14018717$ & %247383101994907530182104$ & 
        $-0.12525399$ & %327070800495223459400541$ & 
        $\;\;\:0.00546592$ & %294760848330877165563101$ & 
        $0.41751129$ & %397090380052609322096349$ & 
        $0.27805384$ & %113645232493158618511920$ & 
        $0$ \\[3pt]
        $1$ & 
        $\;\;\:0.09317511$ & %019158532269929086458649$ & 
        $0$ & 
        $-0.03410445$ & %733053458230064096824696$ &
        $0.06284142$ & %612304550131324344021923$ & 
        $0.60003407$ & %987945143335652047832205$ & 
        $0.27805384$ %113645232493158618511920$ 
        \\[3pt]
        \hline\\[-6pt]
        & 
        $\;\;\:0.09317511$ & %019158532269929086458649$ & 
        $0$ & 
        $-0.03410445$ & %733053458230064096824696$ &
        $0.06284142$ & %612304550131324344021923$ & 
        $0.60003407$ & %987945143335652047832205$ & 
        $0.27805384$ %113645232493158618511920$ 
        \\[3pt]
    \end{tabular}
	\caption{\textbf{Butcher tableaux of \orderfourLDpone\ scheme.} Explicit method (top) and implicit method (bottom).}
	\label{tab:ACO_order4_LDp_1}
\end{table}

\begin{table}[H]
    \small
	\begin{tabular}{c|cccccc}
        $0$ & $0$ & $0$ & $0$ & $0$ & $0$ & $0$ 
        \\[3pt]
        $0.55610768$ & %2272905$ & 
        $0.55610768$ & %227290464986317237023839$ & 
        $0$ & $0$ & $0$ & $0$ & $0$ 
        \\[3pt]
        $0.35255814$ & %3321336$ & 
        $0.04247560$ & %080447506060855927639016$ & 
        $\;\;\:0.31008254$ & %251686135107080368376196$ & 
        $0$ & $0$ & $0$ & $0$ 
        \\[3pt]
        $0.66303019$ & %6713514$ & 
        $0.40040468$ & %265288505091002364829786$ & 
        $-0.03002200$ & %675489878595472852530716$ & 
        $0.29264752$ & %081552805793238726192826$ & 
        $0$ & $0$ & $0$ 
        \\[3pt]
        $0.89332539$ & %9490450$ & 
        $0.52440469$ & %239104730072560485858219$ & 
        $-0.72698360$ & %570313783936100410627791$ & 
        $0.23689320$ & %144933913998869762808758$ & 
        $0.85901111$ & %135320130721548832080942$ & 
        $0$ & $0$ 
        \\[3pt]
        $1$ & 
        $0.23510769$ & %305610948771312618193231$ & 
        $-1.18401116$ & %044701764413359872134580$ & 
        $0.30313013$ & %299242794353175910958896$ & 
        $1.81785990$ & %774140857445327787034050$ & 
        $-0.17208657$ & %334292836156456444051600$ & 
        $0$ 
        \\[3pt]
        \hline\\[-6pt]
        & 
        $0.14889680$ & %421719961224557532327388$ & 
        $0$ & 
        $0.29346325$ & %451908912498901770627694$ & 
        $0.57004255$ & %532867482896979984417250$ & 
        $-0.29045645$ & %520141589113597905884251$ & 
        $0.27805384$ %113645232493158618511920$
        \\[3pt]
    \end{tabular}
	\vskip 0.3cm
	\begin{tabular}{c|cccccc}
        $0$ & $0$ & $0$ & $0$ & $0$ & $0$ & $0$ \\[3pt]
        $0.55610768$ & %2272905$ & 
        $\;\;\:0.27805384$ & %113645232493158618511920$ & 
        $\;\;\:0.27805384$ & %113645232493158618511920$ & 
        $0$ & $0$ & $0$ & $0$ \\[3pt]
        $0.35255814$ & %3321336$ & 
        $\;\;\:0.13902692$ & %056822616246579309255960$ & 
        $-0.06452261$ & %838334207571801631752668$ & 
        $0.27805384$ & %113645232493158618511920$ & 
        $0$ & $0$ & $0$ 
        \\[3pt]
        $0.66303019$ & %6713514$ & 
        $\;\;\:0.09279114$ & %760305742474026378876776$ & 
        $-0.33193809$ & %667149103216397513690069$ & 
        $0.62412330$ & %464549560537980754793269$ & 
        $0.27805384$ & %113645232493158618511920$ & 
        $0$ & $0$ 
        \\[3pt]
        $0.89332539$ & %9490450$ & 
        $-0.06208062$ & %888338897526444414917835$ & 
        $-0.96963964$ & %122181663070068393977856$ & 
        $1.29531364$ & %295374075045234286681320$ & 
        $0.35167818$ & %550546243914998573822581$ & 
        $\;\;\:0.27805384$ & %113645232493158618511920$ & 
        $0$ \\[3pt]
        $1$ & 
        $\;\;\:0.14889680$ & %421719961224557532327388$ & 
        $0$ & 
        $0.29346325$ & %451908912498901770627694$ & 
        $0.57004255$ & %532867482896979984417250$ & 
        $-0.29045645$ & %520141589113597905884251$ & 
        $0.27805384$ %113645232493158618511920$ 
        \\[3pt]
        \hline\\[-6pt]
        & 
        $0.14889680$ & %421719961224557532327388$ & 
        $0$ & 
        $0.29346325$ & %451908912498901770627694$ & 
        $0.57004255$ & %532867482896979984417250$ & 
        $-0.29045645$ & %520141589113597905884251$ & 
        $0.27805384$ %113645232493158618511920$
        \\[3pt]
    \end{tabular}
	\caption{\textbf{Butcher tableaux of \orderfourLDptwo\ scheme.} Explicit method (top) and implicit method (bottom).}
	\label{tab:ACO_order4_LDp_2}
\end{table}

\begin{table}[H]
    \small
	\begin{tabular}{c|cccccc}
        $0$ & $0$ & $0$ & $0$ & $0$ & $0$ & $0$ \\[3pt]
        $0.55610768$ & %2272905$ & 
        $\;\;\:0.55610768$ & %227290464986317237023839$ & 
        $0$ & $0$ & $0$ & $0$ & $0$ 
        \\[3pt]
        $0.35255814$ & %3321336$ & 
        $\;\;\:0.31554253$ & %093184639303084259824532$ & 
        $\;\;\:0.03701561$ & %238949001864852036190680$ & 
        $0$ & $0$ & $0$ & $0$ 
        \\[3pt]
        $0.66303019$ & %6713514$ & 
        $\;\;\:0.04611287$ & %117258900095246292712830$ & 
        $-0.49427752$ & %986851948700005800583809$ & 
        $\;\;\:1.11119485$ & %540944480893527746362880$ & 
        $0$ & $0$ & $0$ \\[3pt]
        $0.03422808$ & %4897008$ & 
        $\;\;\:0.10180581$ & %495807785690414967551824$ & 
        $\;\;\:0.05442689$ & %501481922858747545101440$ & 
        $-0.14018574$ & %586005685769611326629230$ & 
        $0.01818112$ & %078416786981773742624129$ & 
        $0$ & $0$ \\[3pt]
        $1$ & 
        $-2.05008320$ & %478711729352722665447620$ & 
        $-1.19976062$ & %005794354477278735970680$ & 
        $\;\;\:0.13177740$ & %572165126417939861672168$ & 
        $1.55817443$ & %165718761681359260677570$ & 
        $2.55989198$ & %746622195730702279068560$ & 
        $0$ \\[3pt]
        \hline\\[-6pt]
        &         
        $4.03283234$ & %478929511218199434039990$ & 
        $0$ & 
        $1.76674922$ & %958315221718413968940680$ & 
        $-0.36122323$ & %523760461364042905812805$ & 
        $-4.71641218$ & %027129504065729115679790$ & 
        $0.27805384$ %113645232493158618511920$
        \\[3pt]
    \end{tabular}
	\vskip 0.3cm
	\begin{tabular}{c|cccccc}
        $0$ & $0$ & $0$ & $0$ & $0$ & $0$ & $0$ \\[3pt]
        $0.55610768$ & %2272905$ & 
        $\;\;\:0.27805384$ & %113645232493158618511920$ & 
        $\;\;\:0.27805384$ & %113645232493158618511920$ &
        $0$ & $0$ & $0$ & $0$ \\[3pt]
        $0.35255814$ & %3321336$ & 
        $\;\;\:0.13902692$ & %056822616246579309255960$ & 
        $-0.06452261$ & %838334207571801631752668$ & 
        $\;\;\:0.27805384$ & %113645232493158618511920$ & 
        $0$ & $0$ & $0$ \\[3pt]
        $0.66303019$ & %6713514$ & 
        $\;\;\:0.0927911$ & %4760305742474026378876776$ & 
        $-0.33193809$ & %667149103216397513690069$ & 
        $\;\;\:0.62412330$ & %464549560537980754793269$ & 
        $\;\;\:0.27805384$ & %113645232493158618511920$ & 
        $0$ & $0$ \\[3pt]
        $0.03422808$ & %4897008$ & 
        $-0.20769213$ & %693465886655372655076388$ & 
        $-0.01432779$ & %094615112420605594178011$ & 
        $-0.04346364$ & %689299883053218719099570$ & 
        $\;\;\:0.02165781$ & %853436459397363278490212$ & 
        $\;\;\:0.27805384$ & %113645232493158618511920$ & 
        $0$ \\[3pt]
        $1$ & 
        $\;\;\:4.03283234$ & %478929511218199434039990$ & 
        $0$ & 
        $\;\;\:1.76674922$ & %958315221718413968940680$ & 
        $-0.36122323$ & %523760461364042905812805$ & 
        $-4.71641218$ & %027129504065729115679790$ & 
        $0.27805384$ %113645232493158618511920$ 
        \\[3pt]
        \hline\\[-6pt]
        & 
        $4.03283234$ & %478929511218199434039990$ & 
        $0$ & 
        $1.76674922$ & %958315221718413968940680$ & 
        $-0.36122323$ & %523760461364042905812805$ & 
        $-4.71641218$ & %027129504065729115679790$ & 
        $0.27805384$ %113645232493158618511920$ 
    \end{tabular}
	\caption{\textbf{Butcher tableaux of \orderfourLDpthree\ scheme.} Explicit method (top) and implicit method (bottom).}
	\label{tab:ACO_order4_LDp_3}
\end{table}

%%%%%%%%%%%% Supplementary index-1 DAE order conditions %%%%%%%%
\section{Derivation of index-1 DAE supplementary conditions up to order 4}
\label{app:index1_DAE_order4}

We consider the singularly perturbed system
\begin{equation}
\begin{aligned}
    y'(t) &= f(y(t),z(t)), \\
    \varepsilon z'(t) &= g(y(t),z(t)),
\end{aligned}
\end{equation}
with $0 < \varepsilon \ll 1$ and $g_{z}(y,z)$ invertible on the region of interest (\emph{index-1} DAE case). We focus on schemes of type I and we assume $\mathbf{A}$ invertible, so that $\mathbf{w} = \mathbf{A}^{-1}$. The derived relations will be valid for schemes of type II by considering $\boldsymbol{\omega} = \boldsymbol{\mathcal{A}}^{-1}$~\cite{boscarino:2009_third_order}, that is assumed to be invertible. An IMEX-RK scheme of type I applied to the reduced DAE reads~\cite{boscarino:2007}
\begin{subequations}
\begin{align}
    y^{(\timesuperscript,l)} &= y^{\timesuperscript} + \Delta t \sum_{m=1}^{s}\hat{a}_{lm}f(y^{(\timesuperscript,m)},z^{(\timesuperscript,m)}), \\
    0 &= g(y^{(\timesuperscript,m)},z^{(\timesuperscript,m)}), \\
    y^{\timesuperscript+1} &= y^{\timesuperscript} + \Delta t \sum_{l=1}^{s}\hat{b}_{l}f(y^{(\timesuperscript,l)},z^{(\timesuperscript,l}) \\
    z^{\timesuperscript+1} &= \rho z^{\timesuperscript} + \Delta t \sum_{l,m=1}^{s}b_{l}w_{lm}z^{(\timesuperscript,m)} \label{eq:znp1_index_1_order0},
\end{align}
\end{subequations}
where $\rho = R(\infty) < 1$. We need to compute the local truncation error $\Delta z^{\timesuperscript+1}$. We drop the explicit dependence on time, all the quantities are evaluated at time $t^{\timesuperscript}$ (subscript $n$). First, we recall the Taylor expansion of the exact solution
\begin{equation}\label{eq:z_ex_Taylor_simp}
    z_{n+1} = z_{n} + \Delta t z'_{n} + \frac{\Delta t^{2}}{2} z_{n}'' + \frac{\Delta t^{3}}{6}z_{n}''' + \mathcal{O}(\Delta t^{4}).
\end{equation}
We introduce for the sake of convenience
\begin{subequations}
\begin{align}
    \mathcal{Z}_{1,n} &= -g_{z,n}^{-1}g_{y,n} \\
    \mathcal{Z}_{2,n} &= -g_{z,n}^{-1}\rpth{g_{yy,n} + 2g_{yz,n}\mathcal{Z}_{1,n} + g_{zz,n}\mathcal{Z}_{1,n}^{2}} \\
    \mathcal{Z}_{3,n} &= -g_{z,n}^{-1}\rpth{g_{yyy,n} + 3g_{yyz,n}\mathcal{Z}_{1,n} + 3g_{yzz,n}\mathcal{Z}_{1,n}^{2} + g_{zzz,n}\mathcal{Z}_{1,n}^{3} + 3g_{yz,n}\mathcal{Z}_{2,n} + 3g_{zz,n}\mathcal{Z}_{1,n}\mathcal{Z}_{2,n}},
\end{align} 
\end{subequations}
so that
{
\small
$$
\begin{aligned}
    z_{n}' &= \mathcal{Z}_{1,n}y'_{n} = \mathcal{Z}_{1,n}f_{n}, \\
    z_{n}'' &= \mathcal{Z}_{1,n}y''_{n} + \mathcal{Z}_{2,n}(y'_{n})^{2} = \mathcal{Z}_{1,n}y''_{n} + \mathcal{Z}_{2,n}f_{n}^{2} \\
    &= -g_{z,n}^{-1}\Big(\underbrace{g_{y,n}\rpth{f_{y,n}f_{n} + f_{z,n}z_{n}'}}_{\text{first-order derivatives block}} \\
    &\hphantom{= -g_{z,n}^{-1}\Big(}
    + \underbrace{g_{yy,n}f_{n}^{2} + 2g_{yz,n}f_{n}z_{n}' + g_{zz,n}(z_{n}')^{2}}_{\text{second-order derivatives block}}\Big) \\
\end{aligned}
$$
$$
\begin{aligned}
    z_{n}''' &= \mathcal{Z}_{1,n}y_{n}''' + 3\mathcal{Z}_{2,n}y'_{n}y''_{n} + \mathcal{Z}_{3}(y'_{n})^{3} = \mathcal{Z}_{1,n}y_{n}''' + 3\mathcal{Z}_{2,n}f_{n}y''_{n} + \mathcal{Z}_{3,n}f_{n}^{3} \\
    &= -g_{z,n}^{-1}\Big(\underbrace{g_{y,n}\rpth{f_{y,n}(f_{y,n}f_{n} + f_{z,n}z_{n}') + f_{z,n}z_{n}'' + f_{yy,n}f_{n}^{2} + 2f_{yz,n}f_{n}z_{n}' + f_{zz,n}z_{n}'^{2}}}_{\text{first-order derivatives block}} \\
    &\hphantom{= -g_{z,n}^{-1}\Big(}
    + \underbrace{3\rpth{g_{yy,n}f_{n}\rpth{f_{y,n}f_{n} + f_{z,n}z_{n}'} + g_{yz,n}\rpth{z_{n}'\rpth{f_{y,n}f_{n} + f_{z,n}z_{n}'} + f_{n}z_{n}''} + g_{zz,n}z_{n}'z_{n}''}}_{\text{second-order derivatives block}} \\
    &\hphantom{= -g_{z,n}^{-1}\Big(}
    + \underbrace{g_{yyy,n}f_{n}^{3} + 3g_{yyz,n}f_{n}^{2}z_{n}' + 3g_{yzz,n}f_{n}(z_{n}')^{2} + g_{zzz,n}(z_{n}')^{3}}_{\text{third-order derivatives block}}\Big) \\ 
    &= -g_{z,n}^{-1}\Big(\underbrace{g_{y,n}\rpth{f_{y,n}(f_{y,n}f_{n} + f_{z,n}z_{n}') + f_{z,n}z_{n}'' + f_{yy,n}f_{n}^{2} + 2f_{yz,n}f_{n}z_{n}' + f_{zz,n}z_{n}'^{2}}}_{\text{first-order derivatives block}} \\
    &\hphantom{= -g_{z,n}^{-1}\Big(}
    + \underbrace{3\rpth{g_{yy,n}f_{n}\rpth{f_{y,n}f_{n} + f_{z,n}z_{n}'} + g_{yz,n}\rpth{2z_{n}'\rpth{f_{y,n}f_{n} + f_{z,n}z_{n}'} + \mathcal{Z}_{2,n}f_{n}^{3}} + g_{zz,n}z_{n}'z_{n}''}}_{\text{second-order derivatives block}} \\
    &\hphantom{= -g_{z,n}^{-1}\Big(}
    + \underbrace{g_{yyy,n}f_{n}^{3} + 3g_{yyz,n}f_{n}^{2}z_{n}' + 3g_{yzz,n}f_{n}(z_{n}')^{2} + g_{zzz,n}(z_{n}')^{3}}_{\text{third-order derivatives block}}\Big).
\end{aligned}
$$
}
Notice in particular that
$$
\begin{aligned}
    g_{yz,n}\rpth{z_{n}'\rpth{f_{y,n}f_{n} + f_{z,n}z_{n}'} + f_{n}z_{n}''} &= g_{yz,n}\rpth{z_{n}'y_{n}'' + f_{n}z_{n}''} \\
    &= g_{yz,n}\rpth{\mathcal{Z}_{1,n}f_{n}y''_{n} + f_{n}\rpth{\mathcal{Z}_{1,n}y_{n}'' + \mathcal{Z}_{2,n}f_{n}^{3}}} \\
    &= g_{yz,n}\rpth{2\mathcal{Z}_{1,n}f_{n}y''_{n} + \mathcal{Z}_{2,n}f_{n}^{3}}
\end{aligned}
$$
and that
$$z_{n}'z_{n}'' = z_{n}'\rpth{\mathcal{Z}_{1,n}y''_{n} + \mathcal{Z}_{2,n}f_{n}^{2}} = \mathcal{Z}_{1,n}^{2}f_{n}y''_{n} + \mathcal{Z}_{1,n}\mathcal{Z}_{2,n}f_{n}^{3}.$$

In order to compute the Taylor expansion of the numerical method at time $t^{n}$ we need to compute the Taylor expansion of the intermediate stages $z^{(n,m)}$. These are defined implicitly, through the relation $g(y^{(n,m)}, z^{(n,m)}) = 0$, so that
\begin{equation}
\begin{aligned}
    0 &= g(y^{(n,m)}, z^{(n,m)}) = g(y_{n} + \Delta y^{(n,m)}, z_{n} + \Delta z^{(n,m)}) \\
    &= \underbrace{g_n}_{= g(y_{n},z_{n}) = g(y(t^{n}),z(t^{n})) = 0} + g_{y,n}\Delta y^{(n,m)} + g_{z,n}\Delta z^{(n,m)} \\
    &+ \frac{1}{2}g_{yy,n}\rpth{\Delta y^{(n,m)}}^{2} + g_{yz,n}\Delta y^{(n,m)}\Delta z^{(n,m)} + \frac{1}{2}g_{zz,n}\rpth{\Delta z^{(n,m)}}^{2} \\
    &+ \frac{1}{6}g_{yyy,n}\rpth{\Delta y^{(n,m)}}^{3} + \frac{1}{2}g_{yyz,n}\rpth{\Delta y^{(n,m)}}^{2}\Delta z^{(n,m)} + \frac{1}{2} g_{yzz,n}\Delta y^{(n,m)}\rpth{\Delta z^{(n,m)}}^{2} \\
    &+ \frac{1}{6} g_{zzz,n}\rpth{\Delta z^{(n,m)}}^{3} + \mathcal{O}\rpth{\Delta y^{(n,m)}}^{4} + \mathcal{O}\rpth{\Delta z^{(n,m)}}^{4},
\end{aligned}
\end{equation}
so that, by defining
$$\Delta z^{(n,m)} = \Delta z^{(n,m),1} + \Delta z^{(n,m),2} + \Delta z^{(n,m),3} + \mathcal{O}(\Delta t^{4})$$
one obtains
{
\small
\begin{subequations}
\begin{align}
    \Delta z^{(n,m),1} &= -g_{z,n}^{-1}g_{y,n}  \Delta y^{(n,m),1}, \\[1mm]
    \Delta z^{(n,m),2} &= -g_{z,n}^{-1}\rpth{g_{y,n}\Delta y^{(n,m),2} + \frac{1}{2}g_{yy,n}\rpth{\Delta y^{(n,m),1}}^{2} + g_{yz,n}\Delta y^{(n,m),1}\Delta z^{(n,m),1} + \frac{1}{2}g_{zz,n}\rpth{\Delta z^{(n,m),1}}^{2}}, \\[1mm]
    \Delta z^{(n,m),3} &= -g_{z,n}^{-1}\Bigg(g_{y,n}\Delta y^{(n,m),3} + g_{yy,n}\Delta y^{(n,m),1}\Delta y^{(n,m),2} + g_{zz,n}\Delta z^{(n,m),1} \Delta z^{(n,m),2} \nonumber \\
    &\hphantom{= -g_{z,n}^{-1}\Bigg(} + g_{yz,n}\Delta y^{(n,m),1}  \Delta z^{(n,m),2} + g_{yz,n}\Delta y^{(n,m),2}\Delta z^{(n,m),1} \nonumber \\
    &\hphantom{= -g_{z,n}^{-1}\Bigg(} + \frac{1}{6}g_{yyy,n}(\Delta y^{(n,m),1})^{3} + \frac{1}{6}g_{zzz,n}(\Delta z^{(n,m),1})^{3} \nonumber \\
    &\hphantom{= -g_{z,n}^{-1}\Bigg(} + \frac{1}{2}g_{yyz,n}(\Delta y^{(n,m),1})^{2}\Delta z^{(n,m),1} + \frac{1}{2}g_{yzz,n}\Delta y^{(n,m),1}(\Delta z^{(n,m),1})^{2}\Bigg).
\end{align}
\end{subequations}
}

Now we need to compute the expansions of $\Delta y^{(n,m)}$, so that
$$\Delta y^{(n,m)} = \Delta y^{(n,m),1} + \Delta y^{(n,m), 2} + \Delta y^{(n,m), 3} + \mathcal{O}(\Delta t^{4}).$$
In order to achieve this goal, we need to consider the expansion of $f$
\begin{equation}
\begin{aligned}
    f(y^{(n,m)}, z^{(n,m)}) &= f_{n} + f_{y,n}\Delta y^{(n,m)} + f_{z,n}\Delta z^{(n,m)} \\
    &+ \frac{1}{2}f_{yy,n}(\Delta y^{(n,m)})^{2} + f_{yz,n}\Delta y^{(n,m)}\Delta z^{(n,m)} + \frac{1}{2}f_{zz,n}(\Delta z^{(n,m)})^{2} \\
    &+ \frac{1}{6}f_{yyy,n}(\Delta y^{(n,m)})^{3} + \frac{1}{2}f_{yyz,n} (\Delta y^{(n,m)})^{2}\Delta z^{(n,m)} \\
    &+ \frac{1}{2}f_{yzz,n}\Delta y^{(n,m)}(\Delta z^{(n,m)})^{2} + \frac{1}{6}f_{zzz,n}(\Delta z^{(n,m)})^{3}
\end{aligned}
\end{equation}
More specifically, we obtain
\paragraph{Order 1 ($\Delta t$)}
\begin{equation}
    \Delta y^{(n,m),1} = \Delta t\sum_{k}\hat{a}_{mk}f_{n} = \Delta t\hat{c}_{m}f_{n},
\end{equation}    
\paragraph{Order 2 ($\Delta t^{2}$)}
\begin{equation}
    \Delta y^{(n,m),2} = \Delta t\sum_{k}\hat{a}_{mk} \rpth{f_{y,n}\Delta y^{(n,k),1} + f_{z,n}\Delta z^{(n,k),1}},
\end{equation}

\paragraph{Order 3 ($\Delta t^{3}$)}
{
\footnotesize
\begin{equation}
\begin{aligned}
    \Delta y^{(n,m),3} &= \Delta t \sum_{k} \hat{a}_{mk}\Bigg(f_{y,n}\Delta y^{(n,k),2} + f_{z,n}\Delta z^{(n,k),2} \\
    &\hphantom{= \Delta t \sum_{k} \hat{a}_{mk}\Bigg(} + \frac{1}{2}f_{yy,n}(\Delta y^{(n,k),1})^{2} + f_{yz,n}\Delta y^{(n,k),1}\Delta z^{(n,k),1} + \frac{1}{2}f_{zz,n}(\Delta z^{(n,k),1})^{2}\Bigg).
\end{aligned}
\end{equation}
}

By successive substitution, we obtain
\paragraph{Order 1 ($\Delta t$)}
\begin{equation}
    \Delta z^{(n,m),1} = \mathcal{Z}_{1,n}\Delta y^{(n,m),1} = -\Delta t g_{z,n}^{-1}g_{y,n}\hat{c}_{m}f_{n} = \Delta t\hat{c}_{m}z_{n}'.
\end{equation}

\paragraph{Order 2 ($\Delta t^{2}$)}
{
\footnotesize
\begin{subequations}
\begin{align}
  \Delta y^{(n,m),2} &= \Delta t\sum_{k}\hat{a}_{mk}\rpth{f_{y,n}\Delta y^{(n,k),1} + f_{z,n}\mathcal{Z}_{1,n}\Delta y^{(n,k),1}} = \Delta t^{2}\sum_k \hat{a}_{mk}\hat{c}_{k}(f_{y,n}f_{n} + f_{z,n}z_{n}') \\
  \Delta z^{(n,m),2} &= \mathcal{Z}_{1,n}\Delta y^{(n,m),2} + \frac{1}{2}\mathcal{Z}_{2,n}(\Delta y^{(n,m),1})^{2} \\
  &= -g_{z,n}^{-1}\spth{g_{y,n}\Delta y^{(n,m),2} + \frac{1}{2}\rpth{g_{yy,n} + 2g_{yz,n}\mathcal{Z}_{1,n} + g_{zz,n}\mathcal{Z}_{1,n}^{2}}(\Delta y^{(n,m),1})^{2}} \nonumber \\
  &= -\Delta t^{2}g_{z,n}^{-1} \rpth{\underbrace{\sum_{k}\hat{a}_{mk}\hat{c}_{k}g_{y,n}(f_{y,n}f_{n} + f_{z,n}z_{n}')}_{\text{first-order derivatives block}} + \underbrace{\frac{\hat{c}_{m}^{2}}{2}\rpth{g_{yy,n}f_{n}^{2} + 2g_{yz,n}f_{n}z_{n}' + g_{zz,n}(z_{n}')^{2}}}_{\text{second-order derivatives block}}}. \nonumber
\end{align}
\end{subequations}
}

\paragraph{Order 3 ($\Delta t^{3}$)}
{
\scriptsize
\begin{subequations}
\begin{align}
    \Delta y^{(n,m),3} &= \Delta t \sum_{k}\hat{a}_{mk}\Bigg(f_{y,n}\Delta y^{(n,k),2} + f_{z,n}\rpth{\mathcal{Z}_{1,n}\Delta y^{(n,k),2} + \frac{1}{2}\mathcal{Z}_{2,n}(\Delta y^{(n,k),1})^{2}} \nonumber \\
    &\hphantom{= \Delta t \sum_{k} \hat{a}_{mk}\Bigg(} + \frac{1}{2}(\Delta y^{(n,k),1})^{2}\rpth{f_{yy,n} + 2f_{yz,n}\mathcal{Z}_{1,n} + f_{zz,n}\mathcal{Z}_{1,n}^{2}}\Bigg) \\
    &= \Delta t \sum_{k} \hat{a}_{mk}\Bigg(f_{y,n}\Delta y^{(n,k),2} + f_{z,n}\rpth{\mathcal{Z}_{1,n}\Delta y^{(n,k),2} + \frac{1}{2}\mathcal{Z}_{2,n}(\Delta y^{(n,k),1})^{2}} \nonumber \\
    &\hphantom{= \Delta t \sum_{k} \hat{a}_{mk}\Bigg(} + \Delta t^{2}\frac{\hat{c}_{k}^{2}}{2}\rpth{f_{yy,n}f_{n}^{2} + 2f_{yz,n}f_{n}z_{n}' + f_{zz,n}z_{n}'^{2}}\Bigg) \nonumber \\
    &= \Delta t^{3}\sum_{k}\hat{a}_{mk}\Big[ 
    \underbrace{\sum_{j}\hat{a}_{kj}\hat{c}_{j}\rpth{f_{y,n} + f_{z,n}\mathcal{Z}_{1,n}}\rpth{f_{y,n}f_{n} + f_{z,n}z_{n}'} + \frac{\hat{c}_{k}^{2}}{2} f_{z,n}\mathcal{Z}_{2,n}f_{n}^{2}}_{\text{first-order derivatives block}} \nonumber \\
    &\hphantom{= \Delta t^{3}\sum_{k}\hat{a}_{mk}\Big[} + \underbrace{\frac{\hat{c}_{k}^{2}}{2}\rpth{f_{yy,n}f_{n}^{2} + 2 f_{yz,n}f_{n}z'_{n} + f_{zz,n}z_{n}'^{2}}}_{\text{second-order derivatives block}}\Big]. \nonumber
\end{align}
\end{subequations}
}

Before deriving $\Delta z^{(n,m),3}$, we notice that
{
\footnotesize
$$
\begin{aligned}
    &\frac{1}{6}g_{yyy,n}(\Delta y^{(n,m),1})^{3} + \frac{1}{6}g_{zzz,n}(\Delta z^{(n,m),1})^{3} + \frac{1}{2}g_{yyz,n}(\Delta y^{(n,m),1})^{2}\Delta z^{(n,m),1} + \frac{1}{2}g_{yzz,n}\Delta y^{(n,m),1}(\Delta z^{(n,m),1})^{2} = \\
    &\frac{1}{6}(\Delta y^{(n,m),1})^{3}\rpth{g_{yyy,n} + g_{zzz,n}\mathcal{Z}_{1,n}^{3} + 3g_{yyz,n}\mathcal{Z}_{1,n} + 3g_{yyz,n}\mathcal{Z}_{1,n}^{2}}.
\end{aligned}
$$
}
Moreover,
$$
\begin{aligned}
    \Delta z^{(n,m),1}\Delta z^{(n,m),2} &= \mathcal{Z}_{1,n}\Delta y^{(n,m),1}\rpth{\mathcal{Z}_{1,n}\Delta y^{(n,m),2} + \frac{1}{2}\mathcal{Z}_{2,n}(\Delta y^{(n,m),1})^{2}} \\
    \Delta y^{(n,m),1}\Delta z^{(n,m),2} &= \Delta y^{(n,m),1}\rpth{\mathcal{Z}_{1,n}\Delta y^{(n,m),2} + \frac{1}{2}\mathcal{Z}_{2,n}(\Delta y^{(n,m),1})^{2}} \\
    \Delta y^{(n,m),2}\Delta z^{(n,m),1} &= \mathcal{Z}_{1,n}\Delta y^{(n,m),1}\Delta y^{(n,m),2},
\end{aligned}
$$
so that
{
\footnotesize
$$
\begin{aligned}
    &g_{yy,n}\Delta y^{(n,m),1}\Delta y^{(n,m),2} + g_{zz,n}\Delta z^{(n,m),1}\Delta z^{(n,m),2} + g_{yz,n}\Delta y^{(n,m),1}\Delta z^{(n,m),2} + g_{yz,n}\Delta y^{(n,m),2}\Delta z^{(n,m),1} = \\
    &\Delta y^{(n,m),1}\Delta y^{(n,m),2}\rpth{g_{yy,n} + g_{zz,n}\mathcal{Z}_{1,n}^{2} + 2g_{yz,n}\mathcal{Z}_{1,n}} + \frac{1}{2}\rpth{g_{yz,n} + g_{zz,n}\mathcal{Z}_{1,n}}\mathcal{Z}_{2,n}(\Delta y^{(n,m),1})^{3}.
\end{aligned}
$$
}
Hence
\begin{equation}\label{eq:delta_z_3}
    \Delta z^{(n,m),3} = \mathcal{Z}_{1,n}\Delta y^{(n,m),3} + \mathcal{Z}_{2,n}\Delta y^{(n,m),1}\Delta y^{(n,m),2} + \frac{1}{6}\mathcal{Z}_{3,n}(\Delta y^{(n,m),1})^{3}.
\end{equation}
We progressively expand $\Delta z^{(n,m),3}$ in~\eqref{eq:delta_z_3}. We obtain
{
\scriptsize
$$
\begin{aligned}
    \Delta z^{(n,m),3} &= \mathcal{Z}_{1,n}\Delta y^{(n,m),3} + \mathcal{Z}_{2,n}\Delta y^{(n,m),1}\Delta y^{(n,m),2} \\
    &-g_{z,n}^{-1}\spth{\Delta t^{3}\frac{\hat{c}_{m}^{3}}{6}f_{n}^{3}\rpth{g_{yyy,n} + g_{zzz,n}\mathcal{Z}_{1,n}^{3} + 3g_{yyz,n}\mathcal{Z}_{1,n} + 3g_{yyz,n}\mathcal{Z}_{1,n}^{2} + 3\rpth{g_{yz,n} + g_{zz,n}\mathcal{Z}_{1,n}}\mathcal{Z}_{2,n}}} \\
    &= \mathcal{Z}_{1,n}\Delta y^{(n,m),3} + \mathcal{Z}_{2,n}\Delta t^{2}\hat{c}_{m}f_{n}\sum_{k}\hat{a}_{mk} \rpth{f_{y,n}\Delta y^{(n,k),1} + f_{z,n}\Delta z^{(n,k),1}} \\
    &-g_{z,n}^{-1}\spth{\Delta t^{3}\frac{\hat{c}_{m}^{3}}{6}f_{n}^{3}\rpth{g_{yyy,n} + g_{zzz,n}\mathcal{Z}_{1,n}^{3} + 3g_{yyz,n}\mathcal{Z}_{1,n} + 3g_{yyz,n}\mathcal{Z}_{1,n}^{2} + 3\rpth{g_{yz,n} + g_{zz,n}\mathcal{Z}_{1,n}}\mathcal{Z}_{2,n}}} \\
    &= \mathcal{Z}_{1,n}\Delta y^{(n,m),3} + \mathcal{Z}_{2,n}\Delta t^{3}\hat{c}_{m}f_{n}\sum_{k}\hat{a}_{mk}\hat{c}_{k}\rpth{f_{y,n}f_{n} + f_{z,n}z_{n}'} \\
    &-g_{z,n}^{-1}\spth{\Delta t^{3}\frac{\hat{c}_{m}^{3}}{6}f_{n}^{3}\rpth{g_{yyy,n} + g_{zzz,n}\mathcal{Z}_{1,n}^{3} + 3g_{yyz,n}\mathcal{Z}_{1,n} + 3g_{yyz,n}\mathcal{Z}_{1,n}^{2} + 3\rpth{g_{yz,n} + g_{zz,n}\mathcal{Z}_{1,n}}\mathcal{Z}_{2,n}}} \\ 
    &= -\Delta t^{3}g_{z,n}^{-1}\Big[\sum_{k}\hat{a}_{mk}\sum_{j}\hat{a}_{kj}\hat{c}_{j}g_{y,n}\rpth{f_{y,n} + f_{z,n}\mathcal{Z}_{1,n}}\rpth{f_{y,n}f_{n} + f_{z,n}z_{n}'} + \sum_{k}\hat{a}_{mk}\frac{\hat{c}_{k}^{2}}{2}g_{y,n}f_{z,n}\mathcal{Z}_{2,n}f_{n}^{2} \\
    &\hphantom{-\Delta t^{3}g_{z,n}^{-1}\Big[+} + \sum_{k}\hat{a}_{mk}\frac{\hat{c}_{k}^{2}}{2}g_{y,n}\rpth{f_{yy,n}f_{n}^{2} + 2f_{yz,n}f_{n}z_{n}' + f_{zz,n}z_{n}'^{2}} \\
    &\hphantom{-\Delta t^{3}g_{z,n}^{-1}\Big[+} + \rpth{g_{yy,n} + 2g_{yz,n}\mathcal{Z}_{1,n} + g_{zz,n}\mathcal{Z}_{1,n}^{2}}\hat{c}_{m}f_{n}\sum_{k}\hat{a}_{mk}\hat{c}_{k}\rpth{f_{y,n}f_{n} + f_{z,n}z_{n}'} \\
    &\hphantom{-\Delta t^{3}g_{z,n}^{-1}\Big[+} + \frac{\hat{c}_{m}^{3}}{6}f_{n}^{3}\rpth{g_{yyy,n} + g_{zzz,n}\mathcal{Z}_{1,n}^{3} + 3g_{yyz,n}\mathcal{Z}_{1,n} + 3g_{yyz,n}\mathcal{Z}_{1,n}^{2} + 3\rpth{g_{yz,n} + g_{zz,n}\mathcal{Z}_{1,n}}\mathcal{Z}_{2,n}}\Big] \\
    &= -\Delta t^{3}g_{z,n}^{-1}\Big[\underbrace{\sum_{k}\hat{a}_{mk}\sum_{j}\hat{a}_{kj}\hat{c}_{j}g_{y,n}\rpth{f_{y,n}(f_{y,n}f_{n} + f_{z,n}z_{n}') + f_{z,n}\mathcal{Z}_{1,n}y_{n}''} + \sum_{k}\hat{a}_{mk}\frac{\hat{c}_{k}^{2}}{2}g_{y,n}f_{z,n}\mathcal{Z}_{2,n}f_{n}^{2}}_{\text{first-order derivatives block}} \\
    &\hphantom{-\Delta t^{3}g_{z,n}^{-1}\Big[+} \underbrace{+ \sum_{k}\hat{a}_{mk}\frac{\hat{c}_{k}^{2}}{2}g_{y,n}\rpth{f_{yy,n}f_{n}^{2} + 2f_{yz,n}f_{n}z_{n}' + f_{zz,n}z_{n}'^{2}}}_{\text{first-order derivatives block}} \\
    &\hphantom{-\Delta t^{3}g_{z,n}^{-1}\Big[+} \underbrace{+ \hat{c}_{m}\sum_{k}\hat{a}_{mk}\hat{c}_{k}\rpth{g_{yy,n}f_{n}\rpth{f_{y,n}f_{n} + f_{z,n}z_{n}'} + 2g_{yz,n}z_{n}'\rpth{f_{y,n}f_{n} + f_{z,n}z_{n}'} + g_{zz,n}\mathcal{Z}_{1,n}^{2}f_{n}\rpth{f_{y,n}f_{n} + f_{z,n}z_{n}'}}}_{\text{second-order derivatives block}} \\
    &\hphantom{-\Delta t^{3}g_{z,n}^{-1}\Big[+} + \frac{\hat{c}_{m}^{3}}{2}g_{yz,n}\mathcal{Z}_{2,n}f_{n}^{3} + \frac{\hat{c}_{m}^{3}}{2}g_{zz,n}\mathcal{Z}_{1,n}\mathcal{Z}_{2,n}f_{n}^{3} \\
    &\hphantom{-\Delta t^{3}g_{z,n}^{-1}\Big[+} + \underbrace{\frac{\hat{c}_{m}^{3}}{6}\rpth{g_{yyy,n}f_{n}^{3} + 3g_{yyz,n}f_{n}^{2}z_{n}' + 3g_{yzz,n}f_{n}(z_{n}')^{2} + g_{zzz,n}(z_{n}')^{3}}}_{\text{third-order derivatives block}}\Big].
\end{aligned}
$$
}
Recall indeed that $z_{n}'' = \mathcal{Z}_{1,n}y_{n}'' + \mathcal{Z}_{2,n}f_{n}^{2}$ and that $z_{n}'z_{n}'' = \mathcal{Z}_{1,n}^{2}f_{n}y_{n}'' + \mathcal{Z}_{1,n}\mathcal{Z}_{2,n}f_{n}^{3}$. Comparing the numerical and exact expansions, we obtain the supplementary order conditions specific to order 4 reported in~\eqref{eq:index_1_DAE_conditions_order_4}. An alternative, more elegant way to derive such conditions is by means of bicolour rooted trees~\cite{boscarino:2009_third_order, hairer:1996}. As discussed in~\cite{boscarino:2009_third_order}, in order to distinguish derivatives with respect to $y$ (the differential variable) and $z$ (the algebraic variable), we need two kinds of vertices: meagre ($\circ$) and fat ($\bullet$). We associate the expression $(-g_{z,n}^{-1}g)$ with a fat vertex and $f$ with a meagre one. Hence, the trees associated to the supplementary fourth-order conditions are

$$\sum_{l,m} b_{l}w_{lm}\hat{c}_{m}^{3} = 1, \qquad
\begin{tikzpicture}[baseline=-0.5ex]
    \node[circle, draw, fill=white, inner sep=2pt] (root) {};
    \node[circle, draw, fill=black, inner sep=2pt, above left=10pt of root] (l) {};
    \node[circle, draw, fill=black, inner sep=2pt, above=7pt of root] (c) {};
    \node[circle, draw, fill=black, inner sep=2pt, above right=10pt of root] (r) {};
    \draw (root) -- (l);
    \draw (root) -- (c);
    \draw (root) -- (r);
\end{tikzpicture}
$$

$$\sum_{l,m,k} b_{l}w_{lm}\hat{c}_{m}\hat{a}_{mk}\hat{c}_{k} = \frac{1}{2}, \qquad
\begin{tikzpicture}[baseline=-0.5ex]
    \node[circle, draw, fill=white, inner sep=2pt] (root) {};
    \node[circle, draw, fill=black, inner sep=2pt, above left=10pt of root] (l) {};
    \node[circle, draw, fill=black, inner sep=2pt, above right=10pt of root] (r) {};
    \node[circle, draw, fill=black, inner sep=2pt, above=10pt of r] (rr) {};
    \draw (root) -- (l);
    \draw (root) -- (r);
    \draw (r) -- (rr);
\end{tikzpicture}
$$

$$\sum_{l,m,k} b_{l}w_{lm}\hat{a}_{mk}\hat{c}_{k}^{2} = \frac{1}{3}, \qquad
\begin{tikzpicture}[baseline=-0.5ex]
    \node[circle, draw, fill=white, inner sep=2pt] (root) {};
    \node[circle, draw, fill=black, inner sep=2pt, above=7pt of root] (c) {};
    \node[circle, draw, fill=black, inner sep=2pt, above left=10pt of c] (cl) {};
    \node[circle, draw, fill=black, inner sep=2pt, above right=10pt of c] (cr) {};
    \draw (root) -- (c);
    \draw (c) -- (cl);
    \draw (c) -- (cr);
\end{tikzpicture}
$$

$$\sum_{l,m,k,j} b_{l}w_{lm}\hat{a}_{mk}\hat{a}_{kj}\hat{c}_{j} = \frac{1}{6}, \qquad
\begin{tikzpicture}[baseline=-0.5ex]
    \node[circle, draw, fill=white, inner sep=2pt] (root) {};
    \node[circle, draw, fill=black, inner sep=2pt, above=7pt of root] (c) {};
    \node[circle, draw, fill=black, inner sep=2pt, above=7pt of c] (cc) {};
    \node[circle, draw, fill=black, inner sep=2pt, above=7pt of cc] (ccc) {};
    \draw (root) -- (c);
    \draw (c) -- (cc);
    \draw (cc) -- (ccc);
\end{tikzpicture}
$$

The coefficients appearing in the previous relations are related to the so-called \emph{density} of a tree $\tau$~\cite{butcher:2021, ketcheson:2023}. We recall that the density $\tilde{\gamma}(\tau)$ is defined recursively by
\begin{equation}
    \tilde{\gamma}(\bullet) = 1, \qquad \tilde{\gamma}(\tau) = |\tau|\prod_{i=1}^{m}\tilde{\gamma}(\tau_{i}),
\end{equation}
where $|\tau|$ is the \emph{order} of the tree, i.e. the number of nodes. Here the tree $\tau$ is expressed in terms of its children, i.e. $\tau= [\tau_{1}, \dots, \tau_{m}]$, where $\tau_{1}, \dots, \tau_{m}$ are the subtrees obtained by removing the root of $\tau$. As an example, consider the tree
\begin{tikzpicture}[baseline=-0.5ex]
    \node[circle, draw, fill=white, inner sep=2pt] (root) {};
    \node[circle, draw, fill=black, inner sep=2pt, above left=10pt of root] (l) {};
    \node[circle, draw, fill=black, inner sep=2pt, above=7pt of root] (c) {};
    \node[circle, draw, fill=black, inner sep=2pt, above right=10pt of root] (r) {};
    \draw (root) -- (l);
    \draw (root) -- (c);
    \draw (root) -- (r);
\end{tikzpicture}
which has density $\tilde{\gamma}(\tau) = 4$. The coefficients appearing in the order conditions required to formally match the Taylor expansion of the exact solution are given by $1/\tilde{\gamma}(\tau)$~\cite{ketcheson:2023}. In the present context, however, one node plays a special role, namely the meagre node. Consequently, the coefficients reported above have to be computed by considering the tree obtained after removing the meagre root.

%%%%%%%%%%% Bibliography %%%%%%%%%%%%
\bibliographystyle{plain}
\bibliography{bibliography}

\end{document}